\documentclass[10pt]{article}
\date{23.5.12} %
\usepackage{hyperref}
\usepackage{amssymb, amsmath}
\input{liemacs10.sty} 

\newcommand\uline{\underline}
 
\newcommand\cR{\mathcal{R}}

\newcommand\AU{\mathop{\rm AU}\nolimits}
\newcommand\PAU{\mathop{\rm PAU}\nolimits}

\newtheorem{theori}{\theoremname}
\newenvironment{thmi}{\begin{theori}\it}{\end{theori}}

\renewcommand{\mlabel}{\label}

\begin{document} 

\title{Semibounded Unitary Representations of\\ Double Extensions of 
Hilbert--Loop Groups}
\author{K.-H. Neeb\begin{footnote}{
Department  Mathematik, FAU Erlangen-N\"urnberg, Cauerstrasse 11, 
91058 Erlangen, Germany, Email: karl-hermann.neeb@math.uni-erlangen.de}
\end{footnote}
\begin{footnote}{Supported by DFG-grant NE 413/7-1, Schwerpunktprogramm 
``Darstellungstheorie''.} 
\end{footnote}}

\maketitle

\begin{abstract}  A unitary representation 
of a, possibly infinite dimensional, Lie group 
$G$ is called semibounded if the corresponding 
operators $i\dd\pi(x)$ from the derived representation  
are uniformly bounded from above on some non-empty open subset 
of the Lie algebra $\g$ of $G$. We classify all 
irreducible semibounded representations of the 
groups $\hat\cL_\phi(K)$ which are double extensions 
of the twisted loop group $\cL_\phi(K)$, where 
$K$ is a simple Hilbert--Lie group 
(in the sense that the scalar product on its Lie algebra 
is invariant) and $\phi$ is a finite order automorphism of $K$ 
which leads to one of the $7$ irreducible locally affine root systems 
with their canonical $\Z$-grading. To achieve this goal, we extend the method  
of holomorphic induction to certain classes of Fr\'echet--Lie groups and 
prove an infinitesimal characterization of analytic operator-valued 
positive definite functions on Fr\'echet--BCH--Lie groups. 
\end{abstract}

\section*{Introduction} \mlabel{sec:0}

This paper is part of a project concerned with a 
systematic approach to unitary representations 
of infinite-dimensional Lie groups in terms of semiboundedness conditions 
on spectra in the derived representation. For the derived 
representation to carry significant information, we have to 
impose a suitable smoothness condition: Let $G$ be a Lie group with Lie algebra $\g$ 
and exponential function $\exp \: \g \to G$. A unitary 
representation $\pi \: G \to \U(\cH)$ 
is said to be {\it smooth} if the subspace 
$\cH^\infty \subeq \cH$ of smooth vectors is dense. 
This is automatic for continuous 
representations of finite-dimensional Lie groups, but not 
for Banach--Lie groups (\cite{Ne10a}). For any smooth 
unitary representation, the {\it derived representation} 
\[ \dd\pi \: \g \to \End(\cH^\infty), \quad 
\dd\pi(x)v := \derat0 \pi(\exp tx)v\]  
carries significant information in the sense that the closure of the 
operator $\dd\pi(x)$ coincides with the infinitesimal generator of the 
unitary one-parameter group $\pi(\exp tx)$. We call $(\pi, \cH)$ 
{\it semibounded} if the function 
\[ s_\pi \: \g \to \R \cup \{ \infty\}, \quad 
s_\pi(x) 
:= \sup\big(\Spec(i\dd\pi(x))\big) \]
is bounded on a neighborhood of some point $x_0 \in \g$.  
Then the set $W_\pi$ of all such points $x_0$ 
is an open $\Ad(G)$-invariant convex cone. 
We say that  $\pi$ is {\it bounded} if $W_\pi = \g$. 
All finite-dimensional continuous unitary representations are bounded 
and most of the unitary representations appearing in physics are semibounded 
or satisfy similar spectral conditions (cf.\ \cite{Se58, Se78}, \cite{SeG81}, 
\cite{Ca83}, \cite{PS86}, \cite{CR87}, \cite{Mi89}, 
\cite{Ot95},  \cite{FH05}, \cite{Bak07}, \cite{Ne10b}).  

For finite-dimensional Lie groups, the irreducible semibounded 
representations are precisely the unitary highest weight representations 
and one has unique direct integral decompositions into irreducible ones 
\cite[X.3/4, XI.6]{Ne00}. 
For many other classes of groups such as 
the Virasoro group and affine Kac--Moody groups 
(double extensions of loop groups with compact target groups),  
the irreducible highest weight representations 
are semibounded by Theorem~\ref{thm:6.1} below, but to prove the converse 
is more difficult and requires a thorough understanding 
of invariant cones in the corresponding Lie algebras as well 
as of convexity properties of coadjoint orbits (\cite[Sect.~8]{Ne10b}). 

The closest infinite-dimensional relatives of compact Lie algebras are 
{\it Hilbert--Lie algebras}. These are real Lie algebras which are Hilbert spaces 
on which the adjoint group acts by isometries (cf. \cite[Def.~6.3]{HoMo98}). 
We call a Lie group $K$ whose Lie algebra $\fk = \L(K)$ is a Hilbert--Lie algebra 
a {\it Hilbert--Lie group}.\begin{footnote}{In the literature one also finds a 
weaker concept of a ``Hilbert--Lie group,'' namely Lie groups whose Lie algebra 
is a Hilbert space, but no compatibility between the Lie bracket and the 
scalar product is required. }  
\end{footnote}
The finite-dimensional Hilbert--Lie algebras are the compact Lie algebras. 
The main goal of this paper is the classification of all 
semibounded unitary representations of groups which are double extensions of loop groups with 
values in a Hilbert--Lie group. 

For a Hilbert--Lie group $K$, we write $\Aut(K)$ for the group of all 
Lie group automorphisms acting by isometries with respect to the 
scalar product $\la \cdot,\cdot \ra$ on~$\fk$. 
For an automorphism 
$\phi \in \Aut(K)$ of order $N$, 
\[  \cL_\phi(K) 
:= \Big\{ f \in C^\infty(\R,K) \: (\forall t \in \R)\ 
f\Big(t + \frac{2\pi}{N}\Big) = \phi^{-1}(f(t))\Big\} \] 
is called the corresponding {\it twisted loop group}. 
It is a Fr\'echet--Lie group with Lie algebra 
\[ \cL_\phi(\fk) 
:= \Big\{ \xi \in C^\infty(\R,\fk) \: (\forall t \in \R)\ 
\xi\Big(t + \frac{2\pi}{N}\Big) = \L(\phi)^{-1}(\xi(t))\Big\} \] 
(\cite[App. A]{NeWo09}). 
The subgroup $\cL_\phi(K) \subeq C^\infty(\R,K)$ is 
translation invariant, so that we obtain for each $T \in \R$ an 
automorphism of $\cL_\phi(K)$ by 
\begin{equation}
  \label{eq:translat}
\alpha_T(f)(t) := f(t + T) \quad \mbox{ with } \quad 
\L(\alpha_T)(\xi)(t) := \xi(t + T). 
\end{equation}
Our assumption $\phi^N = \id$ implies that 
$\alpha_{2\pi} = \id$, which leads to a smooth action 
of the circle group $\T \cong \R/2\pi \Z$ on $\cL_\phi(K)$. 
The Fr\'echet--Lie algebra $\cL_\phi(\fk)$ carries the positive definite form 
$$ \la \xi, \eta \ra := \frac{1}{2\pi}\int_0^{2\pi} \la \xi(t), \eta(t)\ra \, dt,  $$
which is invariant under the $\R$-action \eqref{eq:translat}. Therefore the derivation 
$D\xi := \xi'$ is skew-symmetric and thus 
$$ \omega(\xi, \eta) := \la D\xi,\eta\ra = \frac{1}{2\pi} 
\int_0^{2\pi} \la \xi'(t), \eta(t)\ra\, dt $$
defines a continuous Lie algebra cocycle on $\cL_\phi(\fk)$ 
(cf.\ \cite[Sect.~3.2]{NeWo09}). 
Let 
$$\tilde\cL_\phi(\fk) := \R \oplus_\omega \cL_\phi(\fk) $$
denote the corresponding central extension and observe that 
$\omega$ is $D$-invariant, so that we obtain the double extension 
$$\g := \hat\cL_\phi(\fk) := (\R \oplus_\omega \cL_\phi(\fk)) \rtimes_D \R. $$
Here we extend $D$ to $\tilde\cL_\phi(\fk)$ by 
$D(z,\xi) := (0,\xi')$ to obtain the Lie bracket on~$\g$: 
$$ [(z_1, \xi_1, t_1), (z_2, \xi_2, t_2)] 
:= (\la \xi_1',\xi_2\ra, t_1 \xi_2' - t_2 \xi_1' + [\xi_1, \xi_2], 0). $$

To formulate our main result, we start with a 
$1$-connected simple Hilbert--Lie group $K$, i.e., $\fk$ contains no proper closed 
ideal. Corresponding 
twisted loop groups $\cL_\phi(K)$ and $\cL_\psi(K)$ 
are isomorphic if $\phi$ and $\psi$ define the same 
conjugacy class in the group $\pi_0(\Aut(K))$ of connected components of the 
Lie group $\Aut(K)$. Since 
every infinite-dimensional simple Hilbert--Lie algebra 
is isomorphic to the algebra $\fu_2(\cH)$ of skew-hermitian 
Hilbert--Schmidt operators on an infinite-dimensional real, 
complex or quaternionic Hilbert space, it is possible 
to determine $\Aut(\fk)$ explicitly (Theorem~\ref{thm:aut-grp}). 
This in turn leads to a complete classification of the corresponding twisted 
loop groups. We thus obtain four classes 
of loop algebras: the untwisted loop algebras  
$\cL(\fu_2(\cH))$, where $\cH$ is an infinite-dimensional 
real, complex or quaternionic Hilbert space, and a twisted type 
$\cL_\phi(\fu_2(\cH))$, where $\cH$ is a complex Hilbert space 
and $\phi(x) = \sigma x \sigma$ holds for an antilinear isometric involution
$\sigma \: \cH \to \cH$ (this corresponds to complex conjugation of the corresponding 
matrices). Our main result is the classification 
of the semibounded unitary representations of the 
$1$-connected Lie groups $G := \hat\cL_\phi(K)$ corresponding to the respective 
double extensions $\g = \hat\cL_\phi(\fk)$. 

To describe this classification, we choose a 
maximal abelian subspace 
\[ \ft \subeq \fk^\phi := \{ x \in \fk \: \L(\phi)x = x \}. \] 
Then $\ft_\g := \R \oplus \ft \oplus \R$ is maximal 
abelian in $\g = \hat\cL_\phi(\fk)$. We write 
$T_G := \exp(\ft_\g) \subeq G$ for the corresponding 
subgroup and identify its character group 
$\hat T_G$ with a subgroup of $i\ft_\g'$, where $'$ denotes the topological 
dual space. 
In the following, we assume that the corresponding $\Z$-graded root system 
$\Delta_\g = \Delta(\g,\ft_\g)$ is one of the seven irreducible locally affine 
reduced root systems of infinite rank 
$A_J^{(1)}, B_J^{(1)}, C_J^{(1)}, D_J^{(1)}, 
B_J^{(2)}, C_J^{(2)}$ or $BC_J^{(2)}$ (cf.~Definition~\ref{def:locaffroot}).

\begin{thmi} \mlabel{thm:1} Irreducible semibounded representations 
$\pi_\lambda$ of $G = \hat\cL_\phi(K)$ are 
characterized by their $\ft_\g$-weight set 
\[ \cP_\lambda = \conv(\hat\cW\lambda) \cap (\lambda + \hat\cQ) \subeq i \ft_\g'
\quad \mbox{ with } \quad 
\Ext(\conv(\cP_\lambda)) = \hat\cW\lambda,\]  
where $\hat\cW$ is the Weyl group of the pair $(\g,\ft_\g)$ and 
$\hat\cQ \subeq i\ft_\g'$ the corresponding root group.
The set of occurring extremal weights $\lambda$ 
is $\pm \cP^+$, where 
\[  \cP^+ := \{ \mu \in \hat T_G  \: 
\inf (\hat\cW\mu)(d) > -\infty \} \quad \mbox{ for  } \quad 
d := (0,0,-i)\in \g.\] 
Let $\cP^+_d \subeq \cP^+$ denote the set of those elements 
$\mu$ for which $\mu(d) = \min(\hat\cW\mu)(d)$. 
For $\mu_c := \mu(i,0,0)$, 
the elements $\mu \in \cP^+$ contained in $\cP_d^+$ are characterized by: 
\[ \mu_c \geq 0, \quad 
|\mu(\check\alpha)| 
\leq \frac{2\mu_c}{(\alpha,\alpha)}, \quad 
|\mu(\check\beta)| \leq \frac{4\mu_c}{(\beta,\beta)} \quad \mbox{ for } \quad 
(\alpha,1),(\beta,0) \in \Delta_\g, \] 
where $\check\alpha \in i \ft$ is the associated coroot. 
The parameter space of the equivalence classes of semibounded representations 
is given by 
\[ \pm \cP^+/\hat\cW \cong \pm \cP_d^+/\cW,\] 
where $\cW\subeq \hat\cW$ is the Weyl group of the pair $(\fk^\phi,\ft)$. 
\end{thmi}

In all cases we obtain an explicit description of the set 
$\cP_d^+/\cW$ of $\cW$ orbits in the set 
$\cP_d^+$ of $d$-minimal integral weights 
which is based on a characterization of the $d$-minimal weights 
(Theorem~\ref{thm:dmin}) and the quite elementary 
classification of $\cW$-orbits (Proposition~\ref{prop:classworb}). 
The remarkable observation that the intersection of a $\hat\cW$-orbit 
with the set $\cP_d^+$ coincides with a $\cW$-orbit is drawn 
from preliminary work on convex hulls of Weyl group orbits 
(\cite{HN12}). 

{\bf Structure of the paper:} 
We start in Section~\ref{sec:1} with the introduction of the 
simple Hilbert--Lie algebra and their root decompositions which leads to 
the four locally finite root systems $A_J, B_J, C_J$ and $D_J$ 
(cf.\ \cite{LN04}). 
Our first main result is the determination of the  
full automorphism groups of the simple Hilbert--Lie algebras (Theorem~\ref{thm:aut-grp}). 
In Section~\ref{sec:2} we introduce the double extensions 
$\g = \hat\cL_\phi(\fk)$ for the twisted loop algebras 
$\cL_\phi(\fk)$, where we restrict our attention to those automorphisms 
$\phi$ for which the corresponding root systems $\Delta_\g$ 
are, as $\Z$-graded root systems, equal to one of the seven locally 
affine root systems 
$A_J^{(1)}, B_J^{(1)}, C_J^{(1)}, D_J^{(1)}, 
B_J^{(2)}, C_J^{(2)}$ or $BC_J^{(2)}$ (cf.~Definition~\ref{def:locaffroot}). 
In Section~\ref{sec:3} we mount to the global level by showing 
that, for every $1$-connected simple Hilbert--Lie group $K$,
 there exists a $1$-connected Fr\'echet--Lie group 
$\hat\cL_\phi(K)$ which is a central $\T$-extension 
$\tilde\cL_\phi(K) \rtimes_\alpha \R$ of $\cL_\phi(K) \rtimes_\alpha \R$. 
Section~\ref{sec:4} focuses on the action of the Weyl 
group $\hat\cW$ on the integral weights. Here the main result is 
the explicit classification of the $d$-minimal weights 
in Theorem~\ref{thm:dmin}. 
After these preparations, we attack our goal of classifying the 
irreducible semibounded representations of 
$G = \hat\cL_\phi(K)$. The first major step is 
Theorem~\ref{thm:5.5}, asserting that for a semibounded 
representation $(\pi, \cH)$, the operator 
$\dd\pi(d)$ is either bounded from below (positive energy representations) 
or from above. 
Up to passing to the dual representation, we may therefore 
assume that we are in the first case. 
Then the minimal spectral value of $\dd\pi(d)$ turns out to be an eigenvalue 
and the group $Z_G(d)$ acts on the corresponding eigenspace, 
which leads to a bounded representation $(\rho,V)$ of this group. 
To proceed further, we rely on some general results concerning 
holomorphic induction. This framework has been developed for 
Banach--Lie groups in \cite{Ne12a} and in Appendix~\ref{app:c} we briefly  
explain how it can be carried over to certain 
Fr\'echet--Lie groups, containing in particular groups such as 
$\hat\cL_\phi(K)$. This permits us to conclude that 
the representation $(\rho,V)$ is irreducible and that it 
determines $(\pi, \cH)$ uniquely. Since an explicit 
classification of the bounded irreducible representations 
of the groups $Z_G(d)_0$ is available from 
\cite{Ne98, Ne11c} (Theorem~\ref{thm:classif-hilbert}) 
in terms of $\cW$-orbits of extremal weights, it remains to characterize 
those weights $\lambda$ for which the corresponding representation 
$(\rho_\lambda, V_\lambda)$ corresponds to a unitary representation 
of $G$. This is achieved in Theorem~\ref{thm:5.9} asserting that 
this is equivalent to $\lambda$ being $d$-minimal, and the final step 
consists in showing that the irreducible $G$-representation 
$(\pi_\lambda, \cH_\lambda)$ corresponding to a $d$-minimal 
weight is actually semibounded (Theorem~\ref{thm:6.1}). 
Compared to related arguments in other contexts 
(cf.\ \cite{Ne98, Ne00, Ne10b}), the argument we give here 
is rather direct and does not require any convexity results on projections 
of coadjoint orbits, such as \cite{AP83, KP84}. 
This brings us full circle and completes the proof of Theorem~\ref{thm:1}. 

For untwisted loop groups $\hat\cL(K)$ and compact groups $K$, 
the corresponding class of representations is well-known from 
the context of affine Kac--Moody algebras 
(cf.\ \cite{Ka90}, \cite{PS86}). In this context one thus obtains 
the class of positive energy representations ($\dd\pi(d)$ bounded from below), 
but this requirement is too weak for infinite-dimensional $K$. We therefore 
 work with the semiboundedness condition which has the additional advantage 
that it is invariant under twisting with arbitrary automorphisms. 
Compared with the classical situation where $K$ is finite-dimensional, we 
thus obtain the new insight that every positive energy 
representation is actually semibounded. 
In various respects our techniques are simpler than the ones used 
in the classical case to prove the existence of the unitary 
representation $(\pi_\lambda, \cH_\lambda)$ for a $d$-minimal weight~$\lambda$ 
(cf.\ \cite{PS86}, \cite{GW84}, \cite{TL99}). Instead of using 
ad hoc operator estimates for the corresponding Lie algebra 
representation, we combine the technique of holomorphic induction and 
some general results on analytic positive definite functions 
(cf.\ Appendix~\ref{app:b}) to see that the 
$d$-minimality of $\lambda$, which is already known to lead to a 
unitary Lie algebra module on the algebraic level (\cite{Ne10a}), 
to integrate to an analytic representation of the Lie group $\hat\cL_\phi(K)$. 
This is done by using the following new characterization: 
An operator-valued function 
$\phi \: V \to B(\cK)$, $\cK$ a Hilbert space, $V$ an identity neighborhood of any 
Fr\'echet--BCH--Lie group $G$, is positive definite in an identity neighborhood 
if and only if the corresponding linear map 
$U(\g_\C) \to B(\cK)$ obtained by derivatives in~$\1$ is positive definite 
(Theorem~\ref{thm:charaposdef}). 

Although it does not appear on the surface of our arguments, 
it is crucial that we deal with the Fr\'echet--Lie group 
$\hat\cL_\phi(K)$ through its Banach analog
$\hat\cL_\phi^H(K)$ constructed similarly from $H^1$-maps instead of 
smooth ones. This is a topological group which is a Banach manifold 
and a semidirect product $\tilde\cL_\phi^H(K) \rtimes_\alpha \R$, where 
the factor on the left is a Banach--Lie group but the translation 
action of $\R$ is not smooth. As a byproduct, our techniques 
imply that the representations $\pi_\lambda$ extend to continuous 
representations of $\hat\cL_\phi^H(K)$ which are analytic on 
$\tilde\cL_\phi(K)$ in the sense that the space of analytic vectors is dense 
(Remark~\ref{rem:5.9}). For the convenience of the reader, we collect 
some basic information on groups of $H^1$-maps in Appendix~\ref{app:a}, 
including the existence of the central Lie group extension 
$\tilde\cL_\phi^H(K)$.

In view of the classification of semibounded irreducible representations 
in terms of extremal weights, it is natural to ask why we first pass to the double 
extension of the group $\cL_\phi(K)$ to study unitary representations. 
Without the double extension, the representation theory of 
loop groups is much less interesting: 
From Theorem~\ref{thm:centriv} it follows 
that all semibounded unitary representations of the central 
extension  $\tilde\cL_\phi(\fk)$  
are trivial on the center and factor through bounded representations 
of $\cL(\fk)$, which in turn are finite-dimensional and tensor products of evaluation representations 
(see \cite{NS11} for the case of 
Banach--Lie algebras of maps, the Fr\'echet case will be dealt with elsewhere). 
We also show in Theorem~\ref{thm:semdirtriv}  that all semibounded 
representations of $\cL_\phi(K) \rtimes_\alpha \R$ are trivial 
on $\cL_\phi(K)$. These two results clearly demonstrate that the double extension 
of $\cL_\phi(\fk)$ is crucial to get hold of the interesting 
class of semibounded representations. 

To put our results into perspective, it is instructive to recall 
that if $X$ is a compact space and $K$ a semisimple compact Lie group, 
then all irreducible bounded unitary representations 
of $C(X,K)_0$ are  finite tensor products of evaluation 
representations, hence in particular finite-dimensional (\cite{NS11}). 
Other irreducible representations $(\pi, \cH)$ of the loop groups 
$\cL(K)\rtimes_\alpha \R$ (twisted loop modules) 
constructed by Chari and Pressley in \cite{CP86} have the property that 
the spectrum of $\dd\pi(d)$ is unbounded 
from below and above and their restrictions to $\tilde\cL(K)$ are 
not irreducible. For any, not necessary compact, Lie group $K$, 
the group $C(X,K)$ has unitary representations obtained 
as finite tensor products of evaluation representations. 
However, for some non-compact groups, such as $K = \tilde \SU_{1,n}(\C)$, 
one even has ``continuous'' tensor product representations which 
are irreducible (cf.\ \cite{JK85}, \cite{Be79}, \cite{CP87}, 
\cite{VGG74, VGG80}). 
In the algebraic context of loops group, these representations also appear in 
\cite{JK89} which contains a classification of various  types of 
unitary representations generalizing highest weight representations. 
In addition to these representations which actually extend to groups 
of measurable maps, there exist irreducible representations of mapping 
groups defined most naturally on maps of Sobolev $C^1$-maps, the so-called 
energy representations (cf.\ \cite{AH78}, \cite{Alb93}) 
and certain variants of positive energy representations of 
gauge groups of tori (cf.\ \cite{To87}). The problem to classify all smooth 
(projective) irreducible unitary representations of gauge groups is still wide 
open, although the classification of their central extensions 
by Janssens and Wockel (\cite{JW10}) is a major step towards this goal.

\subsection*{Notation:} We collect some basic notational conventions used below. 
We write $\N = \{1,2,\ldots\}$ for the natural numbers. 

Hilbert spaces over $\K \in \{\R,\C,\H\}$ are mostly denoted $\cH$. 
We write $B(\cH)$ for the algebra of all bounded operators on $\cH$, 
$B_2(\cH)$ for the ideal of {\it Hilbert-Schmidt
operators}, $B_1(\cH)$ for the ideal of {\it trace class operators} and 
$K(\cH)$ for the ideal of {\it compact operators}. 
Accordingly we write 
\[ \U_2(\cH) := \U(\cH) \cap (\1 + B_2(\cH)) \] 
for the Hilbert--Lie group of unitary operators $u$ for which 
$u - \1$ is Hilbert--Schmidt. 

Let $G$ be a Lie group (modeled on a locally convex space) and unit element~$\1$. 
Then we write $\g = \L(G)$ for its Lie algebra, which is identified with 
the tangent space $T_\1(G)$. The Lie bracket is obtained by identification with the 
Lie algebra of left invariant vector fields. 
A smooth map $\exp_G \: \g \to G$  is called an {\it exponential function} 
if each curve $\gamma_x(t) := \exp_G(tx)$ is a one-parameter group 
with $\gamma_x'(0)= x$. Not every infinite-dimensional Lie group has an 
exponential function (\cite[Ex.~II.5.5]{Ne06}), but exponential functions 
are unique whenever they exist, and this is in particular the case for all 
Banach--Lie groups.

\tableofcontents

\section{Hilbert--Lie groups} \mlabel{sec:1}

In this section we briefly introduce the class of Hilbert--Lie algebras, 
the closest infinite-dimensional relatives of compact Lie algebras. 

\subsection{Hilbert--Lie algebras} 

\begin{definition} \mlabel{def:1.1} 
(a) A {\it Hilbert--Lie algebra} $\fk$ is a real Lie algebra 
endowed with the structure of a real Hilbert space such that the 
scalar product is invariant under the adjoint action, i.e., 
\[ \la [x,y],z \ra = \la x,[y,z]\ra \quad \mbox{ for } \quad x,y,z \in \fk.\]
From the Closed Graph Theorem and the Uniform Boundedness Principle 
one derives that the bracket $\fk \times \fk \to \fk$ is continuous 
with respect to the norm topology on $\fk$ (cf.\ \cite[p.~70]{Sch60}). 
A Hilbert--Lie algebra  $\fk$ is called {\it simple} if 
$\{0\}$ and $\fk$ are the only closed ideals. 
\end{definition}

\begin{ex}
  \mlabel{nex:1.3}
(a) A finite-dimensional Lie algebra $\fk$ carries the structure of a Hilbert--Lie algebra  
if and only if it is compact. 

(b) For any Hilbert space $\cH$ over $\K \in \{\R,\C,\H\}$, the Lie algebra 
\[ \fu_2(\cH) := \{ x \in B_2(\cH) : x^* = - x\} \] 
is a Hilbert--Lie algebra with respect to the scalar product 
$\la x,y \ra := \tr_\R(xy^*) = - \tr_\R(xy)$. It is simple if $\dim \cH = \infty$. 
\end{ex} 

\begin{thm} \mlabel{thm:1.3} {\rm(Schue)} Each Hilbert--Lie algebra 
$\fk$ is an orthogonal direct sum 
$\fk = \z(\fk) \oplus \hat\bigoplus_{j \in J} \fk_j$, where 
each $\fk_j$ is a simple ideal. Each simple infinite-dimensional 
Hilbert--Lie algebra is isomorphic to 
$\fu_2(\cH)$ for an infinite-dimensional real, complex or 
quaternionic Hilbert space $\cH$. 
\end{thm}

\begin{prf} The orthogonal decomposition into center and simple ideals follows from 
\cite[1.2, Th.\ 1]{Sch60}. 
The classification of the simple Hilbert algebras $\fk$ follows immediately from 
the classification of the complex  $L^*$-algebras because 
$\fk_\C$ is a complex $L^*$-algebra. For the separable case, the 
classification was obtained in \cite[3.7, Th.\ 3]{Sch60} 
under the assumption of the existence of a root decomposition whose 
existence was shown in \cite{Sch61}. The classification was extended to the
non-separable case in \cite{CGM90}, \cite{Neh93} and \cite[Thm.~19.28]{St99}. 
\end{prf}

\begin{defn} If $\cH$ is a real Hilbert space, then we also write 
\[ \OO(\cH) := \U(\cH), \quad \OO_2(\cH) := \U_2(\cH), \quad 
\fo(\cH) := \fu(\cH), \quad  \quad \mbox{ and } \quad 
\fo_2(\cH) := \fu_2(\cH).\] 
For a quaternionic Hilbert space $\cH$, we write
\[ \Sp(\cH) := \U(\cH), \quad \Sp_2(\cH) := \U_2(\cH), \quad 
\fsp(\cH) := \fu(\cH) \quad \mbox{ and } \quad 
\fsp_2(\cH) := \fu_2(\cH).\] 
\end{defn}

\begin{thm}\mlabel{thm:1.5} (cf.\ \cite[Sect.~II.4]{Ne02a}) 
For an infinite-dimensional Hilbert space $\cH$, over $\R$, $\C$, resp., $\H$, 
the homotopy groups of $\OO_2(\cH)$, $\U_2(\cH)$, resp., $\Sp_2(\cH)$ 
are given by: 
\[ 
\begin{array}{c||c|c|c|}
& \OO_2(\cH) &  \U_2(\cH)  & \Sp_2(\cH)  \\  \hline\hline
\pi_0 & \Z/2 &  \{\1\}  &  \{\1\}  \\ \hline
\pi_1   &  \Z/2  &  \Z & \{\1\}  \\ \hline
\pi_2 & \{\1\} & \{\1\}& \{\1\}  \\ \hline
\pi_3 & \Z & \Z& \Z  \\
\hline
\end{array}
\]   
\end{thm}

\begin{rem} \mlabel{rem:1.6} 
If $\cH$ is infinite-dimensional, Schur's Lemma implies that the 
center of the groups $\OO_2(\cH), \U_2(\cH)$ and $\Sp_2(\cH)$ is trivial. 
Therefore their fundamental group is isomorphic to the center of the simply 
connected covering group, so that 
\[ Z(\tilde \OO_2(\cH)_0) \cong \pi_1(\OO_2(\cH)) \cong \Z/2 \quad \mbox{ and } \quad 
Z(\tilde \U_2(\cH)) \cong \pi_1(\U_2(\cH)) \cong \Z.\] 
\end{rem}

\subsection{Root decomposition} 

Our parametrization of irreducible semibounded representations is based on weights w.r.t.\ 
a maximal abelian subalgebra. In this subsection we recall some basics on roots and root 
space decompositions. 

\begin{defn} \mlabel{def:basic2} 
(a) Let $\g$ be a real topological Lie algebra and $\g_\C$ be its complexification. 
If $\sigma \: \g_\C \to \g_\C, z = x + iy \mapsto \oline z = x - i y,$ denotes the complex conjugation with 
respect to $\g$, we write $x^* := -\sigma(x)$ for $x \in \g_\C$, so that 
$\g = \{ x \in \g_\C \: x^* = -x\}$. 
Let $\ft \subeq \g$ be a maximal abelian subalgebra 
and $\ft_\C\subeq \g_\C$ be its complexification. 
For a linear functional $\alpha \in \ft_\C'$ (the space of $\C$-valued continuous 
linear functionals on $\ft$ which is identified with the space of $\C$-linear continuous functionals 
on $\ft_\C$),  
\[ \g_\C^\alpha = \{ x \in \g_\C \: (\forall h \in \ft_\C)\ 
[h,x]= \alpha(h)x\}\] 
is called the corresponding {\it root space}, and 
\[ \Delta :=  \Delta(\g,\ft) := \{ \alpha \in \ft_\C^* \setminus \{0\} \: \g_\C^\alpha
\not= \{0\}\}\] 
is the {\it root system} of the pair $(\g,\ft)$. 
We then have $\g_\C^0 = \ft_\C$ and $[\g_\C^\alpha, \g_\C^\beta] 
\subeq \g_\C^{\alpha+\beta}$, hence in particular 
$[\g_\C^\alpha, \g_\C^{-\alpha}] \subeq \ft_\C$. 

(b) If $\g$ is the Lie algebra of a group $G$ with an exponential function, 
then we call $\ft$ {\it elliptic} if the subgroup 
$e^{\ad \ft} = \Ad(\exp \ft) \subeq \Aut(\g)$ is equicontinuous. 
We then have 
\begin{itemize}
\item[\rm(I1)] $\alpha(\ft) \subeq i \R$ for $\alpha \in \Delta$, 
and therefore 
\item[\rm(I2)] $\sigma(\g_\C^\alpha) = \g_\C^{-\alpha}$ 
for $\alpha \in \Delta$. 
\end{itemize}
\end{defn}

\begin{lem}
  \mlabel{lem:e.1} Suppose that $\ft \subeq\g$ is elliptic. 
For $0 \not= x_{\alpha} \in \g_\C^{\alpha}$, the subalgebra 
$\g_\C(x_\alpha) := \Spann_\C\{x_\alpha, x_{\alpha}^*,
[x_\alpha, x_{\alpha}^*]\}$
is $\sigma$-invariant and of one of the following types: 
\begin{description}
\item[\rm(A)] The abelian type: $[x_\alpha, x_{\alpha}^*] = 0$, i.e., 
$\g_\C(x_\alpha)$ is two-dimensional abelian. 
\item[\rm(N)] The nilpotent type: $[x_\alpha, x_{\alpha}^*] \not= 0$ 
and $\alpha([x_\alpha, x_{\alpha}^*]) = 0$, i.e., 
$\g_\C(x_\alpha)$ is a three-dimensional Heisenberg algebra. 
\item[\rm(S)] The simple type: $\alpha([x_\alpha, x_{\alpha}^*]) \not= 0$,
i.e., $\g_\C(x_\alpha) \cong \fsl_2(\C)$. In this case we distinguish 
two cases: 
\begin{description}
\item[\rm(CS)] $\alpha([x_\alpha, x_{\alpha}^*]) > 0$, i.e., 
$\g_\C(x_\alpha) \cap \g \cong \su_2(\C)$, and 
\item[\rm(NS)] $\alpha([x_\alpha, x_{\alpha}^*]) < 0$, i.e., 
$\g_\C(x_\alpha) \cap \g \cong \su_{1,1}(\C) \cong \fsl_2(\R)$. 
\end{description}
\end{description}
\end{lem} 

\begin{prf} (cf.\ \cite[App.~C]{Ne10b})
First we note that, in view of $x_\alpha^* \in \g_\C^{-\alpha}$,  
\cite[Lemma~I.2]{Ne98} applies, and we see that $\g_\C(x_\alpha)$ is of one of
the three types (A), (N) and (S). We note that 
$\alpha([x_\alpha, x_\alpha^*]) \in \R$ 
because of (I2) and $[x_\alpha, x_\alpha^*] \in i\ft$. 
Now it is easy to check that 
$\g_\C(x_\alpha) \cap \g$ is of type (CS), resp., (NS), according to the sign
of this number. 
\end{prf} 

\begin{defn} \mlabel{def:1.8} 
(a) Assume that $\g_\C^\alpha = \C x_\alpha$ is one-dimensional 
and that $\g_\C(x_\alpha)$ is of type (S). 
Then there exists a unique element 
$\check \alpha \in \ft_\C \cap [\g_\C^\alpha, \g_\C^{-\alpha}]$ with 
$\alpha(\check \alpha) = 2$. It is called the {\it coroot of $\alpha$}. 
The root $\alpha \in \Delta$ is said to
be {\it compact} if, for $0\not= x_\alpha \in \g_\C^\alpha$, we have 
$\alpha([x_\alpha,x_\alpha^*]) > 0$ and {\it non-compact} otherwise. 
We write $\Delta_c$ for the set of compact roots. 
With the notation $\R_+ := [0,\infty[$, 
Lemma~\ref{lem:e.1} implies that 
\begin{equation}
  \label{eq:corootsign}
\check \alpha \in \R_+ [x_\alpha,x_\alpha^*] \quad \mbox{ for } \quad 
\alpha \in \Delta_c. 
\end{equation}

(b) The {\it Weyl group} $\cW = \cW(\g,\ft) \subeq \GL(\ft_\C)$ is the subgroup generated 
by all reflections 
\begin{equation}
  \label{eq:ref1}
 r_\alpha(x) := x - \alpha(x) \check \alpha \quad \mbox{ for compact roots} \quad 
\alpha \in \Delta_c. 
\end{equation}
It acts on the dual space $\ft_\C^*$ by the dual maps 
$r_\alpha^*(\beta) := \beta - \beta(\check \alpha) \alpha.$ 

(c) A linear functional $\lambda \in i \ft'$ is said to be an  {\it integral
 weight} 
if $\lambda(\check \alpha) \in \Z$ holds for every compact root 
$\alpha \in \Delta_c$. We write $\cP = \cP(\g,\ft) \subeq i\ft'$ for the group of all 
integral weights. 
\end{defn}

Let $\fk$ be a Hilbert--Lie algebra and 
$\ft \subeq \fk$ be a maximal abelian subalgebra. 
According to \cite{Sch61}, $\ft_\C \subeq \fk_\C$ defines a 
root space decomposition 
$$ \fk_\C = \ft_\C \oplus \hat\bigoplus_{\alpha \in \Delta} \fk_\C^\alpha $$
which is a Hilbert space direct sum with respect to the hermitian extension of the scalar product 
to~$\fk_\C$. 
We now describe the relevant root data for the three types of 
simple Hilbert algebras $\fu_2(\cH)$, where $\cH$ is a 
Hilbert space over $\K \in \{\R,\C,\H\}$. 

\begin{ex} \mlabel{ex:d.1a} (cf.\ \cite[Ex.~C.4]{Ne11c}) 
(Root data of unitary Lie algebras)  
Let $\cH$ be a complex Hilbert space with 
orthonormal basis $(e_j)_{j \in J}$ 
and $\ft \subeq \fk := \fu_2(\cH)$ be the 
subalgebra of all diagonal operators with respect to the $e_j$. 
Then $\ft$ is elliptic and maximal abelian, 
$\ft_\C \cong \ell^2(J,\C)$. The set of 
roots of $\fk_\C\cong \gl_2(\cH)$ with respect to $\ft_\C$  is given by the root 
system 
\[ \Delta = \{ \eps_j - \eps_k \: j\not= k \in J \} =: A_J.\]
Here the operator 
$E_{jk} e_m := \delta_{km} e_j$ is a $\ft_\C$-eigenvector 
in $\gl_2(\cH)$ generating the corresponding eigenspace 
and $\eps_j(\diag(h_k)_{k \in J}) = h_j$. 
From $E_{jk}^* = E_{kj}$ it follows that 
\[  (\eps_j - \eps_k)\,\check{}= E_{jj} - E_{kk} 
= [E_{jk}, E_{kj}] = [E_{jk}, E_{jk}^*],\] 
so that $\Delta = \Delta_c$, i.e., all roots are compact.

The Weyl group $\cW$ is isomorphic to the group $S_{(J)}$ of finite 
permutations of $J$, acting in the canonical  way on $\ft_\C \cong \ell^2(J,\C)$. 
It is generated by the reflections $r_{jk} := r_{\eps_j - \eps_k}$ 
corresponding to the transpositions of $j \not= k \in J$. 
The Weyl group acts transitively on the set of roots and, in particular, 
all roots have the same length~$2$ w.r.t.~the scalar product induced by 
\ $\la x,y \ra = \tr(xy^*)$ on the dual space.
\end{ex}

\begin{rem}
(a) In many situations it is convenient to describe 
real Hilbert spaces as pairs 
$(\cH,\sigma)$, where $\cH$ is a complex Hilbert space 
and $\sigma \: \cH\to \cH$ is a {\it conjugation}, i.e., 
an antilinear isometry with $\sigma^2 = \id_\cH$. Then we write 
$A^\top := \sigma A^* \sigma$, which corresponds to the transposition of matrices 
with respect to any ONB contained in $\cH^\sigma$. 

(b) A quaternionic Hilbert space $\cH$ can be considered as a 
complex Hilbert space $\cH^\C$ (the underlying complex Hilbert 
space), endowed with an {\it anticonjugation} $\sigma$, i.e., 
$\sigma$ is an antilinear isometry with $\sigma^2 = -\1$. 

(c) That all conjugations and anticonjugations on a complex Hilbert space are conjugate under 
the unitary group $\U(\cH)$ has been shown in \cite{Ba69} by describing them in terms 
of orthonormal bases. 
\end{rem}

\begin{ex} \mlabel{ex:d.1b} (cf.\ \cite[Ex.~C.5]{Ne11c}) 
(Root data of symplectic Lie algebras)  
For a complex Hilbert space $\cH$ with a conjugation $\sigma$, 
we consider the quaternionic Hilbert space $\cH_\H := \cH^2$, where the 
quaternionic structure is defined by the anticonjugation 
$\tilde\sigma(v,w) := (\sigma w, -\sigma v)$. 
Then $\fk := \sp_2(\cH_\H) 
= \{ x \in \fu_2(\cH^2) \: \tilde\sigma x = x \tilde\sigma \}$ 
and 
\[ \sp_2(\cH_\H)_\C = \left\{ \pmat{ A & B \cr C & -A^\top \cr} 
\in B_2({\cal H}^2) \: B^\top = B,  C^\top  = C\right\}. \] 
Let $(e_j)_{j \in J}$ be an orthonormal basis of $\cH$ with 
$\sigma(e_j) = e_j$ for every $j$, and 
$\ft \subeq \fk$ be the 
subalgebra of all diagonal operators with respect to the basis elements  
$(e_j,0)$ and $(0,e_k)$ of $\cH^2$. Then $\ft$ 
is elliptic and maximal abelian in $\fk$. 
Moreover, $\ft_\C \cong \ell^2(J,\C)$  
consists of diagonal operators of the form 
$h = \diag((h_j), (-h_j))$, and the set of 
roots of $\fk_\C$ with respect to $\ft_\C$  is given by 
\[ \Delta = \{ \pm 2 \eps_j, \pm (\eps_j \pm \eps_k) \: j \not= k, j,k
\in J \} =: C_J,\] 
where $\eps_j(h)  = h_j$. If we write 
$E_j = \pmat{E_{jj} & 0 \\ 0 & -E_{jj}}\in \ft_\C$ for the element defined by 
$\eps_k(E_j) = \delta_{jk}$, then the coroots are given by 
\begin{equation}
  \label{eq:CJcoroot}
 (\eps_j \pm  \eps_k)\,\check{}= E_j \pm  E_k 
\quad \mbox{ for } \quad j \not=k 
\quad \mbox{ and } \quad 
(2\eps_j)\,\check{}= E_j.
\end{equation}
Again, all roots are compact, and the 
Weyl group $\cW$ is isomorphic to the group $\{\pm 1\}^{(J)} \rtimes S_{(J)}$,
where $\{\pm 1\}^{(J)}$ is the group of finite 
sign changes on $\ell^2(J,\R)$. In fact, 
the reflection $r_{\eps_j - \eps_k}$ acts as a transposition and 
the reflection $r_{2\eps_j}$ changes the sign of the $j$th component. 
The Weyl group has two orbits in $C_J$, the short roots form a root system 
of type $D_J$ and the second orbit is the set 
$\{ \pm 2 \eps_j \: j \in J\}$ of long roots. 
\end{ex}

\begin{ex} \mlabel{ex:d.1c} (cf.\ \cite[Ex.~C.6]{Ne11c}) 
(Root data of orthogonal Lie algebras)  
Let $\cH_\R$ be an  infinite-dimensional real Hilbert space 
and $\fk := \fo_2(\cH_\R)$ be the corresponding simple Hilbert--Lie 
algebra. Let $\ft \subeq \fk$ be maximal abelian. 
The fact that $\ft$ is maximal 
abelian implies that the common kernel $\ker \ft$ is at most 
one-dimensional. Since $\ft$ consists of compact skew-symmetric operators, 
the Lie algebra $\ft_\C$ is diagonalizable on the complexification 
$\cH :=(\cH_\R)_\C$. We conclude that, on 
the space $(\ker \ft)^\bot$, we have an orthogonal complex 
structure $I$ commuting with $\ft$ and there exists an orthonormal 
subset $(e_j)_{j \in J}$ of $(\ker \ft)^\bot$ such that 
$\{ e_j, Ie_j\: j \in J\}$ is an orthonormal basis of $(\ker \ft)^\bot$ 
and all the planes $\R e_j + \R I e_j$ are $\ft$-invariant. 
If $\cH_\R^{\ft}$ is non-zero, we write $e_{j_0}$ for a unit vector 
in this space. For $j \in J$ put 
\[ f_j := \frac{1}{\sqrt 2}(e_j - iI e_j) \quad \mbox{ and } 
\quad f_{-j} := \frac{1}{\sqrt 2}(e_j + iI e_j).\] 
If $\ker \ft\not=\{0\}$, then we also put $f_{j_0} := e_{j_0}$. 
Then the $f_j$ form an orthonormal basis of $\cH$ consisting 
of $\ft$-eigenvectors. 

We conclude that $\ft_\C$ is precisely the set of all those elements 
in $\fk_\C = \fo_2(\cH_\R)_\C$ which are diagonal with respect to the ONB consisting 
of the $f_j$. This implies that $\ft_\C \cong \ell^2(J,\C)$, where 
$x \in \ft_\C$ corresponds to the element $(x_j)_{j \in J} \in \ell^2(J,\C)$ 
defined by $x f_j = x_j f_j$, $j \in J$. 
Writing $\eps_j(x) := x_j$, we see that 
$\{ \pm \eps_j : j \in J\}$, together with $\eps_{j_0}$ if 
$\ker \ft \not=\{0\}$, is the set of $\ft_\C$-weights of $\cH$. 
Accordingly, the set of 
roots of $\fk_\C$ with respect to $\ft_\C$  is given by 
\[ \Delta = \{ \pm \eps_j \pm \eps_k \: j \not= k, j,k\in J \} =: D_J
\quad \mbox{ if } \quad \ker \ft =\{0\}, \] 
and 
\[ \Delta = \{ \pm \eps_j \pm \eps_k \: j \not= k, j,k\in J \} 
\cup \{ \pm\eps_j \: j \in J\} =: B_J \quad \mbox{ otherwise}.\] 

As in Example~\ref{ex:d.1b}, all roots are compact. 
For $B_J$ we obtain the same Weyl group 
$\{ \pm 1\}^{(J)} \rtimes S_{(J)}$ as for $C_J$. 
For $D_J$ the reflection $r_{\eps_j + \eps_k}$ changes the sign of the 
$j$- and the $k$-component, so that the 
Weyl group $\cW$ is isomorphic to the group $\{\pm 1\}^{(J)}_{\rm even} \rtimes S_{(J)}$,
where $\{\pm 1\}^{(J)}_{\rm even}$ is the group of finite even sign changes. 
For $D_J$ the Weyl group acts transitively on the set of roots and 
all roots have the same length. 
For $B_J$ we have two $\cW$-orbits,  the roots $\pm \eps_j$ are short 
and the roots $\pm \eps_j \pm \eps_k$, $j \not=k$, are long. 
\end{ex}

\begin{rem} \mlabel{rem:1.12} 
In a simple Hilbert--Lie algebra, two maximal abelian 
subalgebras are conjugate under the full automorphism group if and only 
if the corresponding root systems are isomorphic 
(see \cite[Prop.~19.24, Rem.~19.25]{St99} and \cite{Ba69}). 
Up to conjugacy by automorphisms, the classification of locally finite 
root systems thus implies that we have only four 
types of pairs $(\fk,\ft)$ and that they correspond to the root 
systems $A_J$, $B_J$, $C_J$ and $D_J$. 

For a real Hilbert space $\cH$, the Hilbert--Lie algebra $\fo_2(\cH)$ 
contains two conjugacy classes of maximal abelian 
subalgebras $\ft$, distinguished by 
$\dim \cH^{\ft} \in \{0,1\}$. 
For the classification purposes in this 
paper, we only need one maximal abelian subalgebra 
to set up the parametrization 
of the equivalence classes of unitary representations (cf.\ Theorem~\ref{thm:1}). 
Passing to a 
different conjugacy class of maximal abelian subalgebras leads to 
a different parameter space for the same class of representations. 
\end{rem}

\subsection{Automorphism groups} 

For a complex Hilbert space $\cH$, we write 
$\AU(\cH)$ for the group of unitary or antiunitary isometries 
of $\cH$ and 
\[ \PAU(\cH) := \AU(\cH)/\T\1 \cong \PU(\cH) \rtimes \{\id,\sigma\}, \] 
where  $\sigma$ is an anticonjugation of $\cH$. 

\begin{thm}
  \mlabel{thm:aut-grp}
The automorphism groups of the simple infinite-dimensional 
Hilbert algebras are given by 
\[  \Aut(\fu_2(\cH)) \cong \PAU(\cH) \] 
for a complex Hilbert space, and for the real and quaternionic case we have 
\[  \Aut(\fo_2(\cH)) \cong \OO(\cH)/\{\pm \1\} \quad \mbox{ and } \quad 
 \Aut(\fsp_2(\cH)) \cong \Sp(\cH)/\{\pm \1\}.\] 
\end{thm} 

\begin{prf} We know from Schue's Theorem~\ref{thm:1.3} that 
any simple infinite-dimensional Hilbert--Lie algebra $\fk$ is isomorphic 
to $\fu_2(\cH)$ for some infinite-dimensional 
Hilbert space $\cH \cong \ell^2(J,\K)$ with $\K \in \{\R,\C,\H\}$. 
Let $\ft \subeq \fk$ be a maximal 
abelian subalgebra and $\Delta$ be the corresponding root system. 
As we have seen in Examples~\ref{ex:d.1a}, \ref{ex:d.1b} and \ref{ex:d.1c}, 
it is of type $A_J$ (for $\K = \C$), 
$C_J$ (for $\K = \H$)  or $B_J, D_J$ (for $\K = \R$).

For any $\phi \in \Aut(\fk)$, the subspace 
$\phi(\ft) \subeq \fk$ is also maximal abelian 
with isomorphic root system. 
It follows from \cite[Thm.~2]{Ba69} (and its proof) 
that for real, complex and quaternionic Hilbert spaces, 
the group $\UU(\cH)$ acts transitive 
on the set of all maximal abelian subalgebras of $\fk = \fu_2(\cH)$ 
whose root system is 
of a given type (see also \cite[Thm.~19.24]{St99}). This implies the existence of $\psi \in 
\UU(\cH)$ for which the corresponding automorphism 
$c_\psi(x) := \psi x \psi^{-1}$ satisfies 
$c_{\psi}(\ft) = \phi(\ft)$. Then $c_{\psi}^{-1}\circ \phi$ fixes $\ft$, 
hence induces an automorphism of the root system $\Delta$. 

For each root system, the automorphism group is known  
from \cite[Props.~5.1-5.4]{St01}: 
$$ \Aut(A_J) \cong S_J \times \{\pm \id\}, \quad
\Aut(B_J) \cong \Aut(C_J) \cong \Aut(D_J) \cong (\Z/2)^J \rtimes S_J. $$
From this description it easily follows that, for 
$\Delta$ of type $B_J$, $C_J$ or $D_J$, each automorphism of the 
root system is implemented by conjugation with an element 
of the corresponding full group $\UU(\cH)$. 
For $A_J$, the elements of $S_J$ are obtained by conjugation with a unitary operator permuting 
the elements $(e_j)_{j \in J}$ of an orthonormal basis, 
and $-\id$ is obtained by $\phi(x)= \sigma x \sigma$, where 
$\sigma$ is a conjugation fixing each $e_j$, $j \in J$. 

For the realization of $\sp_2(\cH_\H)$ as in Example~\ref{ex:d.1b}, we obtain the 
elements of $S_J$ by conjugation with unitary operators of the form 
$\pmat{U & 0 \\ 0 & U}$, where $U \in \U(\cH)$ permutes the elements of the 
ONB $(e_j)_{j \in J}$. To implement the sign changes corresponding to the 
element $\chi \in \{\pm 1\}^J$ taking the value $-1$ on the subset $M \subeq J$ 
and $1$ elsewhere, we consider the projection $P$ on $\cH$ with 
\[ Pe_j =
\begin{cases}
 e_j & \text{ for } j \in M \\ 
0  & \text{ for } j \in J \setminus M. 
\end{cases}\] 
Then conjugation with $u := \pmat{ \1 - P & P \\ -P & \1 - P} \in \Sp(\cH_\bH)$ 
induces on $\Delta = B_J$ the automorphism corresponding to $\chi$. 

For the realization of $\fo_2(\cH_\R)$ as in Example~\ref{ex:d.1c}, we obtain the 
elements $\pi \in S_J$ by conjugation with orthogonal operators 
$U$ satisfying $Ue_j = e_{\pi(j)}$ for $j \in J$ and commuting with the complex structure~$I$. 
To implement the sign changes represented by $\chi \in \{\pm 1\}^J$ as above, 
we conjugate with $U \in \OO(\cH_\R)$ satisfying $U e_j = e_j$ for $j \in J$ and 
\[ U I e_j =
\begin{cases}
 -I e_j & \text{ for } j \in M \\ 
I e_j  & \text{ for } j \in J \setminus M. 
\end{cases}\]

This reduces the problem to show that every automorphism is given by 
conjugation with a unitary (or antiunitary operator in the complex case) 
to the special case where it preserves $\ft$ and 
induces the trivial automorphism on $\Delta$. In view of 
\cite[Lemma~6.2(b)]{St01}, any such automorphism $\phi$ satisfies 
$\phi(x_\alpha) = \chi(\alpha) x_\alpha$ for $x_\alpha \in \fk_\C^\alpha$ and a 
homomorphism $\chi \: \cQ := \la \Delta \ra_{\rm grp} 
\to \C^\times$. As $\phi$ is supposed to 
be isometric, we have $\im(\chi) \subeq \T$. Conversely, the orthogonality 
of the root decomposition of $\fk_\C$ implies that every homomorphism 
$\chi \: \cQ \to \T$ occurs. 
We now show that any such automorphism is given by conjugation with an element of 
$U \in \U(\cH)$. For $\Delta = A_J$ we pick an element $j_0 \in J$ and put 
\[ u_{j_0} := 1 \quad \mbox{ and } \quad u_j := \chi(\eps_j - \eps_{j_0}) \quad \mbox{ for } \quad 
j \not=j_0 \] 
to find the required element $U = \diag((u_j)) \in \U(\cH)$. 
For $C_J$ we first extend $\chi$ to the $\Z$-span of $\{ \eps_j \: j \in J\}$ 
($\T$ is divisible) and put $U := \diag(\chi(\eps_j), \chi(-\eps_j))$. 
For $B_J$ we use the same element $U \in \U((\cH_\R)_\C)$, and 
for $D_J$ we first extend $\chi$ to the $\Z$-span of $B_J$ and proceed with 
the diagonal operator $U$ with $u_{j_0} = 1$ and 
$u_{\pm j} = \chi(\eps_j)^{\pm 1}$ for $j \in J$ 
(cf.\ \cite[Sect.~6]{St01} for a similar argument in the  algebraic context). 

Finally, we note that, if $g \in \UU(\cH)$ induces the trivial automorphism 
$c_g = \id$ on $\fu_2(\cH)$, 
then $g \in \T \1$ in the complex case and 
$g \in \{ \pm \1\}$ in the real and quaternionic case. 
We thus obtain  
$\Aut(\fo_2(\cH)) \cong \OO(\cH)/\{\pm \1\},$
$\Aut(\fsp_2(\cH)) \cong \Sp(\cH)/\{\pm \1\},$ and 
$\Aut(\fu_2(\cH)) \cong \PAU(\cH)$. 
\end{prf}

\begin{cor} \mlabel{cor:2.7} 
The automorphism groups of 
$\fo_2(\cH)$ and $\fsp_2(\cH)$ are connected, 
whereas the automorphism group of 
$\fu_2(\cH)$ ($\cH$ a complex Hilbert space), 
has two connected components. 
\end{cor}

\begin{prf} It only remains to recall that the groups 
$\OO(\cH)$, $\Sp(\cH)$ and $\U(\cH)$ are connected for any 
infinite-dimensional real, quaternionic, resp., complex Hilbert space 
(cf.\ \cite[Thm.~II.6]{Ne02a}). 
\end{prf}

\begin{rem}\mlabel{rem:2.8}
The preceding corollary implies in particular that 
for $\fk = \fo_2(\cH), \fsp_2(\cH)$, each automorphism of $\fk$ acts trivially 
on the homotopy groups $\pi_j(K), j \in \N,$ for 
any connected Lie group $K$ with Lie algebra $\fk$. 
\end{rem}

\begin{prop} \mlabel{prop:1.18} Let $\cH$ be a complex Hilbert space 
and  $\sigma$ be a conjugation of $\cH$. 
For $\fk = \fu_2(\cH)$ and the automorphism $\phi(x) = \sigma x\sigma$ of $\U_2(\cH)$, we then  
have 
\[  \pi_{2k-1}(\phi) = (-1)^k \id \quad \mbox{ for } \quad k \in \N. \] 
For any connected Lie group $K$ with the Lie algebra $\fk = \fu_2(\cH)$,  
$\pi_1(K)$ either is trivial or isomorphic to $\Z$, there 
exists a $\phi_K \in \Aut(K)$ with 
$\L(\phi_K) = \phi$, and this automorphism satisfies $\pi_1(\phi_K) = -\id$. 
\end{prop} 

\begin{prf} We pick an orthonormal basis $(e_j)_{j \in J}$ in $\cH$ fixed pointwise by 
$\sigma$, and represent operators accordingly as matrices. Then the 
involution $\phi$ is given by component-wise 
conjugation $\phi(x_{ij}) = (\oline{x_{ij}})$.  
Using the approximation techniques described in 
\cite[Thm.~II.14, Cor.~II.15]{Ne02a}, it suffices to study the action of 
$\phi$ on the subgroup $\UU_n(\C)$, fixing all but $n$ basis vectors. 
Therefore it follows from \cite[Prop.~19]{Kue06} that 
$\pi_{2k-1}(\phi) = (-1)^k \id.$

If $K$ is a connected Lie group with $\L(K) \cong \fu_2(\cH)$ which 
is not simply connected, then it is a quotient 
of $\tilde\UU_2(\cH)$ by an infinite cyclic group 
because $Z(\tilde\UU_2(\cH)) \cong \Z$ by Remark~\ref{rem:1.6}. 
The automorphism $\phi$ of $\fu_2(\cH)$ induces an automorphism 
$\tilde\phi_K$ of $\tilde\UU_2(\cH)$, and this automorphism 
preserves the center. In view of $\Aut(\Z) = \{\pm \id_\Z\}$, it also 
preserves all subgroups of the center. We conclude that 
it also induces an automorphism $\phi_K$ on $K$. Since $\tilde\phi_K$ acts on
 $Z(\tilde\UU_2(\cH))\cong \pi_1(\U_2(\cH))$ by inversion, 
we obtain  $\pi_1(\phi_K) = -\id$. 
\end{prf}

Combining Remark~\ref{rem:2.8} with the preceding proposition, 
we obtain in particular 

\begin{cor} \mlabel{cor:2.10} If $K$ is a Hilbert--Lie group for which 
$\fk$ is  simple and $\phi \in \Aut(K)$, then 
$\pi_3(\phi) = \id$. 
\end{cor}

\section{Double extensions of twisted loop algebras} \mlabel{sec:2}

In this section  we introduce the double extensions 
$\g = \hat\cL_\phi(\fk)$ for the twisted loop algebras 
$\cL_\phi(\fk)$, where we restrict our attention to those automorphisms 
$\phi$ for which the corresponding root systems $\Delta_\g$ 
are, as $\Z$-graded root systems, equal to one of the seven locally 
affine root systems 
$A_J^{(1)}, B_J^{(1)}, C_J^{(1)}, D_J^{(1)}, 
B_J^{(2)}, C_J^{(2)}$ or $BC_J^{(2)}$ (cf.~Definition~\ref{def:locaffroot}). 

\subsection{Root decomposition of double extensions} 

\begin{defn} A {\it quadratic (topological) Lie algebra} is a pair $(\g,\kappa)$, 
consisting of a topological Lie algebra $\g$ and a non-degenerate invariant symmetric 
continuous bilinear form $\kappa$ on $\g$. 
Suppose that $\ft \subeq \g$ 
is maximal abelian and elliptic and that $\sum_\alpha \g_\C^\alpha$ is dense in $\g_\C$. 
We extend $\kappa$ to a hermitian 
form on $\g_\C$ which is also denoted~$\kappa$.
Then the root spaces satisfy  
$\kappa(\g_\C^\alpha, \g_\C^\beta) = \{0\}$ for $\alpha \not= \beta$ (cf.\ Definition~\ref{def:basic2}), and 
in particular, $\kappa$ is non-degenerate on $\ft_\C = \g_\C^0$ and all the root spaces 
$\g_\C^\alpha$. We thus obtain an injective antilinear map $\flat \: \ft_\C \to \ft_\C', h \mapsto 
h^\flat, h^\flat(x) := \kappa(x,h)$, where $\ft_\C'$ denotes the space of continuous linear 
functionals on $\ft_\C$. For $\alpha\in 
\ft_\C^\flat := \flat(\ft_\C)$ 
we put $\alpha^\sharp := \flat^{-1}(\alpha)$ and define a hermitian 
form on on $\ft_\C^\flat$ by 
\begin{equation}
  \label{eq:rootrel}
 (\alpha,\beta) := \kappa(\beta^\sharp, \alpha^\sharp) 
= \alpha(\beta^\sharp) = \oline{\beta(\alpha^\sharp)}.
\end{equation}
For $h \in \ft_\C$ and $x_\alpha \in \g_\C^\alpha$ we then have  
\begin{equation}
  \label{eq:brack-rel}
\alpha(h)\kappa(x_\alpha, x_{\alpha})
= \kappa([h,x_\alpha], x_{\alpha})
= \kappa(h,[x_\alpha, x_{\alpha}^*]). 
\end{equation}
Since $\kappa$ is non-degenerate on $\g_\C^\alpha$, we may choose $x_\alpha$ such that 
$\kappa(x_\alpha, x_\alpha) \not=0$. Then  the non-degeneracy of $\kappa$ on $\ft_\C$ 
leads to $\alpha\in \ft_\C^\flat$ and 
\begin{equation} \label{eq:8} 
[x_\alpha, x_{\alpha}^*] = \kappa(x_\alpha, x_{\alpha}) \alpha^\sharp.   
\end{equation}
This shows that $\Delta \subeq \ft_\C^\flat$, so that 
$(\alpha,\beta)$ is defined for $\alpha,\beta \in \Delta$ by \eqref{eq:rootrel}.
\end{defn}  

\begin{rem}
If $\alpha$ is compact and $\check \alpha = [x_\alpha, x_{\alpha}^*]$  (cf.\ 
Definition~\ref{def:1.8}), 
then \eqref{eq:brack-rel} 
and $\alpha(\check \alpha) =2$ imply 
$\kappa(\check \alpha, \check \alpha) = 2\kappa(x_\alpha, x_{\alpha})$, 
which leads for $\beta \in \ft_\C^\flat$ to 
\begin{equation}
  \label{eq:intflat}
\alpha^\sharp = \frac{2 \check \alpha}{\kappa(\check \alpha, 
\check \alpha)}, \quad 
(\alpha, \alpha) = \frac{4}{\kappa(\check \alpha, \check \alpha)} \quad 
\mbox{ and } \quad 
(\beta, \alpha) = \frac{2\beta(\check \alpha)}{\kappa(\check\alpha, 
\check \alpha)}.
\end{equation}
\end{rem}

\begin{defn} \label{def:doubext} (Double extensions) Let $(\g,\kappa)$ 
be a  real quadratic Fr\'echet--Lie algebra 
and $D \in \der(\g,\kappa)$ 
be a derivation which is skew-symmetric with respect to $\kappa$. 
Then $\omega_D(x,y) := \kappa(Dx,y)$ defines a continuous $2$-cocycle 
on $\g$, and $D$ extends to a derivation $\tilde D(z,x) := (0,Dx)$ 
of the corresponding central extension 
$\R \oplus_{\omega_D} \g$. The Lie algebra 
$$ \hat\g := \g(\kappa, D) 
:= (\R \oplus_{\omega_D} \g) \rtimes_{\tilde D} \R$$ 
with the Lie bracket 
$$ [(z,x,t), (z',x',t')] = (\kappa(Dx,x'), [x,x'] + tDx' - t'Dx,0) $$
is called the corresponding {\it double extension}. 
It carries a non-degenerate invariant symmetric bilinear form 
$$ \hat\kappa((z,x,t), (z',x',t)) = \kappa(x,x') + zt'+z't, $$
so that $(\hat\g,\hat\kappa)$ also is a quadratic Fr\'echet--Lie algebra 
(cf.~\cite{MR85}). 
In terms of the hermitian extension of $\kappa$ to $\g_\C$, the Lie bracket on 
$\hat\g_\C$ is given by 
\begin{equation}
  \label{eq:brackcomplex}
[(z,x,t), (z',x',t')] 
= (\kappa(Dx,\oline{x'}), [x,x'] + tDx' - t'Dx,0).
\end{equation}
\end{defn}

\begin{ex} \mlabel{ex:loop} 
Let $(\fk,\la\cdot, \cdot \ra)$ be a Hilbert--Lie algebra 
and let $\ft \subeq \fk$ be a maximal abelian subalgebra, so that 
$\fk$ has a root decomposition with respect to $\ft$ (cf.\ \cite{Sch61}) 
and all roots in $\Delta(\fk, \ft)$ are compact. 
We consider the corresponding loop algebra 
$\cL(\fk)$ of $2\pi$-periodic smooth functions $\R \to \fk$. 
Then 
\[ \la \xi,\eta\ra := \frac{1}{2\pi} 
\int_0^{2\pi} \la \xi(t),\eta(t)\ra\, dt \] 
defines a non-degenerate invariant symmetric 
bilinear form on $\cL(\fk)$. We use the same notation for the unique hermitian extensions 
of $\la \cdot, \cdot \ra$ to $\fk_\C$ and to $\cL(\fk)_\C$. 
Further, $D\xi := \xi'$ is a $\la \cdot, \cdot \ra$-skew symmetric 
derivation on $\cL(\fk)$, so that we may form the associated double extension 
\[  \g := \hat\cL(\fk) := 
(\R \oplus_{\omega_D} \cL(\fk)) \rtimes_{\tilde D} \R,\] 
where $\omega_D(\xi,\eta) = \la \xi',\eta\ra$ and 
$\tilde D(z,\xi) := (0,\xi')$ is the canonical extension of $D$ to 
the central extension $\R \oplus_{\omega_D} \cL(\fk)$ 
(cf.\ Definition~\ref{def:doubext}). This Lie algebra is called the 
{\it affinization of the quadratic Lie algebra 
$(\fk, \la \cdot, \cdot \ra)$.} 
Now 
\[  \kappa((z_1,\xi_1,t_1),(z_2,\xi_2,t_2)) := z_1 t_2 + z_2 t_1 + 
\la \xi_1, \xi_2\ra \] 
is a continuous invariant Lorentzian symmetric bilinear form on $\g$ and 
$\ft_\g := \R \oplus \ft \oplus \R$
is maximal abelian and elliptic in $\g$ (cf.\ Definition~\ref{def:basic2}(b)). 
In the following we identify $\ft$ with the subspace $\{0\} \times \ft \times \{0\}$ 
of $\ft_\g$. We put 
\begin{equation}
  \label{eq:cd}
c := (i,0,0)\in i \g \subeq \g_\C, \quad 
d := (0,0,-i)  \quad \mbox{ and } \quad 
e_n(t) := e^{int}. 
\end{equation}
Then  $c$ is central and the eigenvalue of $\ad d$ on $e_n \otimes \fk_\C$ is $n$. 
It is now easy to verify that 
the set $\Delta_\g$ of roots of $(\g,\ft_\g)$ 
can be identified with the set 
$$ (\Delta(\fk,\ft) \times \Z) \cup (\{0\} \times (\Z \setminus \{0\})) \subeq 
i\ft' \times \R, \ \ \mbox{ where } \ \ 
(\alpha,n)(z,h,t) := (0,\alpha,n)(z,h,t) = \alpha(h) + i t n. $$
The roots $(0,n)$, $0\not=n \in \Z$, corresponding to the 
root spaces $e_n \otimes \ft_\C$, $n \not=0$, are of nilpotent type 
and $(\Delta_\g)_c = \Delta(\fk,\ft) \times \Z$ is the set of compact roots.

For $\alpha \in \Delta(\fk, \ft)$ pick $x_\alpha \in \fk_\C^\alpha$ with 
$[x_\alpha, x_\alpha^*] = \check \alpha$, so that we obtain with equations 
\eqref{eq:8} and \eqref{eq:intflat} the relation 
\begin{equation}
  \label{eq:cor}
\la x_\alpha, x_\alpha \ra 
= \frac{\la \check \alpha, \check \alpha\ra}{2} = \frac{2}{(\alpha, \alpha)}, 
\end{equation}
which leads with  \eqref{eq:brackcomplex} to 
\[ \la D(e_n \otimes x_\alpha), \oline{e_{-n} \otimes x_\alpha^*}\ra 
=   \la in e_n \otimes x_\alpha, - e_{n} \otimes x_\alpha\ra 
=  -in \la x_\alpha, x_\alpha \ra 
= - \frac{2in}{(\alpha, \alpha)}. \]  
For the corresponding root vectors 
$x_{(\alpha,n)} = e_n \otimes x_\alpha \in \g_\C^{(\alpha, n)}$,  
we thus obtain with \eqref{eq:brackcomplex} 
\[  
[e_n \otimes x_\alpha, (e_{n} \otimes x_\alpha)^*] 
= [e_n \otimes x_\alpha, e_{-n} \otimes x_\alpha^*] 
= \Big(-\frac{2in}{(\alpha,\alpha)}, \check \alpha,0\Big). \]
Since, by definition, 
$(\alpha,n)$ takes the value $2$ 
on this element, it follows that 
\begin{equation}
  \label{eq:coroot}
(\alpha,n)\check{} 
= \Big(-\frac{2in}{(\alpha,\alpha)}, \check \alpha\Big)
= \check \alpha 
- \frac{2n}{(\alpha,\alpha)} c 
= \check \alpha 
- \frac{n\|\check \alpha\|^2}{2} c.
\end{equation}

For a linear functional $\lambda \in i \ft_\g' \cong 
i\R \times i\ft'\times i\R$ and 
$\lambda_c := \lambda(c)$, we conclude that $\lambda$ is an integral weight  
if and only if 
\[ \lambda\big((\alpha,n)\check{}\big) = 
 \lambda(\check \alpha) 
- \frac{2n}{(\alpha,\alpha)} \lambda_c \in \Z 
\quad \mbox{ for } \quad 0 \not= n \in \Z, \alpha \in \Delta(\fk,\ft)\] 
(cf.\ Definition~\ref{def:1.8}(c)). 
This means that 
\begin{equation}
  \label{eq:gweight} 
\cP(\g,\ft_\g) 
= \Big\{ \lambda \in i\ft_\g' \: 
\lambda\res_{\ft}  \in \cP(\fk,\ft),\ 
(\forall \alpha \in \Delta(\fk,\ft))\ \lambda_c \in \frac{(\alpha,\alpha)}{2} 
\Z\Big\}.
\end{equation}
\end{ex}

\begin{ex} \mlabel{ex:twistloop} 
In addition to the setting of the preceding example, we assume that 
$\fk$ is semisimple and let $\phi \in \Aut(\fk)$ be an automorphism of order $N$. 
Then 
\[ \cL_\phi(\fk) 
:= \Big\{ \xi \in C^\infty(\R,\fk) \: (\forall t \in \R)\ 
\xi\Big(t + \frac{2\pi}{N}\Big) = \phi^{-1}(\xi(t))\Big\} \] 
is a closed Lie subalgebra of $\cL(\fk)$. 
Accordingly, we obtain a Lie subalgebra 
\[  \g := \hat\cL_\phi(\fk) := 
(\R \oplus_{\omega_D} \cL_\phi(\fk)) \rtimes_{\tilde D} \R 
\subeq \hat\cL(\fk),\] 
called the {\it $\phi$-twisted affinization of 
$(\fk, \la \cdot, \cdot \ra)$.} 

Let  $\ft \subeq \fk^\phi$ be a maximal abelian 
subalgebra, so that $\ft_\fk := \fz_\fk(\ft)$ is maximal abelian in $\fk$ 
by Lemma~\ref{lem:d.2}. Then $\ft_\g = \R \oplus \ft \oplus \R$ 
is maximal abelian in $\cL_\phi(\fk)$ and 
$\Delta_\g := \Delta(\g,\ft_\g)$ can be identified with the set of pairs 
$(\alpha,n)$, where 
\[ (\alpha,n)(z,h,t) := (0,\alpha,n)(z,h,t) = \alpha(h) + i t n, \quad 
n \in \Z, \alpha \in \Delta_n, \] 
where $\Delta_n \subeq i\ft'$ is the set of 
$\ft$-weights in $\fk_\C^n  = \{ x \in \fk_\C \: \phi^{-1}(x) = e^{2\pi in/N} x\}$. 
For $(\alpha,n) \not=(0,0)$, the corresponding root space is 
\[  \g_\C^{(\alpha,n)} = e_n \otimes \fk_\C^{(\alpha,n)} = e_n \otimes (\fk_\C^\alpha \cap \fk_\C^{n}), 
\quad \mbox{ where } \quad 
e_n(t) = e^{i n t}. \] 
The discussion in Appendix~\ref{app:d} implies that 
\[ (\Delta_\g)_c = \{ (\alpha,n) \: 0 \not= \alpha \in \Delta_n, n \in \Z \}. \]
This leads to the $N$-fold layer structure  
\[ (\Delta_\g)_c  = \bigcup_{n = 0}^{N-1} \Delta_n^\times \times (n+ N \Z),  
\quad \mbox{ where } \quad 
\Delta_n^\times := \Delta_n \setminus \{0\}.\] 
For $n \in \Z$ and $x \in \fk_\C^{(\alpha,n)}$ with 
$[x,x^*] = \check \alpha$  (cf.\ Appendix~\ref{app:d}), the element 
$e_n \otimes x \in \g_\C^{(\alpha,n)}$ satisfies 
$(e_n \otimes x)^* = e_{-n} \otimes x^*$, 
which leads to the coroot 
\[ [e_n \otimes x, (e_n \otimes x)^*] = (\alpha,n)\check{} 
= \Big(-\frac{2in}{(\alpha,\alpha)}, \check \alpha,0\Big)
= \check \alpha - \frac{2n}{(\alpha,\alpha)} c. \] 
As $(\alpha,n) \in \Delta_\g$ implies 
$(\alpha,n + k N) \in \Delta_\g$ for every $k \in \Z$, 
we obtain the following description of the integral weights
\begin{equation}
  \label{eq:gweighttwist}
\cP(\g,\ft_\g) 
= \Big\{ \lambda \in i\ft_\g' \: 
(\forall \alpha \in \Delta_n^\times, n \in \Z)\ 
\lambda_c \in \frac{(\alpha,\alpha)}{2N} \Z, 
\lambda(\check \alpha) 
\in \Z + \frac{2n}{(\alpha,\alpha)}\lambda_c\Big\}.
\end{equation}

Suppose that $\fk$ is simple. 
Then Lemma~\ref{lem:d.3} implies that $(\Delta_\g)_c$ does not decompose
into two mutually orthogonal proper subsets, so that 
\[  \hat\cL_\phi(\fk)^{\rm alg}_\C := 
\C c + \C d + \Spann_\C\check\Delta_\g + \sum_{(\alpha,n)\in \Delta_\g} \g_\C^{(\alpha,n)} \] 
is a locally extended affine Lie algebra in the sense of \cite[Def.~1.2]{Ne10a} 
(see also \cite{MY06}). 
Since only the roots of the form $(0,n)$, $n \in \Z$, are isotropic, 
\cite[Prop.~2.5]{Ne10a} further implies that $(\Delta_\g)_c$ is a locally affine 
root system. 
\end{ex}

\begin{defn} \mlabel{def:locaffroot} 
Instead of going into the axiomatics of locally affine root systems 
developed in \cite{YY10}, 
we only recall that a locally affine root system is in particular 
a subset $\Delta$ of a 
vector space $V$ endowed with a positive semidefinite form. 
For a reduced locally finite root system $\Delta$ in the euclidean space $V$ 
(such as $A_J$, $B_J$, $C_J$ or $D_J$ from Examples~\ref{ex:d.1a}, \ref{ex:d.1b}, \ref{ex:d.1c}),  
we put $\Delta^{(1)} 
:= \Delta \times \Z\subeq V \times \R$, where the scalar product on 
$V \times \R$ is defined 
by $((\alpha,t), (\alpha', t')) := (\alpha,\alpha')$. 
According to Yoshii's classification (\cite[Cor.~13]{YY10}), 
there exist $7$ isomorphy classes of irreducible reduced 
locally affine root system of infinite rank: the four untwisted 
reduced root  systems $A_{J}^{(1)}, B_J^{(1)}, C_J^{(1)}, D_J^{(1)}$, and, for 
$BC_J := B_J \cup C_J$,  
the three twisted root systems 
\begin{align*}
B_J^{(2)} &:= \big(B_J\times 2\Z\big) 
\cup \big(\{ \pm \eps_j \: j \in J\}\times (2\Z + 1)\big), \\
C_J^{(2)} &:= (C_J \times 2\Z) 
\cup \big(D_J  \times (2 \Z+1)\big) \\
(BC_J)^{(2)} &:= (B_J \times 2\Z) \cup \big(BC_J\times (2\Z+1)\big). 
\end{align*}
\end{defn}

\begin{rem}  \mlabel{rem:weylgrp} 
Let $q \: V \times \R \to V, (v,t) \mapsto v$ denote the projection. Then the root systems 
$\Delta^{(2)}$ satisfy $q(\Delta^{(2)}) = \Delta$, but 
$\Delta_0 := \{ \alpha \in \Delta \: (\alpha,0) \in \Delta^{(2)}\}$ 
may be smaller, as the example $BC_J$ shows. However, 
in all $7$ cases the subsystem $\Delta_0 \subeq \Delta$ 
contains enough elements to obtain all generating 
reflections of the Weyl group of $\Delta$. Hence 
the subgroup $\cW_0 := \la r_{(\alpha,0)} \: \alpha \in \Delta_0\ra \subeq \hat\cW$ 
is isomorphic to the Weyl group $\cW$ of $\Delta$. 
\end{rem}

\subsection{Three involutive automorphisms} 

In this subsection we introduce three involutive automorphisms whose significance lies 
in the fact that these automorphism lead to the three twisted affine 
root systems of infinite rank $B_J^{(2)}$, $C_J^{(2)}$ and $BC_J^{(2)}$. 
In the following we call these three automorphism {\it standard}.

\begin{ex}
  \mlabel{ex:auto-bj2}
Consider $\fk = \fo_2(\cH)$, where $\cH$ is a real infinite-dimensional 
Hilbert space. We consider the automorphism 
$\phi(x) := g x g^{-1}$, where $g$ is the orthogonal reflection in 
the hyperplane $v_0^\bot$ for some unit vector~$v_0 \in \cH$. 

Then $\fk^\phi \cong \fo_2(v_0^\bot) \cong \fo_2(\cH)$ is a simple 
Hilbert--Lie algebra. Pick 
\[ \ft \subeq \fk^\phi  = \{ x \in \fk \: x v_0 = 0\} \] 
maximal abelian with $\dim(\ker(\ft) \cap v_0^\bot) = 1$. 
Then $\ft_\fk = \fz_\fk(\ft)$ is maximal abelian (Lemma~\ref{lem:d.2}) with 
$\ker(\ft_\fk) = \{0\}$ and $\ft$ is a hyperplane in $\ft_\fk$. Hence the root system  
$\Delta(\fk,\ft_\fk)$ is of type $D_J$, and 
$\Delta(\fk^\phi,\ft)$ is of type $B_{J_0}$, where 
$J_0 = J\setminus \{j_0\}$ for some $j_0 \in J$ (cf.\ Example~\ref{ex:d.1c}). 
It follows from \cite[Thm.~5.7(ii)]{Ne10a} that the set of compact roots 
of $\hat\cL_\phi(\fk)$ is of type $B_J^{(2)}$ 
(cf.\ Example~\ref{ex:twistloop}). 
\end{ex}

\begin{ex}  \mlabel{ex:auto-cj2} 
To obtain a root system  of type 
$C_J^{(2)}$, we start with $\fk = \fu_2(\cH)$ 
for a complex Hilbert space $\cH$. 
We write $\cH$ as $\cH_0 \oplus  \cH_0$ for a complex 
Hilbert space $\cH_0$ endowed with a conjugation~$\sigma_0$ 
and extend this conjugation by 
$\sigma(x,y) := (\sigma_0(x), \sigma_0(y))$ to a conjugation of $\cH$. 
Then 
\[ \phi(g) = S (g^\top)^{-1} S^{-1} = S \sigma g \sigma  S^{-1} \quad \mbox{ for } \quad 
S = \pmat{ 0 & \1 \\ -\1 & 0}\] 
defines an involutive automorphism of $\U_2(\cH)$. 
Let $(e_j)_{j \in J}$ be an ONB of $\cH_0$, so that 
$\{e_j, Se_j\: j \in J\}$ is an ONB of $\cH$. 
Let $\ft \subeq \fk^\phi$ be those elements which are diagonal 
with respect to this ONB. Then $\ft_\fk = \fz_\fk(\ft)$ consists of all elements in 
$\fu_2(\cH)$ which are diagonal with respect to this 
ONB, and from \cite[Thm.~5.7(ii)]{Ne10a} we know that the set of compact roots 
of $\hat\cL_\phi(\fk)$ is of type $C_J^{(2)}$.
As $\tilde\sigma := S \sigma = \sigma S$ is an anticonjugation of $\cH$, it defines a quaternionic 
structure $\cH_\H = (\cH,\tilde\sigma)$. In this sense we have 
\[ U_2(\cH)^\phi = \{ g \in \U_2(\cH) \: \tilde\sigma g \tilde\sigma^{-1} = g \} 
\cong \Sp_2(\cH_\H).\] 
\end{ex}

\begin{ex}  \mlabel{ex:auto-bcj2} 
To obtain a root system of type 
$BC_J^{(2)}$, we consider $\fk = \fu_2(\cH)$ for a complex Hilbert space 
$\cH$. We write $\cH$ as $\cH_0 \oplus \C \oplus  \cH_0$ for a complex 
Hilbert space $\cH_0$ endowed with a conjugation~$\sigma_0$ 
and extend $\sigma_0$ by  $\sigma(x,y,z) := (\sigma_0(x),\oline y, \sigma_0(z))$ 
to a conjugation of~$\cH$. We consider the automorphism 
\[ \phi(g) = S (g^{-1})^\top S^{-1} \quad \mbox{ for } \quad 
S = \pmat{ 0 & 0 & \1 \\ 0 & 1 & 0 \\ \1 & 0 & 0}, g \in \U_2(\cH).\] 
As $S$ commutes with $\sigma$, $\tau := S \sigma$ is a conjugation 
on $\cH$ with $\phi(g) = \tau g \tau$, so that 
\[ \U_2(\cH)^\phi 
= \{ x \in \U_2(\cH) \:  S (g^{-1})^\top S^{-1}  = g \} 
= \{ x \in \U_2(\cH) \:  \tau g \tau = g \} 
\cong \OO_2(\cH^\tau).\] 
In particular, $\fk^\phi\cong \fo_2(\cH)^\tau$ is a simple Hilbert--Lie algebra. 

Pick an ONB $(e_j)_{j \in J}$ of $\cH_0$, so that the elements 
$(e_j,0,0), (0,0,e_j)$, $j \in J$, together 
with $(0,1,0)$, form an ONB of $\cH$. 
Let $\ft \subeq \fk^\phi$ be those elements which are diagonal 
with respect to this ONB. Then $\ft_\fk = \fz_\fk(\ft)$ consists of all elements in 
$\fu_2(\cH)$ which are diagonal with respect to this 
ONB, and from \cite[Thm.~5.7(ii)]{Ne10a} we know that the set of compact roots 
of $\hat\cL_\phi(\fk)$ is of type $BC_J^{(2)}$. 
\end{ex}

\begin{rem} In Examples~\ref{ex:auto-cj2} and 
\ref{ex:auto-bcj2} the operator $S$ is contained in the 
identity component of $\U(\cH)$ (Theorem~\ref{thm:aut-grp}), so that in both cases, the 
involution $\phi$ is homotopic to the automorphism 
$x \mapsto - x^\top = \sigma x \sigma$ of~$\fu_2(\cH)$. 
\end{rem}

We now record some topological properties of the subgroup $K^\phi$ for the standard involutions. 

\begin{lem} \mlabel{lem:kphicon} 
If $K$ is simply connected, then $K^\phi$ is connected 
for each of the standard involutions.   
It is $1$-connected in all cases except the $BC_J^{(2)}$-case, in which 
$\pi_1(K^\phi) \cong \Z/2$. 
\end{lem}

\begin{prf} (a) First we consider $K := \OO_2(\cH)_0$ and 
$\phi(k) = gkg^{-1}$, where $g \in \OO(\cH)$ is the reflection 
in a hyperplane $v_0^\bot$ (Example~\ref{ex:auto-bj2}). We show that the group 
$\tilde K^{\tilde \phi}$ of fixed points of the lifted automorphism 
$\tilde\phi$ of the $2$-fold covering group $q_K \: \tilde K \to K$ 
is connected. 

Let $\OO_2(\cH)_\pm$ denote the two connected components of 
$\OO_2(\cH)$, where $\1 \in \OO_2(\cH)_+$ (Theorem~\ref{thm:1.5}). 
The group  $\OO_2(\cH)^\phi \cong \{ \pm 1 \} \times \OO_2(v_0^\bot)$ 
has four connected components, two of which lie in $\OO_2(\cH)_0$, so that 
\[ K^\phi \cong  (\{1\} \times \OO_2(v_0^\bot)_+) \times 
(\{-1\} \times \OO_2(v_0^\bot)_-) \] 
has two connected components and 
$(K^\phi)_0 \cong \OO_2(v_0^\bot)_+$. 
We conclude that $\pi_1(K^\phi) \cong \Z/2$ (Theorem~\ref{thm:1.5}) 
and that the inclusion 
$j \: (K^\phi)_0 \to K$ induces an isomorphism of fundamental groups, 
so that $\tilde K^{\tilde\phi}$ has a simply connected identity component which is the 
universal covering group of $K^\phi_0$ and contains 
$\ker(q_K)$. In particular,  $q_K^{-1}(K^\phi_0)$ is connected. 
To see that $\tilde K^{\tilde \phi}$ is connected, it therefore 
suffices to show that $q_K$ maps it into $K^\phi_0$, which is 
equivalent to $(-1,r) \not\in q_K(\tilde K^{\tilde \phi})$ 
for any reflection $r \: v_0^\bot \to v_0^\bot$ in a hyperplane 
$w^\bot$ of $v_0^\bot$. Let $V := \Spann_\R \{ v_0,w \} \cong \R^2$. Then 
the inclusion $j \:\T \cong \SO(V)\to \OO_2(\cH)$ induces a surjection 
$\Z \cong \pi_1(\T) \to \pi_1(\OO_2(\cH))$ (\cite{Ne02a}), so that it 
lifts to an inclusion 
$\tilde j \: \Spin(V) \to \tilde K$, where 
$\Spin(V)$ denotes the unique $2$-fold covering of $\SO(V)$. 
In $\Spin(V)$ the inverse image of $-\id_V$ consists of 
two elements of order $4$ which are exchanged by 
$\tilde\phi$. Therefore $(-\1,r)$ does not lift to a 
$\tilde\phi$-fixed element. This proves the connectedness of 
$\tilde K^{\tilde\phi}$. 

(b) For $K = \U_2(\cH)$ and  
$\phi(g) = \tilde\sigma g \tilde\sigma^{-1}$ as in Example~\ref{ex:auto-cj2}, the subgroup  
$K^\phi \cong \Sp(\cH_\bH)$ is $1$-connected by Theorem~\ref{thm:1.5}. 
Since $(\tilde K^\phi)_0$ is a covering of $K^\phi$, this group 
is also simply connected and isomorphic to $K^\phi$. 
For the universal covering $q_K \: \tilde K \to K$ we thus obtain 
\[ q_K^{-1}(K^\phi) \cong \ker q_K \times K^\phi 
\cong \Z \times K^\phi,\] 
and $\tilde\phi$ acts by $-\id$ on $\pi_1(K)$ (Proposition~\ref{prop:1.18}). 
Therefore $\tilde K^{\tilde\phi} \cong K^\phi$ is connected and hence 
$1$-connected. 

(c) For $K = \U_2(\cH)$ and  
$\phi(g) = \tau g \tau$ as in 
Example~\ref{ex:auto-bcj2}, the group 
$K^\phi \cong \OO_2(\cH^\tau)$ has $2$ connected components  and 
satisfies $\pi_1(K^\phi) \cong \Z/2$ (Theorem~\ref{thm:1.5}). 
Since the inclusions $\SO_n(\R) \to \U_n(\C)$ induces the trivial 
homomorphism $\pi_1(\SO_n(\R)) \cong \Z/2\to \pi_1(\U_n(\C))\cong \Z$ for each $n>2$, 
the same holds for the inclusion 
$(K^\phi)_0 \into K$ (cf.\ \cite{Ne02a}). In particular, this inclusion 
lifts to a homomorphism $(K^\phi)_0 \to \tilde K$, so that 
$\pi_1(\tilde K^{\tilde \phi}) \cong \Z/2$. 

To see that $\tilde K^{\tilde\phi}$ is connected, it suffices to show 
that, for any reflection $r \in \OO_2(\cH^\tau)$ in a hyperplane 
$v_0^\bot \subeq \cH^\tau$, there is 
no $\tilde\phi$-invariant inverse image under the covering map 
\break $q_K \:  \tilde K \to K$. Since the homotopy groups of 
$\U_2(\cH)$ are obtained from the direct limit 
$\indlim \U_n(\C)$ (cf.\ \cite{Ne02a}), it suffices to prove the 
corresponding assertion for the case where $n := \dim \cH < \infty$. 
Then $\tilde \U_n(\C) \cong \SU_n(\C) \rtimes \R$ with 
$\tilde\phi(g,t) = (\oline g, -t)$, so that 
$\tilde\U_n(\C)^{\tilde\phi} \cong \SU_n(\C)^\phi = \SO_n(\R)$ 
is connected. This implies that $\tilde K^{\tilde\phi}$ 
is connected in the general case. 
\end{prf}

\subsection{The adjoint action of twisted loop groups} 

We claim that the following 
formula describes the adjoint action of $\cL_\phi(K)$ on $\g = \hat\cL_\phi(\fk)$: 
\begin{equation}
  \label{eq:adjt}
\Ad_{\g}(g)(z,\xi,t) = 
\Big(z + \la \delta^l(g), \xi\ra 
- \frac{t}{2}\|\delta^r(g)\|^2, \Ad(g)\xi - t \delta^r(g), t\Big)   
\end{equation}
where $\delta^r(g) = g' g^{-1} \in \cL_\phi(\fk)$ 
denotes the {\it right logarithmic derivative} 
and $\delta^l(g) = g^{-1}g' \in \cL_\phi(\fk)$ 
denotes the {\it left logarithmic derivative} 
(cf.\ \cite[Prop.~4.9.4]{PS86} for the case where $K$ is compact and a 
different choice of sign on the center). 
In fact, modulo the center, the formula 
$g.(\xi, t) = (\Ad(g)\xi - t \delta^r(g), t\Big)$ 
defines a smooth action of $\cL_\phi(K)$ on $\cL_\phi(\fk) \rtimes \R$, 
whose derived action is given by 
\[ \eta.(\xi,t) = ([\eta, \xi] - t D \eta, 0) = [(\eta,0), (\xi,t)]\] 
which implies that it is the adjoint action of the group 
$\cL_\phi(K)$. The formula for the central 
component is now obtained from the invariance of 
the Lorentzian form on $\hat\cL_\phi(\fk)$ under the adjoint action. 

\begin{rem} Since we need it below, we take a closer look 
at the affine action of $\cL_\phi(K)$ on $\cL_\phi(\fk)$ by 
\[ g * \xi := 
\Ad(g)\xi - \delta^r(g) = \Ad(g)\xi - g' g^{-1}.\] 
To understand its orbits, for a smooth curve 
$\xi \: \R \to \fk$, let 
$\gamma_\xi \: \R \to K$ be the unique smooth curve 
with $\gamma_\xi(0) = \1$ and 
$\delta^l(\gamma_\xi) = \xi$. 
We consider the smooth holonomy maps 
\[ \Hol_t \: \cL_\phi(\fk) \to K, \quad 
\Hol_t(\xi) := \gamma_\xi(t), \quad t \in \R.\] 
If $\xi \in \cL_\phi(\fk)$, then 
\[ \delta^l(\gamma_\xi)_{t + 2\pi/N} 
= \xi\Big(t + \frac{2\pi}{N}\Big)
= \L(\phi)^{-1}(\xi(t)) 
= \delta^l(\phi^{-1} \circ \gamma_\xi)_t\] 
implies that 
\begin{equation}
  \label{eq:translate}
\gamma_{\xi}\Big(t + \frac{2\pi}{N}\Big)
= \Hol_{2\pi/N}(\xi)\phi^{-1}(\gamma_\xi(t)).
\end{equation}
From 
\begin{equation} \label{eq:quotrule}
\delta^l(\gamma_{\xi} g^{-1}) 
= \delta^l(g^{-1}) + \Ad(g) \xi 
= g * \xi, 
\end{equation}
we derive the following equivariance property 
\[\Hol_s(g * \xi) 
= g(0) \Hol_s(\xi)g(s)^{-1}.\] 
For $s = 2\pi/N$, we obtain in particular 
\begin{equation}\label{eq:evaleq}
\Hol_{2\pi/N}(g * \xi) = g(0) \Hol_{2\pi/N}(\xi) \phi^{-1}(g(0)^{-1}).
\end{equation}
\end{rem}

\begin{prop} \mlabel{prop:2.x} 
The map $\Hol_{2\pi/N} \: \cL_\phi(\fk) \to K$ 
is equivariant with respect to the action of 
$\cL_\phi(K)$ on $K$ for which 
$g$ acts by $c_{g(0)}^\phi$, where 
$c_k^\phi(k') := k k' \phi^{-1}(k^{-1})$ is the $\phi^{-1}$-twisted conjugation map. 
The fibers of $\Hol_{2\pi/N}$ coincide with the 
orbits of the subgroup
\[ \cL_\phi(K)_* := \{ g \in \cL_\phi(K)\: g(0) = \1\} \] 
and the image of the $\cL_\phi(K)$-orbits are the 
$\phi^{-1}$-twisted conjugacy classes in $K$.  
\end{prop}

\begin{prf} The asserted equivariance is \eqref{eq:evaleq}. 
In particular, $\Hol_{2\pi/N}$ is constant 
on the orbits of the subgroup $\cL_\phi(K)_*$. 
If $\Hol_{2\pi/N}(\xi_1) = \Hol_{2\pi/N}(\xi_2)$, then \eqref{eq:translate} implies that the smooth curve 
$g:= \gamma_{\xi_2}^{-1}\gamma_{\xi_1} \: \R \to K$ 
is contained in $\cL_\phi(K)_*$. It satisfies 
$\gamma_{\xi_2} = \gamma_{\xi_1}g^{-1}$, so that 
\eqref{eq:quotrule} implies $g * \xi_1 = \xi_2$. 
This completes the proof. 
\end{prf}

\section{Double extensions of twisted loop groups} \mlabel{sec:3}

In the preceding section we have introduced the double extension 
$\hat\cL_\phi(\fk)$ of the Fr\'echet--Lie algebra 
$\cL_\phi(\fk)$ of smooth $\phi$-twisted loops. 
Since its 
construction involves a central extension, it is not obvious 
that this extension is the Lie algebra of a Lie group.
In this section we show that such a group always exists 
for a simple infinite-dimensional 
Hilbert--Lie algebra $\fk$, in particular, we obtain a corresponding 
$1$-connected Lie group which we denote $\hat\cL_\phi(K)$.

\subsection{The central extension for smooth loops}

\begin{rem} \mlabel{rem:homot} Let $K$ be a connected Lie group for which 
$\fk$ is a simple Hilbert--Lie algebra and write 
$\ev_0 \: \cL_\phi(K) \to K, f \mapsto f(0)$
for the evaluation map. 
In \cite[Prop.~3.5, Rem.~3.6]{NeWo09}, we have seen that the 
long exact homotopy sequence of the Lie group extension 
$$ \1 \to \cL_\phi(K)_* := \ker(\ev_0) 
\to \cL_\phi(K) \smapright{\ev_0} K \to \1 $$
provides crucial information on the homotopy groups of 
$\cL_\phi(K)$. 
For each $j \geq 1$, we obtain a short exact sequence 
$$\1\to \pi_j(K)_\phi := \pi_j(K)/\im(\pi_j(\phi) -\id) 
\into \pi_{j-1}(\cL_\phi(K)) \onto 
\pi_{j-1}(K)^\phi \to\1. $$
As $K$ is connected, $\pi_2(K)$ vanishes (Theorem~\ref{thm:1.5}) and 
$\pi_3(\phi) = \id$ (Corollary~\ref{cor:2.10}), we obtain in particular  
$$ \pi_0(\cL_\phi(K)) \cong \pi_1(K)_\phi, \quad 
\pi_1(\cL_\phi(K)) \cong \pi_1(K)^\phi \quad \mbox{ and } \quad 
\pi_2(\cL_\phi(K)) \cong \pi_3(K) \cong \Z. $$
If $K$ is $1$-connected, these relations imply that 
$\cL_\phi(K)$ is also $1$-connected. 
\end{rem}

\begin{defn} \mlabel{def:circ} 
In the following we shall identify the Lie algebra 
$\L(\T)$ of the circle group $\T \subeq \C^\times$ 
with $i\R$, so that the exponential function 
is given by $\exp_\T(it) = e^{it}$ with 
$\ker(\exp_\T) = 2\pi i \Z$. 
\end{defn}

The following theorem generalizes the corresponding result for 
compact target groups which can be found in \cite[Sect.~4.4]{PS86} 
for the untwisted case. 

\begin{defn} We say that the scalar product $\la \cdot,\cdot \ra$ on 
the simple Hilbert--Lie algebra $\fk$ is {\it normalized} 
if $\|\check \alpha\|^2 = 2$ holds for the coroots 
of all long roots\begin{footnote}{At most three root lengths occur, 
the long roots are those of maximal lenght.} \end{footnote}
$\alpha \in \Delta(\fk,\ft_\fk)$, where $\ft_\fk \subeq \fk$ is maximal abelian.  
In view of \eqref{eq:intflat}, this is equivalent to 
$(\alpha, \alpha) = 2$ for all long roots. 
\end{defn}

\begin{thm}\mlabel{thm:cent-ext} Let $K$ be a $1$-connected 
simple Hilbert--Lie group and suppose that the scalar product on $\fk$ 
is normalized. Then the central Lie algebra extension 
$\tilde\cL_\phi(\fk)$ defined by the cocycle 
\[ \omega(\xi, \eta) = \frac{1}{2\pi} \int_0^{2\pi} \la \xi'(t), \eta(t)\ra\, dt \] 
integrates to a $2$-connected central Lie group extension 
  $$ \1 \to Z \to \tilde\cL_\phi(K) \to \cL_\phi(K) \to \1, $$
where $Z = \exp(\R ic)$ and $\ker(\exp\res_{\R i c}) = 2\pi N \Z i c$. 
\end{thm} 

\begin{prf} Let $\omega^l \in \Omega^2(\cL_\phi(K),\R)$ denote the 
left invariant $2$-form corresponding to the $2$-cocycle 
$\omega$. 
We derive from Remark~\ref{rem:homot} that $\cL_\phi(K)$ 
is $1$-connected because $K$ is $1$-connected. 
Therefore we can use \cite[Thm.~7.9]{Ne02b} to see that 
the required simply connected Lie group extension 
$\tilde\cL_\phi(K)$ exists if and only if 
the range of the period homomorphism 
$$ \per_\omega \: \pi_2(\cL_\phi(K)) \to \R, \quad 
[\sigma] \mapsto \int_{\bS^2} \sigma^*\omega^l $$
is discrete and in this case $Z \cong \R/\im(\per_\omega)$. 

To use the results from \cite{NeWo09}, where the circle is identified
with $\R/\Z$, we have to translate them to our context where 
$\bS^1 \cong \R/\frac{2\pi}{N} \Z$. 
So let 
\[  \cL_\phi^1(K) 
:= \big\{ f \in C^\infty(\R,K) \: (\forall t \in \R)\ 
f(t + 1)= \phi^{-1}(f(t))\big\} \] 
and observe that 
\[ \Phi \: \cL_\phi(K) \to \cL_\phi^1(K), \quad \Phi(\xi)(t) := \xi(2\pi t/N) \] 
is an isomorphism of Fr\'echet--Lie groups. For the cocycle 
$\omega^1(\xi,\eta) :=  \int_0^1 \la \xi'(t), \eta(t)\ra\, dt$ 
on $\cL_\phi^1(\fk)$ we then obtain 
\begin{align*}
\big(\L(\Phi)^*\omega^1\big)(\xi,\eta) 
&= \int_0^1 \la \xi'(2\pi t/N)\frac{2\pi}{N}, \eta(2\pi t/N)\ra\, dt
= \int_0^{2\pi/N} \la \xi'(t),\eta(t)\ra\, dt
= \frac{2\pi}{N}\omega(\xi,\eta).
\end{align*}
This implies that $\im(\per_\omega) = \frac{N}{2\pi} \im(\per_{\omega^1})$. 

According to \cite[Lemma~3.10]{NeWo09}, the 
period homomorphism $\per_{\omega^1}$ of the restriction of $\omega^1$ to the 
ideal $\cL_{\phi}^1(\fk)_*$ coincides with the homomorphism 
$$ \frac{1}{2} \per_{C_\fk} \: 
\pi_2(\cL_\phi^1(K)_*) \cong \pi_3(K) \cong \Z\to \R, $$
where $C_\fk(x,y,z) := \la [x,y],z\ra$ 
is the $3$-cocycle defined by the scalar product on $\fk$ and 
$$ \per_{C_\fk} \: \pi_3(K) \to \R, \quad 
[\sigma] \mapsto \int_{\bS^3} \sigma^*C_\fk^l $$
is the period homomorphism of the corresponding closed 
$3$-form $C_\fk^l$ on $K$. 
This map can be evaluated quite explicitly as follows.  
Let $\alpha \in \Delta(\fk,\ft_\fk)$ be a long root and 
$\fk(\alpha) \subeq \fk$ be the corresponding $\su_2$-subalgebra and 
$\check \alpha \in i \ft$ be the corresponding coroot. 
The associated homomorphism $\gamma_\alpha \: \SU_2(\C) \to K$ 
induces an isometric embedding $\L(\gamma_\alpha) \: \su_2(\C) \to \fk$ 
with respect to the normalized scalar products. 
Hence \cite[Ex.~3.11]{NeWo09} implies that 
\[ \frac{1}{2} \per_{C_\fk}([\gamma_\alpha]) 
= \frac{1}{2} \frac{\|\check \alpha\|^2}{2} 8 \pi^2 
=  4\pi^2. \]
Since $\pi_2(K)$ vanishes (Remark~\ref{rem:homot}), 
\cite[Thm.~3.12]{NeWo09} thus leads to 
\[ \im(\per_{\omega}) = \frac{N}{2\pi} \im(\per_{\omega^1}) 
=   2\pi N \Z.\] 
As $\pi_2(\cL_\phi(K)_*) \cong \Z$, the non-zero homomorphism 
$\per_\omega$ is injective. 
Since $\cL_\phi(K)$ is $1$-connected by 
Remark~\ref{rem:homot}, \cite[Rem.~5.12(b)]{Ne02b} now implies that, 
for $Z := \R/\im(\per_\omega)$, the group 
$\tilde\cL_\phi(K)$ is $2$-connected. 
\end{prf}

\begin{defn}\mlabel{def:3.5a}
Since the rotation action $\alpha$ of $\R$ on 
$\cL_\phi(K)$ lifts uniquely to a smooth action on the central 
extension $\tilde\cL_\phi(K)$ (\cite[Thm.~V.9]{MN03}), 
we obtain a $2$-connected Fr\'echet--Lie group 
\[ \hat\cL_\phi(K) := \tilde\cL_\phi(K) \rtimes_\alpha \R. \] 
\end{defn}

\begin{rem} \mlabel{rem:polar} 
If $K$ is a simple Hilbert--Lie group, 
then it has a universal complexification 
$\eta_K \: K \to K_\C$ which has a polar decomposition, i.e., 
the map 
\[ K \times \fk \to K_\C, \quad (k,x) \mapsto k \exp ix \] 
is a diffeomorphism (cf.\ \cite{Ne02a}). This property is inherited by the group 
$\cL_\phi(K)$, which implies in particular that the inclusion 
$\cL_\phi(K) \to \cL_\phi(K_\C)$ induces isomorphisms of all homotopy 
groups.  Hence the cocycle $\omega$ and its complex 
bilinear extension to $\cL_\phi(\fk_\C) \cong \cL_\phi(\fk)_\C$ 
have the same period group. Now \cite[Thm.~7.9]{Ne02b} 
implies the existence of a central extension of complex 
Lie groups 
\[ \1 \to \C^\times \to \tilde\cL_\phi(K_\C) \to \cL_\phi(K_\C) \to \1 \] 
for which the inclusion 
$\tilde\cL_\phi(K) \into \tilde\cL_\phi(K_\C)$ is a universal complexification 
and a weak homotopy equivalence. 
\end{rem}

In the preceding theorem we have seen the importance of normalizing the 
scalar product. To evaluate the period group in all cases, it is thus important 
to identify the normalized scalar products in all cases. 

\begin{rem}(Normalization of the scalar product)  \mlabel{rem:3.4} 

(a) For $\fk = \fu_2(\cH)$ all roots in $\Delta(\fk,\ft) = A_J$
 have the same length and the coroots correspond to diagonal 
matrices of the form $E_{jj} - E_{kk}$, so that 
$\la x,y \ra = \tr(xy^*)$ is a scalar product with 
$\|\check\alpha\|^2 = 2$ for all roots $\alpha$. 

(b) For $\fk = \sp_2(\cH_\H)$, the long roots are of the form 
$\pm 2 \eps_j$ and their coroots are diagonal matrices of the form 
$(E_{jj}, - E_{jj})$ with respect to the decomposition of the complex 
Hilbert space $\cH_\H = \ell^2(J,\C) \oplus \ell^2(J,\C)$. 
Therefore $\la x,y \ra = \tr_\C(xy^*)$ satisfies 
$\|\check\alpha\|^2 = 2$ for all long roots $\alpha$.
 
(c) For $\fk = \fo_2(\cH_\R)$ and $\Delta(\fk,\ft)$ 
of type $B_J$ or $D_J$, the long roots are 
$\pm \eps_j \pm \eps_k$, $j \not=k$. On the complex Hilbert space 
$\cH := (\cH_\R)_\C$ their coroots correspond to diagonal matrices 
of the form $\pm E_{jj} \mp E_{-j,-j} \pm E_{kk} \mp E_{-k,-k}$, 
$j \not=k$ 
(cf.\ Example~\ref{ex:d.1c}) satisfying 
\[ \tr_\C((\pm E_{jj} \mp E_{-j,-j} \pm E_{kk} \mp E_{-k,-k})^2) = 4,\] 
so that  $\la x,y \ra = \shalf\tr_\C(xy^*) = \shalf \tr_\R(xy^\top)$ satisfies 
$\|\check\alpha\|^2 = 2$ for all long roots $\alpha$.

(d) In Example~\ref{ex:auto-bj2} we have the inclusion 
$\fk^\phi \cong \fo_2(v_0^\bot) \into \fk = \fo_2(\cH_\R)$, and 
(c) shows that this is isometric with respect to the normalized scalar products. 

(e) In Example~\ref{ex:auto-cj2} we have for 
$\cH = \cH_0^2 = (\cH_0)_\H$ the inclusion 
$\fk^\phi\cong \sp_2((\cH_0)_\H) \into \fk = \fu_2(\cH)$, so that 
(a) and (b) imply that it is isometric with respect to the normalized scalar products. 

(f) In Example~\ref{ex:auto-bcj2} we have the inclusion 
$\eta \: \fk^\phi\cong \fo_2(\cH^\tau) \into \fk = \fu_2(\cH)$, so that 
(c) implies that $\la \eta(x), \eta(y) \ra = 2 \la x, y \ra$ for 
$x,y \in \fk^\phi$. In particular, the roots 
$\alpha = \pm \eps_j \pm \eps_k \in \Delta(\fk^\phi, \ft) = BC_J$ satisfy 
$\|\check \alpha\|^2_\fk = 4$. Accordingly we find 
$\|\eps_j \pm \alpha_k\|^2 = 1$ and $\|2 \eps_j\|^2 = 2$. 
\end{rem}

At this point we can also make the description 
of the weight set $\cP(\g,\ft)$ from \eqref{eq:gweighttwist} more explicit 
for all $7$ locally affine root systems.  

\begin{ex} \mlabel{ex:affineweights} 
(a) For the untwisted root systems 
$X_J^{(1)}$, we have seen in \eqref{eq:gweight} 
that $\lambda \in i \ft_\g^*$ is an integral weight 
if and only if $\lambda_c \in \frac{\|\alpha\|^2}{2} \Z$ for 
every root $\alpha \in X_J$ and $\lambda\res_{\ft}$ is a weight of $X_J$. 
As $\|\alpha\|^2 = 2$ for the long roots, the condition on 
$\lambda_c$ is equivalent to $\lambda_c \in \Z$. 

(b) For $B_J^{(2)}$ we find the condition 
$\lambda_c \in \frac{\|\alpha\|^2}{4} \Z$ for every root 
$\alpha \in B_J$, i.e., $\lambda_c \in \shalf \Z$, and that 
$\lambda(\check \alpha) \in \Z$ for $\alpha \in \Delta_0 = B_J$. 
As $\Delta_1 = \{ \pm \eps_j \: j \in J\} \subeq \Delta_0$ and 
$\|\eps_j\| = 1$, this implies 
\[ \lambda(\check \alpha) \in \Z + \frac{2}{\|\alpha\|^2} \lambda_c
\quad \mbox{ for } \quad \alpha \in \Delta_1.\] 
Therefore 
\[ \cP(\g,\ft_\g) 
= \big\{ \lambda \in i\ft_\g' \: 
\lambda_c \in \shalf \Z, \ \lambda \in \cP(B_J)\big\}.\] 

(c) For $C_J^{(2)}$, the long roots of $C_J$ 
are of the form $\pm 2 \eps_j$, $j \in J$, 
so that our normalization leads to 
$\|2\eps_j\|^2 = 2$ (cf.\ Remark~\ref{rem:3.4}(e)), which in turn 
implies $\|\eps_j\|^2 = \shalf$. The condition 
$\lambda_c \in \frac{\|\alpha\|^2}{4} \Z$ for every root 
$\alpha$ leads to $\lambda_c \in \shalf\Z$. 
We further obtain 
$\lambda(\check \alpha) \in \Z$ for $\alpha \in \Delta_0 = C_J$ and this 
already implies 
\[ \lambda(\check \alpha) \in \Z + \frac{2}{\|\alpha\|^2} \lambda_c
= \Z + 2 \lambda_c\quad \mbox{ for } \quad \alpha \in \Delta_1 = D_J.\] 
Therefore 
\[ \cP(\g,\ft_\g) 
= \big\{ \lambda \in i\ft_\g' \: 
\lambda_c \in \shalf \Z, \ \lambda \in \cP(C_J)\big\}.\] 

(d) For $BC_J^{(2)}$ we have seen in Remark~\ref{rem:3.4}(f) that 
$\| \eps_j \pm \eps_k\|^2 = 1$ and $\|2\eps_j\|^2 = 2$. 
Hence the condition 
$\lambda_c \in \frac{\|\alpha\|^2}{4} \Z$ for every root 
$\alpha \in BC_J$ means that $\lambda_c \in \shalf\Z$. 
We further obtain 
$\lambda(\check \alpha) \in \Z$ for $\alpha \in \Delta_0 = B_J$ 
which means that $\lambda_j \in \shalf \Z$ for every $j$ with 
$\lambda_j - \lambda_k \in \Z$ for $j \not=k$. 
An integral weight $\lambda$ also has to satisfy 
\[ \lambda(\check \alpha) \in \Z + \frac{2}{\|\alpha\|^2} \lambda_c
\quad \mbox{ for } \quad \alpha \in \Delta_1 = BC_J.\] 
For $\alpha \in B_J$, this is satisfied because
$\|\alpha\|^2 \in \{\shalf,1\}$. For $\alpha = \pm 2 \eps_j$, it means that 
\[ \pm \lambda_j = \lambda(\check \alpha) \in \Z + \lambda_c. \] 
Therefore the parity of  $2\lambda_c$ equals the parity of $2\lambda_j$. 

If we also take into account that $\lambda$ should be continuous, 
i.e., $(\lambda_j) \in \ell^2(J,\R)$, then only finitely many 
$\lambda_j$ are non-zero, which leads to $\lambda_j \in \Z$ 
and hence to $\lambda_c \in \Z$. We therefore have 
\[ \cP(\g,\ft_\g) 
= \big\{ \lambda \in i\ft_\g' \: 
\lambda_c \in \Z, \ \lambda \in \cP(B_J)\big\}.\] 
\end{ex}

\subsection{The topology of the fixed point group $\hat \cL(K)^{\tilde\gamma}$}

Let $K$ be a $1$-connected simple infinite-dimensional Hilbert--Lie group and 
$\hat\cL(K) \cong \tilde\cL(K) \rtimes \R$ the 
simply connected Fr\'echet--Lie group with 
Lie algebra $\hat\cL(\fk)$ from Definition~\ref{def:3.5a}. 
Let $\hat\gamma \in \Aut(\hat\cL(K))$ be the automorphism 
induced by the automorphism $\gamma$ of $\cL(K)$ given by 
\[ \gamma(f)(t) = \phi\Big(f\big(t + \frac{2\pi}{N}\big)\Big) \quad \mbox{ and } \quad 
\L(\hat\gamma)(z,\xi,t) = (z, \L(\gamma)\xi,t).\]  
Then $\hat\cL(K)^{\hat\gamma} = \tilde\cL(K)^{\hat\gamma} \rtimes\R$ 
is a Lie subgroup 
with the Lie algebra $\tilde\cL(\fk)^{\L(\hat\gamma)} \rtimes \R$. 
Here we use that the central extension 
$\tilde\cL(K)$ of the locally exponential Lie group $\cL(K)$ is again locally 
exponential (cf.\ \cite{GN12} and also \cite[Thm.~IV.2.10]{Ne06}). 

\begin{prop} The inclusion 
$\hat\cL_\phi(\fk) \into \hat\cL(\fk)$ integrates to a Lie group morphism 
$\hat\cL_\phi(K) \to \hat\cL(K)$ whose kernel is the subgroup 
$C_N := \{ z \in Z \cong \T\: z^N = 1\}$ and whose range is 
$\hat\cL(K)^{\hat\gamma}$. In particular, 
$\pi_1(\hat\cL(K)^{\hat\gamma}) \cong C_N.$
\end{prop}

\begin{prf}
With the normal subgroup 
$\cL(K)_* := \{ f \in \cL(K) \: f(0) = \1\} \trile \cL(K)$ 
of based loops, we obtain the  semidirect decomposition 
$\cL(K) \cong \cL(K)_* \rtimes K$  
corresponding to the inclusion $K \into \cL(K)$ as the constant maps. 
Then 
$\tilde\cL(K)\cong \tilde\cL(K)_* \rtimes K,$ where $\tilde\cL(K)_*$ 
is a simply connected central $\T$-extension of $\cL(K)_*$.

As $\gamma$ commutes with the translation action 
of $\R$ on $\cL(K)$, $\hat\cL(K)^{\hat\gamma} \cong \tilde\cL(K)^{\hat\gamma} \rtimes \R,$ 
and we have a central extension 
\[ \1 \to \T \to 
\tilde\cL(K)^{\hat\gamma} \to 
\cL(K)^{\gamma} \to \1.\] 
In Remark~\ref{rem:homot} we have seen that 
$\cL_\phi(K) = \cL(K)^\gamma$ is $1$-connected, 
so that $\tilde\cL(K)^{\hat\gamma}$ is a central 
$\T$-extension of a $1$-connected Lie group. Therefore 
it is connected, and its fundamental group is isomorphic to the cokernel 
of the corresponding period homomorphism 
\[ \per \: \pi_2(\cL(K)^\gamma) \to \pi_1(\T) \cong 2\pi\Z \] 
(cf.\ \cite[Rem.~5.12(b)]{Ne02b}). 

Using Remark~\ref{rem:homot} again, we see that the triviality 
of $\pi_3(\phi)$ (Corollary~\ref{cor:2.10}) implies that we have a commutative diagram 
\[ \matr{ 
\Z \cong \pi_3(K) \cong\pi_2(\cL_\phi(K)_*) &\mapright{\cong} &\pi_2(\cL_\phi(K)) \\ 
\mapdown{} & & \mapdown{} \\
\Z \cong \pi_3(K) \cong \pi_2(\cL(K)_*) &\mapright{\cong} &\pi_2(\cL(K)).}\] 
For the group 
\[ \cL^c_\phi(K) := \Big\{ f \in C(\R,K) \: (\forall t \in \R)\ 
f\Big(t + \frac{2\pi}{N}\Big) = \phi^{-1}(f(t))\Big\} \] 
of continuous maps, it is easy to see that 
\[ \cL^c_\phi(K)_* \cong \Big\{ f \in C(\R,K)_* \: (\forall t \in \R)\ 
f\Big(t + \frac{2\pi}{N}\Big) = f(t)\Big\}, \] 
which coincides with the range of the map 
\[ \Phi \: \cL^c(K)_* \to \cL^c(K)_*, \quad 
\Phi(f)(t) := f(Nt).\] 
Since the inclusion of smooth into continuous maps induced isomorphisms of homotopy 
groups (cf.\ \cite[Cor.~3.4]{NeWo09}), \cite[Lemma~1.10]{MN03} now 
implies that the vertical arrows in the above 
diagram correspond to the endomorphism $\Z \to \Z, n \mapsto Nn$. 
We conclude that $\coker(\per) \cong \Z/N\Z$ and from that the assertion 
follows immediately. 
\end{prf}

\section{$d$-extremal weights} \mlabel{sec:4}

For $\g = \hat\cL_\phi(\fk)$, recall its root decomposition and 
the elements $c$ and $d$ from \eqref{eq:cd} in Example~\ref{ex:loop}. We 
write $\hat\cW  = \cW(\g,\hat\ft)$ for the {\it Weyl group of~$\g$.} 
In this section we derive a characterization of the set of those elements 
$\lambda \in i\ft_\g^*$  (the space of all linear functionals which are 
not necessarily continuous) which are $d$-minimal in the sense that 
$\lambda(d) = \min (\hat\cW\lambda)(d)$. In Sections~\ref{sec:5} and~\ref{sec:6} below, 
we  see that these elements parametrize the irreducible semibounded representations 
of $\hat\cL_\phi(K)$.

\begin{defn}  Let $\hat\cW$ denote the Weyl group of the pair $(\g,\ft_\g)$ 
(cf.\ Definition~\ref{def:1.8}). 
We call $\lambda \in i\ft_\g^*$ is {\it $d$-minimal} if 
$\lambda(d) = \min \la \hat\cW \lambda, d \ra$. 
\end{defn}

\begin{lem} \mlabel{lem:4.2} {\rm(\cite[Lemma~3.8]{HN12})}
If $(\hat\cW\lambda)(d)$ is bounded from 
below, then  $\lambda_c \geq 0$. If, in addition, 
$\lambda_c = 0$, then $\lambda$ is fixed by $\hat\cW$. 
\end{lem}

\begin{prop} \mlabel{prop:4.2} 
 {\rm(\cite[Cor.~3.6, Prop~3.9]{HN12})} 
Suppose that $(\Delta_\g)_c$ is one of the $7$ irreducible 
locally affine root systems with their natural $\Z$-grading. For 
$\lambda \in i\ft_\g^*$ with 
$\lambda_c > 0$, the following are equivalent: 
\begin{description}
\item[\rm(i)] $\lambda$ is $d$-minimal. 
\item[\rm(ii)] $(\forall \alpha \in \Delta(\fk,\ft), n =1,2)\quad 
(\alpha,n) \in(\Delta_\g)_c \Rarrow |\lambda(\check\alpha)| \frac{(\alpha,\alpha)}{2n}\leq \lambda_c$. 
\item[\rm(iii)] For $\uline\alpha = (\alpha,n) \in (\Delta_\g)_c$ with $n > 0$ we have
$\lambda(\check{\uline\alpha}) \leq 0$. 
\end{description}
\end{prop}

\begin{thm} \mlabel{thm:dmin} 
For the seven irreducible locally affine 
root systems $X_J^{(r)} = (\Delta_\g)_c \subeq i\ft_\g'$ of infinite rank, 
a linear functional $\lambda \in i \ft_\g^*$ with $\lambda_c > 0$ is $d$-minimal 
if and only if the following conditions are satisfied by 
the corresponding function $\oline\lambda \: J \to \R, j \mapsto \lambda_j$: 
\begin{description}
\item[\rm($A_J^{(1)}$)] 
$\max\oline\lambda - \min \oline\lambda \leq\lambda_c$. 
\item[\rm($B_J^{(1)}$)] $|\lambda_j| + |\lambda_k| \leq \lambda_c$ 
for $j\not=k$.
\item[\rm($C_J^{(1)}$)] $|\lambda_j| \leq \lambda_c$ 
for $j \in J$. 
\item[\rm($D_J^{(1)}$)] $|\lambda_j| + |\lambda_k| \leq \lambda_c$ 
for $j\not=k$.
\item[\rm($B_J^{(2)}$)] $|\lambda_j| \leq \lambda_c$ 
for $j\in J$. 
\item[\rm($C_J^{(2)}$)] $|\lambda_j| + |\lambda_k| \leq 2\lambda_c$ 
for $j\not=k$.
\item[\rm($BC_J^{(2)}$)] $|\lambda_j| \leq \lambda_c$ 
for $j \in J$. 
\end{description}
\end{thm}

\begin{prf} In all cases where the normalization of the scalar product 
is such that $\|\eps_j\|^2  = 1$ for every $j$, this follows 
immediately from \cite[Thm.~3.10]{HN12}, and this is the 
case for $A_J^{(1)}$, $B_J^{(1)}$ and $D_J^{(1)}$, where 
the long roots are of the form $\pm \eps_j \pm \eps_k$, and 
for $B_J^{(2)}$ it follows from Remark~\ref{rem:3.4}(d). 

In all other cases, the normalization for the scalar product $(\cdot,\cdot)_*$ in 
\cite[Thm.~3.10]{HN12} by $(\eps_j, \eps_j)_*= 1$ is different and we have to take a
closer look at the consequences. 
For $C_J^{(1)}$ the long roots are of the form $\pm 2 \eps_j$, which 
leads to the normalization $\|\eps_j\|^2 = \shalf$. 
For $C_J^{(2)}$ and $BC_J^{(2)}$ we have the same normalization by 
Remark~\ref{rem:3.4}(e),(f), which leads to 
$(\cdot, \cdot)_* = 2 (\cdot, \cdot)$ in these $3$ cases. 
The relation 
$|\lambda(\check \alpha)|\frac{(\alpha,\alpha)}{2n} \leq \lambda_c$
is therefore equivalent to 
$|\lambda(\check \alpha)|\frac{(\alpha,\alpha)_*}{2n} \leq 2\lambda_c.$ 
Replacing $\lambda_c$ in \cite[Thm.~3.10]{HN12} by $2\lambda_c$ now 
leads to the correct inequalities in these $3$ cases. 
\end{prf}

\begin{rem} \mlabel{rem:4.5} (a) The preceding theorem implies that 
$d$-minimal weights $\lambda \in i\ft_\g^*$ define bounded 
functions $\oline\lambda \: J \to \R$ and, moreover, that the boundedness of $\oline\lambda$ is 
equivalent to the existence of a $\lambda_c > 0$ such that 
a corresponding $\lambda \in \ft_\g^*$ is $d$-minimal. 

(b) If $\lambda  \in i\ft_\g^*$ satisfies 
$\lambda(\check \alpha) \in \Z$ for each $\alpha \in (\Delta_\g)_c$, 
then the subset 
$\lambda + \hat\cQ \subeq i \ft_\g^*$, where  
$\hat\cQ = \la \Delta_\g \ra_{\rm grp}$  is the {\it root group}, 
is invariant under the Weyl group $\hat\cW$. 
Therefore $(\hat\cW\lambda)(d) 
\subeq \lambda(d) + \Z$. 
If $(\hat\cW\lambda)(d)$ is bounded from below, we thus 
obtain the existence of a $d$-minimal element in $\hat\cW\lambda$. 
In particular, we obtain 
$\hat \cP^+ = \hat\cW \cP^+_d.$ 
\end{rem}

\section{Semibounded  representations of Hilbert loop groups} \mlabel{sec:5}

After the  preparations in the preceding sections, we now approach our goal of classifying the 
irreducible semibounded representations of 
$G = \hat\cL_\phi(K)$. The first major step is 
Theorem~\ref{thm:5.5}, asserting that for a semibounded 
representation $(\pi, \cH)$, the operator 
$\dd\pi(d)$ is either bounded from below (positive energy representations) 
or from above. 
Up to passing to the dual representation, we may therefore 
assume that we are in the first case. 
Then the minimal spectral value of $\dd\pi(d)$ turns out to be an eigenvalue 
and the group $Z_G(d)$ acts on the corresponding eigenspace, 
which leads to a bounded representation $(\rho,V)$ of this group. 
We then show that  $(\pi,\cH)$ can be reconstructed from $(\rho,V)$  
by holomorphic induction and that $\rho$ is irreducible if and only if $\pi$~is. 
Since an explicit 
classification of the bounded irreducible representations 
of the groups $Z_G(d)_0$ can be given in terms of 
$\cW$-orbits of extremal weights $\lambda$ (Theorem~\ref{thm:classif-hilbert}), 
the final step is to characterize 
those weights $\lambda$ for which the corresponding representation 
$(\rho_\lambda, V_\lambda)$ occurs. 

\subsection{From semibounded to bounded representations} 

Note that $d = (0,0,-i)$ satisfies $\exp(2\pi i \cdot d) \in \ker \Ad = Z(G)$. 
Therefore the following lemma can be used to obtain smooth eigenvectors of 
$\dd\pi(d)$ in irreducible representations of~$G$. The assumption of this 
lemma implies that $\pi(\exp \R x)\T$ is a torus, so that we know a priori 
that the Hilbert space decomposes into eigenspaces of this group. 

\begin{lem} \mlabel{lem:5.1} For a unitary representation $(\pi, \cH)$ of $G$ and 
$x \in \g$ with $\pi(\exp(T x)) = e^{i\mu}\1$ for some $T > 0$ and $\mu \in \R$, the  space 
$\cH^\infty$ of smooth vectors is invariant under the operators
\[ P_n(v) := \frac{1}{T} \int_0^T e^{-(2\pi n+ \mu)it/T} \pi(\exp tx)v\, dt\] 
which are orthogonal projections onto the eigenvectors of 
$-i\dd\pi(x)$ for the eigenvalues \break 
$(\mu + 2 \pi n)/T$, $n \in \Z$. 
\end{lem}

\begin{prf} (cf.\ the proof of \cite[Prop.~4.11]{Ne10b}) 
For $v \in \cH^\infty$, we have 
\[ \pi(g) P_n(v) = \frac{1}{T} \int_0^T e^{-(2\pi n+ \mu)it/T} \pi(g\exp tx)v\, dt,\] 
which is an integral of a smooth function on 
$[0,T] \times G$ over the compact factor $[0,T]$, which results in a smooth function on $G$. 
\end{prf}

\begin{thm} \mlabel{thm:5.5} 
Suppose that either $\phi  = \id_K$  or that 
$\phi$ is one of the three standard involutions. 
If $(\pi, \cH)$ is a semibounded unitary representation of 
$G = \hat\cL_\phi(K)$ for which $\dd\pi(c)$ is bounded, then 
$\dd\pi(d)$ is bounded from below, resp., above. 

If, in addition, $\pi$ is irreducible, 
then this is the case and, accordingly, 
the minimal/maximal spectral value of $\dd\pi(d)$ is an eigenvalue 
and the $K^\phi$-representation on the corresponding eigenspace is bounded. 
\end{thm}

\begin{prf} Since $\alpha_t := e^{t \ad i\cdot d}$ defines a continuous circle action on 
$\g$, the open invariant convex cone $W_\pi$ intersects the fixed point algebra 
$\z_\g(d)$ of this circle action. 
Since $[d,(0,\xi,0)] = -i(0,\xi',0)$,  
an element $(0,\xi,0) \in \hat\cL_\phi(\fk)$ 
commutes with $d$ if and only if $\xi$ is constant, 
i.e., its values are contained in $\fk^\phi$. We thus  have 
\[ \z_\g(d) = \ker (\ad i d) = \R i c \oplus \fk^\phi \oplus \R i \cdot d, \] 
which is a Hilbert--Lie algebra. 

Since every open invariant cone in the Hilbert--Lie algebra 
$\fz_\g(d)$ intersects the center (\cite[Prop.~A.2]{Ne11b}), 
the non-empty open invariant cone 
$W_\pi \cap \fz_\g(d)$ actually intersects the subspace 
$\R i c\oplus \fz(\fk^\phi)\oplus \R i\cdot d$. 
For the three types of standard involutions we have 
$\fz(\fk^\phi) = \{0\}$ 
(cf.\ Examples~\ref{ex:auto-bj2}, \ref{ex:auto-cj2} and \ref{ex:auto-bcj2}), 
so that 
\[ W_\pi \cap (\R i c + \R i \cdot d) \not=\eset.\] 
If, in addition, $\dd\pi(c)$ is bounded, then 
$i c +  W_\pi = W_\pi$ leads to $i \cdot d \in W_\pi \cup - W_\pi$. 
In particular, $\dd\pi(d)$ is bounded from below or above. 

Now we assume that $\pi$ is irreducible and w.l.o.g.\ that 
$i\cdot d \in W_\pi$. Then 
$\pi(\exp 2 \pi i \cdot d) \in \pi(Z(G)) \subeq \T\1$ by Schur's Lemma 
and Lemma~\ref{lem:5.1} implies that the minimal spectral value $\mu$ of 
$\dd\pi(d)$ is an eigenvalue. Let $V := \ker(\dd\pi(d) - \mu \1)$ be the
corresponding eigenspace. It is invariant under the subgroup 
$Z_G(d)$ which leads to a unitary representation  
$(\rho,V)$ of this group. The corresponding open convex cone 
$W_\rho$ satisfies 
\[ i \cdot d \in W_\pi \cap \z_\g(d) \subeq W_\rho,\] 
but since $\dd\rho(d) = \mu \1$, this leads to 
$0 \in W_\rho - i\cdot d = W_\rho$ and thus to $W_\rho = \fz_\g(d)$, i.e., 
$\rho$ is bounded. 
\end{prf}

\begin{rem}
  \mlabel{rem:3.5} 
(a) Since $Z(G)_0 \cong \T$ by Theorem~\ref{thm:cent-ext}, 
any unitary representation of $G$ is a direct sum of $\dd\pi(c)$-eigenspaces, so that we 
can easily reduce to the situation where $\dd\pi(c) \in \R \1$. 

(b) It is a key point in the proof of Theorem~\ref{thm:5.5} that 
$\fz(\fk^\phi) = \{0\}$ which holds for all $3$ standard involutions. 
For any finite order automorphism 
$\phi$ for which $\fz(\fk^\phi) = \{0\}$, the argument in the proof 
of Theorem~\ref{thm:5.5} goes through, and even 
Theorem~\ref{thm:5.4} below remains valid. 
This is in particular the case if $\fk$ is abelian and 
$\fk^\phi = \{0\}$. With automorphisms of the form 
$\phi(g) = U gU^{-1}$, where $U$ is unitary of order $N$ with some finite-dimensional eigenspaces, 
one obtains examples where $\fz(\fk^\phi) \not=\{0\}$. 
\end{rem}

Our next step is to explain why irreducible semibounded representations with 
$i\cdot d \in W_\pi$ are 
uniquely determined by the representations on the minimal eigenspaces of $\dd\pi(d)$. 
This requires the technique of holomorphic induction from Appendix~\ref{app:c}.

\subsection{Holomorphic induction for $\hat\cL_\phi(K)$} 

Let $\g_B := \tilde\cL^H_\phi(\fk)$ be the central extension of 
$H^1$-loop algebra $\cL_\phi^H(\fk)$ 
from  Definition~\ref{def:a1.1} which is a Banach--Lie algebra. 
According to Theorem~\ref{thm:cent-ext-h1}, there exists a corresponding $1$-connected 
Banach--Lie group $G_B := \tilde\cL_\phi^H(K)$ and also a complex group 
$(G_B)_\C = \tilde\cL_\phi^H(K_\C)$ (Remark~\ref{rem:polar-h1}). 
Then 
\[ \fp_B^\pm := \oline{\sum_{\pm n > 0} e_n \otimes \fk_\C} \] 
are closed subalgebras of $(\g_B)_\C$. 
We also put 
\[ \fh := \fh_B  := \R i c + \fk^\phi \quad \mbox{ and } \quad \fq_B := \fh_\C \ltimes \fp_B^+.\] 
The Fourier expansion of $H^1$-loops implies that 
$\g_B$ satisfies the splitting condition 
\eqref{eq:splitcond} from Appendix~\ref{app:c}:  
\[ (\g_B)_\C = \fp_B^+ \oplus \fh_\C \oplus \fp_B^-. \] 
Therefore all assumption of Example~\ref{ex:c.4}(a) are satisfied. 

On $G_B$ we now consider the one-parameter group 
$\alpha \: \R \to \Aut(\tilde\cL_\phi^H(K))$ 
defining the translation action $(\alpha_T f)(t) = f(t+T)$ of $\R$. 
Then $\fp_B^\pm$, $\fh$ and $\fq_B$ are $\alpha$-invariant, so that 
Example~\ref{ex:c.4}(b) applies. 
On the Lie algebra level, the subspace $\g \subeq \g_B$ of smooth vectors 
for $\alpha$ coincides with the Fr\'echet--Lie algebra $\tilde\cL_\phi(\fk)$ 
defined by smooth loops, and $\fh$ corresponds to 
the Lie subalgebra of $\alpha$-invariant elements. Therefore the concepts and results 
 from Appendix~\ref{app:c} 
concerning holomorphic induction are now available for the pair $(G, Z_G(d)_0)$, 
resp., 
the complex homogeneous space 
\[ \hat\cL_\phi(K)/Z_G(d)_0 
\cong \tilde\cL_\phi(K)/Z_{\tilde\cL_\phi(K)}(d)_0 
\cong \cL_\phi(K)/(K^\phi)_0.\] 
In view of Lemma~\ref{lem:kphicon} we know that 
the subgroup $K^\phi$ is connected, so that 
$\cL_\phi(K)/K^\phi  \cong \cL_\phi(K)/(K^\phi)_0$ is actually 
simply connected because $\cL_\phi(K)$ is $1$-connected.

\begin{thm}
  \mlabel{thm:5.4} 
Every irreducible semibounded unitary representation 
$(\pi, \cH)$ of $\hat\cL_\phi(K)$ for which $\dd\pi(d)$ is bounded from below 
is holomorphically induced from the bounded representation $(\rho, V)$ of 
$H = Z_G(d)_0$ on the minimal eigenspace of $\dd\pi(d)$. 
\end{thm}

\begin{prf} We want to apply Theorem~\ref{thm:c.3}.  
We know from Theorem~\ref{thm:5.5} that $(\rho,V)$ 
is a bounded representation of $H$, which implies (HI1). 
So we first use Lemma~\ref{lem:5.1} to see that the projection 
$P_0 \: \cH \to V$ onto the eigenspace $V = \ker(\dd\pi(d)-\mu \1)$ 
maps $\cH^\infty$ onto 
$\cH^\infty \cap V$, and since $P_0$ is continuous and $\cH^\infty$ 
is dense in $\cH$, 
$\cH^\infty \cap V$ is dense in $V$. Let $P_n \: \cH \to \cH_n$, $n \in \Z$, 
denote the projections onto the other eigenspaces of $\exp(\R i d)$ from 
Lemma~\ref{lem:5.1}. As $V$ is the minimal eigenspace 
of the diagonalizable operator $\dd\pi(d)$, the fact that 
$\dd\pi(\g_\C^k) \cH_n \subeq \cH_{n+k}$ for $k,n \in \Z$ implies that 
$V \cap \cH^\infty \subeq (\cH^\infty)^{\fp^-}$ for 
$\fp^- = \fp_B^- \cap \g_\C$. This proves (HI2). 
Finally (HI3) follows from the irreducibility of~$(\pi, \cH)$.
\end{prf}

In view of the preceding theorem, a classification of the irreducible semibounded 
representations of $\hat\cL_\phi(K)$ now consists in a classification of the 
irreducible bounded representations $(\rho,V)$ of $Z_G(d)_0$ which are inducible 
in the sense of Definition~\ref{def:d.1}. It is easy to pinpoint a necessary 
condition for inducibility. 

\begin{defn}
 We say that a representation 
$(\rho,V)$ of $\fz_\g(d)$ is {\it $d$-minimal} if 
$\rho(\ft_\g)$ is diagonalizable and all $\ft_\g$-weights of $\rho$ 
are $d$-minimal. 
\end{defn}

\begin{prop}  \mlabel{prop:5.5} 
For $y \in \g_\C = \fp^+ \oplus \fh_\C \oplus \fp^-$, we 
write $y = y_+ + y_0 + y_-$ for the corresponding decomposition. 
If a bounded unitary representation $(\rho,V)$ of $H = Z_G(d)_0$ 
is holomorphically inducible for $\fq = \fp^+ \rtimes \fh_\C$, 
then \begin{equation}
  \label{eq:inducness}
\dd\rho([z^*,z]_0) \geq 0 \quad \mbox{ for } \quad z \in \fp^+,  
\end{equation}
and this implies that it is $d$-minimal, provided $\rho(\ft_\fg)$ is diagonalizable. 
\end{prop}

\begin{prf} Suppose that $(\pi, \cH)$ is obtained from $(\rho,V)$ by holomorphic induction. 
Then $V \subeq (\cH^\infty)^{\fp^-}$, so that we obtain for $v \in V$ 
and $z \in \fp^+$ 
\begin{align*}
\la \dd\rho([z^*,z]_0)v,v\ra 
&= \la \dd\rho([z^*,z])v,v\ra 
=  \la [\dd\pi(z^*), \dd\pi(z)]v,v\ra \\
&=  \la \dd\pi(z^*)\dd\pi(z)v,v\ra 
= \|\dd\pi(z)v\|^2 \geq 0.
\end{align*}
This proves the necessity of \eqref{eq:inducness}. 

For every weight vector $v_\mu \in V$ with weight $\mu \in i \ft_\g'$ and 
$\uline \alpha = (\alpha, n) \in (\Delta_\g)_c$ with $n > 0$ we pick 
$x \in \g_\C^{\uline \alpha}$ such that 
$[x,x^*] = \check{\uline\alpha}$ (cf.\ Example~\ref{ex:twistloop}). Then 
\eqref{eq:inducness} implies 
$\mu(\check{\uline\alpha}) \leq 0$. In view of 
Proposition~\ref{prop:4.2}(iii), this is equivalent to the $d$-minimality of $\mu$. 
\end{prf} 

\subsection{Bounded representations of $K^\phi$} 

\begin{rem} \mlabel{rem:5.6} 
Let $\Delta = \Delta(\fk,\ft)$ be 
a root system of type $A_J, B_J, C_J$ or $D_J$. 
We represent a corresponding integral weight 
as a function $\lambda \: J \to \R$ and 
observe that 
\[ A_J \subeq D_J  = B_J \cap C_J.\] 
Then the integrality with respect to $A_J$ 
means that $\lambda_j - \lambda_k \in \Z$ for 
$j\not=k \in J$. Since $J$ is infinite in our context, 
the requirement $\lambda \in i \ft' \cong \ell^2(J,\R)$ 
implies that $\lambda$ is finitely supported with values in~$\Z$. 
This in turn implies that $\lambda$ is an integral weight 
for $A_J$, $B_J$, $C_J$ and $D_J$. 
\end{rem}

\begin{prop} \mlabel{prop:classworb} {\rm(Classification of $\cW$-orbits)} 
For $\Delta = \Delta(\fk,\ft)$ of 
type  $A_J, B_J, C_J$ or $D_J$ the corresponding set 
of integral weights $\cP(\fk,\ft) \subeq i\ft' \cong \ell^2(J,\R)$ 
coincides with $\ell^2(J,\Z) \cong \Z^{(J)}$. For the $\cW$-action 
on this set, we have the following set of invariants which is complete 
in the sense that it separates the $\cW$-orbits in $\cP(\fk,\ft)$: 
\begin{itemize}
\item $A_J:\quad m_n(\lambda) := |\{ j \in J \: \lambda_j = n\}|$ 
for $0 \not= n \in \Z$. 
\item $B_J, C_J$ and $D_J$: $\quad m_n(\lambda) 
:= |\{ j \in J \: |\lambda_j| = n\}|$ 
for $n \in \N$. 
\end{itemize}
\end{prop}

\begin{prf} (a) For $A_J$, the functions $m_n$ are constant on the 
orbits of $\cW \cong S_{(J)}$ and, conversely, if 
$m_n(\lambda) = m_n(\lambda')$ for $\lambda, \lambda' \in \cP(\fk,\ft)$, 
then $\lambda' \in \cW \lambda$ follows from the finiteness of the 
support of~$\lambda$. 

(b) For the root systems $B_J$, $C_J$ and $D_J$, the functions 
$m_n$, $n \in \N$, are constant on the $\cW$-orbits, and since 
every $\lambda \in \cP(\fk,\ft)$ is finitely supported, 
its $\cW$-orbit contains a non-negative 
element. Hence 
$m_n(\lambda) = m_n(\lambda')$ for $\lambda, \lambda' \in \cP(\fk,\ft)$ 
and every $n \in \N$ leads to $\lambda' \in \cW \lambda$. 
\end{prf}

\begin{thm} \mlabel{thm:classif-hilbert} 
Let $K$ be a simple Hilbert--Lie group with Lie algebra $\fk$ and 
$\ft \subeq \fk$ maximal abelian with root system $\Delta \subeq i \ft'$.  
Then every bounded unitary representation of $K$ is a direct sum 
of irreducible ones. The irreducible representations 
$(\rho_\lambda, V_\lambda)$ can be parametrized by their extremal weights 
$\lambda \in \cP(\fk,\ft)$  as follows. 
If $\cQ := \la \Delta \ra_{\rm grp} \subeq i \ft'$ is the root group, 
then the weight set $\cP_\lambda$ of $\rho_\lambda$ satisfies 
\[ \cP_\lambda = \conv(\cW \lambda) \cap (\lambda + \cQ) \quad \mbox{ and } \quad 
\Ext(\conv(\cP_\lambda)) = \cW \lambda.\] 
We have $\rho_\lambda \sim \rho_\mu$ if and only if $\mu \in \cW\lambda$, so that the 
irreducible bounded unitary representations of $K$ are classified by 
the set $\cP(\fk,\ft)/\cW$ of $\cW$-orbits in $\cP(\fk,\ft)$. 
All these representations factor through the adjoint group $K/Z(K)$. 
\end{thm}

\begin{prf} In view of the classification of simple Hilbert--Lie algebras, 
the assertion on the classification follows from 
\cite[Thm.~III.14]{Ne98} for $\fk = \fu_2(\cH)$ 
and from \cite[Thms.~D.5, D.6]{Ne11c} for $\fk = \fo_2(\cH)$ and $\sp_2(\cH)$. 

That all these representations factor through the adjoint group is trivial 
for $\fk = \sp_2(\cH)$ because in this case the center of the corresponding 
simply connected group $\Sp_2(\cH)$ is trivial 
(Theorem~\ref{thm:1.5}). For $\fk = \fu_2(\cH)$ it follows 
from \cite[Rem.~D.2]{Ne11c}, and for $\fk = \fo_2(\cH)$ the description of the corresponding 
highest weights (cf.\ Remark~\ref{rem:5.6}) 
implies that they are contained in the root group $\cQ$, 
and hence that the corresponding representation is trivial on the center. 
\end{prf}

For the $3$ standard involutions $\phi$ 
(cf.\ Examples~\ref{ex:auto-bj2}-\ref{ex:auto-bcj2}), 
the Lie algebra $\fk^\phi$ is simple, so that 
we obtain an explicit description of the bounded irreducible representations of the groups 
$K^\phi$, resp., $Z_G(d)_0$ in all seven cases. 
According to Theorem~\ref{thm:5.4}, 
any irreducible semibounded representations of $G = \hat\cL_\phi(K)$ 
for which $\dd\pi(d)$ is bounded from below is holomorphically induced 
from the representation $(\rho,V)$ on the minimal eigenspace, 
hence uniquely determined by this representation 
(cf.\ Definition~\ref{def:d.1}). 
It therefore remains to identify those bounded 
representations of $Z_G(d)_0$ which are holomorphically inducible.
 
\subsection{Characterization of the inducible bounded representations} 

In this subsection we show that any $d$-minimal bounded representation 
$(\rho,V)$ of $Z_G(d)_0$  is inducible. 

\begin{thm} \mlabel{thm:5.9} An irreducible bounded unitary representation 
$(\rho,V)$ of 
$Z_G(d)_0$ is holomorphically inducible 
if and only if it is $d$-minimal. 
\end{thm}

\begin{prf} We have already seen in Proposition~\ref{prop:5.5} 
that $\rho$ is $d$-minimal if it is holomorphically inducible. 
Now we assume that $(\rho,V)$ is a $d$-minimal bounded 
representation of $H = Z_G(d)_0$ of extremal weight 
$\lambda \in \hat T_G = \Hom(T_G,\T)$ (recall $T_G = \exp \ft_\g$). Then $\lambda$ is $d$-minimal, 
so that 
$\lambda(\check{\uline\alpha}) \geq 0$ for 
$\uline\alpha = (\alpha,n)$, $n < 0$. This means that 
\[ \fp_\lambda := (\ft_\g)_\C 
+ \sum_{\lambda(\uline\alpha^\sharp) \geq 0} \g_\C^{\uline\alpha} 
\supeq  \sum_{\alpha \in \Delta, n > 0} \g_\C^{(\alpha,n)}.\] 
If $\lambda_c = 0$, then $\lambda$ vanishes on $\check\Delta_\g$ (Lemma~\ref{lem:4.2}) 
which implies that $G$ has a one-dimensional representation $(\pi, \cH)$ 
for which  $\dd\pi \: \g_\C \to \End(\cH) \cong \C$ extends $\lambda$. 
We may therefore assume that $\lambda_c \not=0$. 

In view of Theorem~\ref{thm:c.2}, 
it suffices to show  that the corresponding linear map 
\[ \beta \: U(\g_\C) \to B(V) \] 
that vanishes on 
$\fp^+ U(\g_\C) + U(\g_\C) \fp^-$ and satisfies 
$\beta\res_{U(\fh_\C)} = \dd\rho$ for $\fh = \fz_\g(d)$, is positive definite on the 
$*$-algebra $U(\g_\C)$ (cf.\ Definition~\ref{def:a.1}(c)). 

Let $\ft_\C^{\rm alg} :=  \C i \cdot d + \Spann_\C \check\Delta_\g$ and 
\[ \g_\C^{\rm alg} := \C  d + \la \g_\C^\alpha \: \alpha \in (\Delta_\g)_c\ra_{\rm Lie\ alg},\] 
and observe that this is a Lie algebra with a root decomposition 
with respect to $\ft_\C^{\rm alg}$. 
It is a coral locally affine complex Lie algebra 
in the sense of \cite[Def.~3.1]{Ne10a} and 
$\lambda$ defines an integral weight of 
$\g_\C^{\rm alg}$ for which $\lambda_c \not=0$. 
Therefore \cite[Thm.~4.11]{Ne10a} implies the existence of a 
unitary extremal weight module 
$(\pi_\lambda, L(\lambda))$ of $\g_\C$ generated by a $\fp_\lambda$-eigenvector 
$v_\lambda$ of weight $\lambda$. 
Note that \cite[Thm.~5.7]{Ne10a} shows that 
$\g_\C^{\rm alg}$ is a locally extended affine Lie algebra with root system $\Delta_\g$ 
in the sense of \cite{MY06} (see also \cite{Ne10a}). 

Now 
$V_\lambda := U(\fh_\C)v_\lambda$ is an $\fh_\C$-module of extremal weight 
$\lambda$, so that we may identify it with a dense subspace 
of $V$. For $\fp^\pm_{\rm alg} := \fp^\pm \cap \g_\C^{\rm alg}$, 
the relation $v_\lambda \in L(\lambda)^{\fp^-_{\rm alg}}$ implies that 
$V_\lambda \subeq L(\lambda)^{\fp^-_{\rm alg}}$, and 
\[ L(\lambda) 
= U(\g_\C^{\rm alg})v_\lambda 
= U(\fp^+_{\rm alg})V_\lambda 
\subeq V_\lambda + \fp^+_{\rm alg}L(\lambda) \] 
shows that $V_\lambda$ is the minimal $d$-eigenspace in 
$L(\lambda)$. Let $p_V \: L(\lambda) \to V_\lambda \subeq V$ denote the 
orthogonal projection. Then 
\[ \gamma(D) := p_V \pi_\lambda(D) p_V^* \] 
satisfies 
$\fp^+_{\rm alg} U(\g_\C^{\rm alg}) + U(\g_\C^{\rm alg}) \fp^-_{\rm alg} \subeq \ker \gamma$ 
and 
$\gamma\res_{U(\fh_\C^{\rm alg})} = \rho_\lambda$. 
Since the representation $\pi_\lambda$ on $L(\lambda)$ is unitary, 
$\gamma$ is positive definite. 
Since all maps 
\[ \g_\C^k \to B(V), \quad 
(x_1, \ldots, x_k) \mapsto \beta(x_1 \cdots x_k) \] 
are continuous and the restriction of 
$\beta$ to the subalgebra $U(\g_\C^{\rm alg})$ coincides with 
$\gamma$, it follows that $\beta$ is also positive definite. 
Now the assertion follows from Theorem~\ref{thm:c.2}. 
\end{prf}

\begin{rem} \mlabel{rem:5.9} (a) Consider the Banach--Lie group $\tilde\cL_\phi^H(K)$ 
from Appendix~\ref{app:a} and the subgroup 
$\T\times K^\phi$ corresponding to the centrally extended Lie algebra 
$\R \times \fk^\phi \subeq \tilde\cL_\phi^H(\fk)$. 
Suppose that $(\rho,V)$ is a bounded representation of 
$H := \T\times K^\phi$ which is holomorphically inducible to the 
Fr\'echet--Lie group $\tilde\cL_\phi(K)$. 
Since $\tilde\cL_\phi(\fk)$ is dense in $\tilde\cL_\phi^H(\fk)$, 
the fact that the conditions in Theorem~\ref{thm:c.2} are satisfied 
for $\tilde\cL_\phi(K)$ immediately implies that they are also satisfied for 
the bigger group $\tilde\cL_\phi^H(K)$. 
Therefore the holomorphically induced representation
 $(\pi,\cH)$ of $\tilde\cL_\phi(K)$ extends to a holomorphically 
induced representation of the Banach--Lie group 
$\tilde\cL_\phi^H(K)$, and we thus obtain a continuous unitary representation 
of the topological group $\hat\cL_\phi^H(K)$. 

(b) The preceding argument also shows that, if $\rho$ is irreducible, then the 
same holds for the corresponding holomorphically induced representation of 
$\tilde\cL_\phi(K)$ resp., $\tilde\cL_\phi^H(K)$. 

(c) Assume that $\phi = \id$. Then we 
can also ask about the restriction of $\pi$ to the subgroup 
$L := \tilde\cL(K)_*$ corresponding to functions vanishing in $\1$. 
As $\tilde\cL(K) = L \rtimes K$ corresponding to functions vanishing in $\1$, 
the group $L$ acts transitively on the complex homogeneous space 
$\cL(K)/K$ which implies that $\pi\res_L$ is holomorphically induced 
from the trivial representation of $L \cap K = \{\1\}$ on $V$. This leads 
to $\pi(L)' \cong B(V)$ (Theorem~\ref{thm:c.1}(ii)) which implies in particular
that $\pi\res_L$ is irreducible if and only if $\dim V = 1$. 
\end{rem}

\subsection{Semibounded representations of one-dimensional extensions} 
\mlabel{sec:d.2.1}

In this subsection we provide a few results supporting the point of view that, 
without the double extension, the representation theory of 
loop groups is much less interesting. We show that all 
semibounded unitary representations of the central extension  $\tilde\cL_\phi(\fk)$  
are trivial on the center and factor through bounded representations 
of $\cL(\fk)$. One can actually show that these 
are finite-dimensional and tensor products of evaluation representations. 
For those representation extending to the Lie algebra $\cL^c(\fk)$ 
of continuous maps, this follows from \cite{NS11}. 
We also show that all semibounded 
representations of $\cL_\phi(K) \rtimes_\alpha \R$ are trivial 
on $\cL_\phi(K)$.

\begin{lem} \mlabel{lem:solvlem} Let $(\pi, \cH)$ be a smooth representation of a 
Lie group $G$. Then the following assertions hold: 
\begin{description}
  \item[\rm(i)] Let $x,y \in \g$ with $[x,[x,y]] = 0$. 
If $-i\dd\pi(y)$ is bounded from below, then $\dd\pi([x,y]) = 0$. 
  \item[\rm(ii)] If $\g$ is $2$-step nilpotent and 
$\pi$ is semibounded, then $[\g,\g] \subeq \ker(\dd\pi)$. 
\end{description}
\end{lem}

\begin{prf} (i) For a smooth unit vector $v \in \cH^\infty$, we consider the 
continuous linear functional $\lambda(z) := \la -i \dd\pi(z)v,v\ra$ 
on $\g$. Then our assumption implies that 
$\lambda$ is bounded from below on $\Ad(G)y$ which contains 
$\Ad(\exp \R x)y = y + \R [x,y]$. This leads to 
$\lambda([x,y]) = 0$. We thus obtain $\dd\pi([x,y]) = 0$.   

(ii) Pick $y \in W_\pi$. For every $x \in \g$ we then have 
$[x,[x,y]] = 0$, so that (i) leads to 
$\dd\pi([x,y]) =0$ and thus $[W_\pi, \g] \subeq \ker(\dd\pi)$. 
As $W_\pi$ is open, the assertion follows. 
\end{prf}

\begin{prop} \mlabel{prop:trivcone} Let $\fk$ be a simple Hilbert--Lie algebra 
and $\phi \in \Aut(\fk)$. 
Then all open invariant cones in $\cL_\phi(\fk)$ are trivial.   
\end{prop}

\begin{prf} Let $\eset \not= W \subeq \cL_\phi(\fk)$ be an open invariant convex cone. 

(a) First we consider the case $\phi = \id$ and show that, for every 
compact manifold $M$ with or without boundary, 
all open invariant cones in $C^\infty(M,\fk)$ are trivial. 
Since $\fk =\fu_2(\cH)$ for an infinite dimensional real, complex or quaternionic 
Hilbert space, and the union of the subalgebras 
$\su(\cH_F)$, where $\cH_F \subeq \cH$ is a finite-dimensional 
subspace, is dense in $\fu(\cH)$, 
the union of the subalgebras 
$C^\infty(M,\su(\cH_F))$ is dense in $C^\infty(M,\fk)$. 
Hence there exists a subspace $\cH_F$ for which 
$\eset \not= W_F := W \cap C^\infty(M,\su(\cH_F))$. 
Then $W_F$ is invariant under conjugation with the compact group $\SU(\cH_F)$, 
hence contains a $\SU(\cH_F)$-fixed point. Since 
$\su(\cH_F)$ has trivial center, $0$ is the only fixed point, 
and thus $0 \in W_F \subeq W$. This in turn implies that 
$W = \cL(\fk)$. 

(b) For the general case, we consider the closed interval 
$I := [-a, a]$ for $0 < a < \frac{\pi}{N}$. 
Then we have a continuous restriction map 
\[ R \: \cL_\phi(\fk) \to C^\infty(I,\fk),\] 
which is also surjective (cf.\ \cite[Cor.~III.7]{Wo06}), hence open by the 
Open Mapping Theorem. Therefore 
$R(W)$ is an open invariant cone in $C^\infty(I,\fk)$, and (a) 
implies that $0 \in R(W)$, which in turn implies that 
$W \cap \ker R \not=\eset.$ 

For $b := \frac{2\pi}{N}$, we have 
\[ \ker R 
= \{ f \in \cL_\phi(\fk) \: f\res_{[-a,a]} = 0\}
\cong \{ f \in C^\infty([0,b],\fk) \: f\res_{[0,a]} = 0 = f_{[b-a,b]}\}.\] 
For every $k \in K$ there exists an element $f \in \cL_\phi(\fk)$ 
restricting to the constant function $k$ on $[a,b-a]$, so that 
$W \cap \ker R$ is invariant under conjugation with constant functions 
in $K$. Passing to a sufficiently large finite dimensional 
subalgebra $\fk_F \subeq \fk$ and averaging over the action of the 
corresponding compact group $K_F$, it follows as in (a) that 
$0 \in W \cap \ker R \subeq W$, so that $W = \cL_\phi(\fk)$.  
\end{prf}

\begin{cor} If $K$ is a simple Hilbert--Lie group, then all 
semibounded unitary representations of $\cL_\phi(K)$ are bounded.   
\end{cor}

\begin{thm} \mlabel{thm:centriv} 
{\rm(Semibounded representations of central extensions)} 
If $\fk$ is a simple Hilbert--Lie algebra, then 
all semibounded unitary representations of the central extension 
$\tilde\cL_\phi(\K)$ are trivial on the center and bounded. 
\end{thm}

\begin{prf} Localization on the center reduces the problem to representations 
which are bounded on the center, so that $\R i c + W_\pi = W_\pi$. 
Hence $W_\pi$ defines an open invariant cone in $\cL_\phi(\fk) 
\cong \tilde\cL_\phi(\fk)/\R i c$, which is trivial by Proposition~\ref{prop:trivcone}. 
Therefore $\pi$ is bounded. 
In particular, the restriction of $\pi$ 
to the $2$-step nilpotent group $\tilde\cL_\phi(\ft_\fk)$ 
is bounded. Since it is $2$-step nilpotent, $\dd\pi$ is trivial on the commutator 
algebra (cf.\ Lemma~\ref{lem:solvlem}(ii)), and thus $\dd\pi(c)=~0$. 
We conclude that $\dd\pi(c)$ vanishes for all semibounded representations, and 
hence also that these representations are bounded. 
\end{prf}

The preceding result shows that the central extension 
$\tilde\cL_\phi(K)$ and $\cL_\phi(K)$ have the same (semi-)bounded representations. 
In a similar vein, extending $\cL_\phi(K)$ to the semidirect product defined by 
the translation actions only leads to trivial semibounded representations. 

\begin{thm} \mlabel{thm:semdirtriv} 
{\rm(Semibounded representations of semidirect products)} 
If $\fk$ is a simple Hilbert--Lie algebra and 
$\phi$ any finite order automorphism of the corresponding 
simply connected group $K$ for which 
$\fz(\fk^\phi) = \{0\}$, then 
every unitary representation $(\pi, \cH)$ of 
\break $\cL_\phi(K) \rtimes_\alpha \R$ for which 
$-i\dd\pi(0,1)$ is bounded from below is trivial on 
$\cL_\phi(K)$. 
\end{thm}

\begin{prf} For any abelian $\phi$-invariant subalgebra $\fa \subeq \fk$, 
we consider the $2$-step solvable Lie algebra 
$\cL_\phi(\fa)$ and note that 
$\cL_\phi(\fa) 
= [d,\cL_\phi(\fa)] \oplus \fz_{\cL_\phi(\fa)}(d) 
= [d,\cL_\phi(\fa)] \oplus \fa^\phi$. 
Therefore Lemma~\ref{lem:solvlem} implies that 
$\cL_\phi(\fa) \subeq \fa^\phi + \ker(\dd \pi)$. 

Applying this observation to one-dimensional subalgebras 
$\fa= \R x \subeq \fk^\phi$, we obtain 
$\cL_\phi(\fk^\phi) = \cL(\fk^\phi) \subeq \fk^\phi + \ker(\dd \pi)$. 
As $\fk^\phi$ is topologically perfect, it is contained in the 
ideal of $\cL(\fk^\phi)$ generated by 
$\fz_d(\cL(\fk^\phi))$. This leads to  
$\cL(\fk^\phi) \subeq \ker(\dd\pi)$. 

For the abelian subalgebra $\ft_\fk \subeq \fk$ 
we likewise obtain 
$\cL_\phi(\ft_\fk) \subeq \ft +\ker(\dd \pi) = \ft_\fk + \ker(\dd\pi)$. 
Hence $\ft \subeq \fk^\phi \subeq \ker(\dd\pi)$ enventually yields 
$\cL_\phi(\ft_\fk) \subeq \ker(\dd \pi)$. 
We finally arrive at $\cL_\phi(\fk) \subeq \ft_\fk + [\ft, \cL_\phi(\fk)] \subeq \ker(\dd\pi)$. 
\end{prf}

\section{Semiboundedness of holomorphically induced representations} 
\mlabel{sec:6}

We are now ready to complete the picture by showing that 
the irreducible $G$-representations  
$(\pi_\lambda, \cH_\lambda)$ obtained by holomorphic induction from 
$d$-minimal representations $(\rho_\lambda, V_\lambda)$ are semibounded. 

\begin{thm} \mlabel{thm:6.1} 
Let $K$ be a $1$-connected Hilbert--Lie group 
and $(\pi,\cH)$ be a unitary representation of $G = \hat\cL_\phi(K)$ 
which is holomorphically induced from the bounded representation $(\rho,V)$ of 
$H = Z_G(d)_0$ for which $\dd\rho(d) = \mu \1$ for some $\mu \in \R$. 
Then $(\pi, \cH)$ is semibounded with $i\cdot d \in W_\pi$. 
\end{thm}

\begin{prf} Recall the subalgebras $\fp_B^\pm$ and $\fp^\pm$ from Section~\ref{sec:5}. 
We note that the representation $\ad_{\fp^+}$ of the Hilbert--Lie algebra 
$\fh = \fz_\g(d)$ on the Hilbert space 
$\fp_B^+ = \oline{\sum_{n > 0} \g_\C^n} \subeq(\g_B)_\C$ 
is unitary, $i c$ acts trivially, $d$ acts by $n \cdot\id$ on $\g_\C^n 
 = e_n \otimes \fk_\C^{n}$, and 
$\fk^\phi$ acts by the adjoint representation on $\g_\C^n \cong \fk_\C^n$. 
Hence an element 
$x = (z,x_0,t) \in \fz_\g(d)$ satisfies 
\[ -i\ad_{\g_\C^n}(z,x_0,t) = t n -i \ad_{\fk_\C}x_0 \geq 0 \quad \mbox{ for every } n > 0\quad \mbox{ if } \quad \|\ad x_0\|\leq t.\]
 Therefore the elements $x \in \fz_\g(d)$ with this property 
form a closed invariant cone with non-empty interior $C$. 

For the $\dd\pi(d)$-eigenspace decomposition 
\[ \cH = \hat\bigoplus_{n \in \N_0} \cH_n \quad 
\mbox{ with } \quad 
\cH_n = \ker(\dd\pi(d) - (\mu + n)\1), \] 
all subspaces $\cH_n$ are invariant under $Z_G(d)_0$. 
For every $v \in V = \cH_0 \subeq \cH^\omega$, the Poincar\'e--Birkhoff--Witt 
Theorem shows that 
\[ U(\g_\C)v = U(\fp^+)U(\fh_\C)U(\fp^-)v 
= U(\fp^+)U(\fh_\C)v \subeq U(\fp^+)V\] 
is dense in $\cH$. Therefore $\cH_n$, as a unitary representation of 
$Z_G(d)_0$ containing a dense subspace which 
is a quotient of of the bounded unitary representation on the Hilbert space 
\[ \Big(\bigoplus_{0 \leq k \leq n} (\fp_B^+)^{\hat\otimes k} \Big)_n \hat\otimes V,\] 
is a bounded representation of $Z_G(d)_0$.
For $x \in C$, the spectrum of $-i x$ on the left hand factor 
is non-negative, so that 
\[ \inf(\Spec(-i\dd\pi(x)) = \inf(\Spec(-i\dd\rho(x)),\] 
resp., 
\begin{equation}
  \label{eq:spirho} 
s_{\pi}(x) = \sup(\Spec(i\dd\pi(x))) = s_\rho(x) 
\quad \mbox{ for } \quad x \in C.
\begin{footnote}{This is  trivial if $\ad x$ is diagonalizable on $V$ 
on each $\g_\C^n$. Let $\pi_n$ denote the representation of $Z_G(d)_0$ on $\cH_n$ 
and $\pi_n^x(t) := \pi_n(\exp tx)$. 
For the general case, it is instructive to think of the 
spectrum of $\dd\pi_n(x)$ in terms of Arveson's spectral theory, where 
$\Spec(-i\dd\pi_n(x))$ is the minimal closed subset $S \subeq \R$ 
with the property that, 
for every Schwartz function $f \: \R \to \R$ with $\supp(f) \cap S = \eset$, 
we have $\pi_n^x(\hat f) = 0$. In this context it is clear that 
we obtain the same spectrum from the representation on any dense 
invariant subspace, and that equivariant bilinear maps 
are compatible with addition of spectra (cf.\ \cite[Prop.~A.14]{Ne12a}). 
}\end{footnote}
\end{equation}

To see that $C \subeq W_\pi$, 
it now remains to see that $\Ad(G)C$ has interior points. 
Let $U^\pm \subeq \fk^{\pm\phi} = \ker(\id \mp \phi)$ be open convex symmetric $0$-neighborhoods 
for which the map 
\[ E \: U^+ \times U^- \to K, \quad (x_+, x_-) \mapsto 
\exp x_- \exp x_+ \exp x_- \] 
is a diffeomorphism onto an open subset 
of $K$. The existence of such a $0$-neighborhood follows from 
the Inverse Function Theorem because  the differential 
of $E$ in $(0,0)$ is given by $(x_+, x_-) \mapsto x_+ + 2 x_-$. 

For the $\phi^{-1}$-twisted conjugation action 
$c_k^\phi(h) := k h \phi(k)^{-1}$ of $K$ on itself we have 
\[ E(x_+, x_-) = \exp(x_-)\exp(x_+) \exp(x_-) 
= c^\phi_{\exp x_-}(\exp x_+).\] 
Therefore each $c^\phi$-orbit meeting $\im(E)$  also meets
$\exp(U^+)$. 

Next we recall the smooth map $\Hol_{2\pi/N} \: \cL_\phi(\fk) \to K$ 
from Proposition~\ref{prop:2.x} and note that 
\[ V := \{ \xi \in \cL_\phi(\fk) \: \Hol_{2\pi/N}(\xi) \in \im(E)\} \] 
is an open $0$-neighborhood. From Proposition~\ref{prop:2.x} 
we derive that 
every element in $V \times \{1\} \subeq \cL_\phi(\fk) \rtimes \R$ is 
conjugate under $\Ad(\cL_\phi(K))$ to an element $\tilde\xi$ with 
$\Hol_{2\pi/N}(\tilde\xi) \in \exp(U^+)$. 
As $\Hol_{2\pi/N}(x_+) = \exp\big(\frac{2\pi}{N}x_+\big)$ for $x_+ \in \fk^\phi$,  
this further implies that $(\xi,1)$ is conjugate to an element 
in $\frac{N}{2\pi} U^+  \times \{1\}$. We conclude that, 
for every $0$-neighborhood $B \subeq \fk^\phi$, \break  
$\Ad(\cL_\phi(K))(\R \times B \times \{1\})$ contains an open subset 
of the hyperplane  $\R \times \cL_\phi(\fk) \times \{1\} \subeq \hat\cL_\phi(\fk)$. 
Eventually this shows that $\Ad(G)C$ has interior points, and hence that 
$W_\pi\not=\eset$. 
\end{prf} 

At this point we are ready to prove Theorem~\ref{thm:1} stated in the 
introduction. 

\begin{prf} {\bf of Theorem~\ref{thm:1}:} 
We have seen in Theorems~\ref{thm:5.5} and \ref{thm:5.4} that a semibounded 
representation $(\pi, \cH)$ for which $\dd\pi(d)$ is bounded from below 
(which is always the case for $\pi$ or its dual representation) 
is holomorphically induced from a bounded representation 
$(\rho, V)$ of the Hilbert--Lie group $Z_G(d)_0$ 
whose Lie algebra is $\R i c \oplus \fk^\phi \oplus \R i \cdot d$, where 
$\fk^\phi$ is a simple Hilbert--Lie algebra in all $7$ cases. 

Theorem~\ref{thm:classif-hilbert} now provides a classification 
of the bounded irreducible representations of the simply connected 
covering group 
$\R \times \tilde K^\phi \times \R$ of $Z_G(d)_0 
= Z \times K^\phi \times \R$ in terms 
of extremal weights $\lambda \in \cP(\fk^\phi,\ft) \subeq i \ft'$. 
From Lemma~\ref{lem:kphicon} we know that $K^\phi$ is connected, so that its 
simply connected covering group is defined. We even know that $K^\phi$ is 
$1$-connected if $(\Delta_\g)_c$ is not of type $BC_J^{(2)}$, and in the latter case 
$\pi_1(K^\phi) \cong \Z/2$. 
Next Theorem~\ref{thm:5.9} characterizes the weights 
$\lambda$ for which $(\rho_\lambda, V_\lambda)$ is holomorphically inducible 
as the $d$-minimal weights and the corresponding $G$-representation 
$(\pi_\lambda, \cH_\lambda)$ is semibounded by Theorem~\ref{thm:6.1}. 

That the $\ft_\g$-weight set $\cP_\lambda$ of $\pi_\lambda$ satisfies 
\[ \cP_\lambda := \conv(\hat\cW\lambda) \cap (\lambda + \hat\cQ) 
\quad \mbox{ with } \quad 
\Ext(\conv(\cP_\lambda)) = \hat\cW\lambda\]   
follows from the corresponding result in \cite[Thm.~4.10]{Ne10a} for the 
highest weight module $L(\lambda)$ of $\cL_\phi(\fk)_\C^{\rm alg}$. 
This description implies in particular that the equivalence of 
$\pi_\lambda$ and $\pi_\mu$ implies that $\mu \in \hat\cW\lambda$. 
We also see that the set of weights occurring as extremal weights 
in this context is contained in set 
 $\cP^+ = \hat\cW \cP_d^+$  of integral weights bounded from below 
(Remark~\ref{rem:4.5}). 
From \cite[Thm.~3.5]{HN12} we further know that
\[\hat\cW \lambda \cap \cP^+_d = \hat\cW_d \lambda 
= \cW \lambda 
\quad \mbox{ for } \quad \lambda \in \cP^\pm_d,\] 
where $\hat\cW_d \cong \cW$ is the stabilizer of $d$ in $\hat\cW$. 
This leads to a bijection
\[ \cP^\pm_d/\cW \to  \cP^\pm/\hat\cW, \quad 
\cW \lambda \mapsto \hat\cW \lambda\] 

To complete the proof, it only remains to show that 
every elements $\lambda \in \cP_d^+$ actually  is an extremal 
weight of a bounded representation $(\rho_\lambda, V_\lambda)$ 
of $Z_G(d)_0 \cong Z \times K^\phi \times \R$.
As $\lambda\res_{\ft}$ is a weight for $\Delta_0$, 
the existence of the corresponding unitary representation of 
$K^\phi$ follows from Theorem~\ref{thm:classif-hilbert}. 
It therefore remains to verify that 
$N \lambda_c = \lambda(N c) \in \Z$ (Theorem~\ref{thm:cent-ext}). 
In the untwisted cases the normalization of the scalar product is such that 
long roots $\alpha$ satisfy $(\alpha, \alpha) = 2$, so that 
$\lambda_c \in \Z$ follows from \eqref{eq:gweight}. 
In the twisted cases $2\lambda_c \in \Z$ follows from Example~\ref{ex:affineweights}.
\end{prf}

\begin{rem} Let us take a closer look at the ambiguities arising 
in our parametrization of irreducible semibounded unitary representations of 
$G = \hat\cL_\phi(K)$ in terms of bounded representations $(\rho,V)$ of 
$Z_G(d)_0 \cong \T \times K^\phi \times \R$. 

If $(\rho,V)$ is $d$-minimal, the corresponding representation 
of $G$ is obtained by holomorphic induction 
with $\fq = \fp^+ \rtimes \fh_\C$ 
and $V = (\cH^\infty)^{\fp^-}$ is the minimal eigenspace of $\dd\pi(d)$. 
If $(\rho,V)$ is $d$-maximal, then $(\pi, \cH)$ 
is obtained by holomorphic induction 
with $\fq = \fp^- \rtimes \fh_\C$, and 
$V = (\cH^\infty)^{\fp^+}$ is the maximal eigenspace of $\dd\pi(d)$. 

Therefore the only ambiguity in the parametrization 
of corresponding irreducible unitary representations of 
$\hat\cL_\phi(K)$ arises for representations for which 
$\dd\pi(d)$ is bounded, which only happens for one-dimensional representations, 
see Proposition~\ref{prop:dpbd} below. Hence the ambiguity of the parametrization 
consists only in twisting with characters of $\hat\cL_\phi(K)$, 
resp., representations vanishing on the codimension $1$ subgroup 
$\tilde\cL_\phi(K)$ of $\hat\cL_\phi(K)$. 

On the level of $d$-minimal/maximal weights, the corresponding 
assertion is that a weight $\lambda \in i \ft_\g'$ is $d$-minimal 
and $d$-maximal at the same time if and only if 
$n \lambda((\alpha,n)\,\check{}\, ) =  0$ holds for 
every root $(\alpha,n) \in (\Delta_\g)_c$, but this implies 
that the corresponding representation of $G$ is one-dimensional. 
\end{rem}

\begin{prop} \mlabel{prop:dpbd} 
If $(\pi, \cH)$ is an irreducible semibounded representation 
of $\hat\cL_\phi(K)$ for which $\dd\pi(d)$ is bounded, 
then it is one-dimensional. 
\end{prop}

\begin{prf} As $i\cdot d \in W_\pi \cup -W_\pi$, the boundedness of 
$\dd\pi(d)$ implies that $0 \in W_\pi$ and hence that $\pi$ is bounded. 
From Theorem~\ref{thm:centriv} we  now obtain 
that $\dd\pi(c) =0$, so that we obtain a positive 
energy representation of the semidirect product $\cL_\phi(K) \rtimes \R$. 
From Theorem~\ref{thm:semdirtriv} we now derive that 
$\tilde\cL_\phi(K) \subeq \ker \pi$, so that the image of $\pi$ 
is an abelian group and the assertion follows from Schur's Lemma.
\end{prf}

\section{Perspectives and open problems} 
\mlabel{sec:7}

\begin{prob} Let $\fk$ be a simple Hilbert--Lie algebra 
and $\phi \in \Aut(\fk)$ be an automorphism of finite order. 
Then 
\[ \cL_\phi(\fk)\cong   \cL_\phi^1(\fk) 
= \big\{ f \in C^\infty(\R,\fk) \: (\forall t \in \R)\ 
f(t + 1)= \phi^{-1}(f(t))\big\} \] 
can be identified with the space of smooth sections 
of the Lie algebra bundle $\bL_\phi \to \bS^1 \cong \R/\Z$ 
obtained as the quotient of the trivial bundle $\R \times \fk$ 
by the equivalence relation generated by $(t+1,x) \sim  (t, \phi(x))$. 
The Lie connections on this bundle lead to covariant derivatives of the form 
\[ D_A \xi = \xi' + A \xi \quad \mbox{ with } \quad 
A \in \cL_{c_\phi}^1(\der(\fk)), \quad c_\phi(B) = \phi^{-1}B \phi, \] 
and these operators are continuous derivations on 
$\cL_\phi^1(\fk)$ which are skew-symmetric with respect to our 
scalar product, so that they can all be used to define double extensions.  

A map of the form 
$\Gamma(\xi)(t) = \gamma(t)\xi(t)$ with $\gamma \in C^\infty(\R,\Aut(\fk))$ 
defines an isomorphism $\cL_\phi^1(\fk) \to \cL_\psi^1(\fk)$ 
if and only if 
\[ \gamma(t+1) = \psi^{-1} \gamma(t) \phi \quad \mbox{ for } \quad t \in \R.\] 
Let $\gamma \: \R \to\Aut(\fk)$ be the unique smooth curve 
with $\delta^l(\gamma) = A$ with $\gamma(0) = \1$. Then $A \in \cL_{c_\phi}(\der(\fk))$ 
implies that 
\[ \gamma(t+1) = \gamma(1) \phi^{-1} \gamma(t) \phi,\] 
so that we obtain for $\psi := \phi \gamma(1)^{-1}$ an isomorphism 
$\Gamma \: \cL_\phi^1(\fk) \to \cL_\psi^1(\fk)$ satisfying 
\[ D_0 \circ \Gamma = \Gamma \circ D_A.\] 
This means that, by changing the automorphism, we can transform 
the covariant derivative $D_A$ into the standard one. This has the 
advantage that the corresponding one-parameter group 
$(\alpha^A_t)_{t \in \R}$ of automorphisms satisfies 
\[ \alpha_t \circ \Gamma = \Gamma \circ \alpha^A_t \quad \mbox{ for } 
\quad t \in \R.\] 
We conclude in particular, that $\alpha^A$ is periodic if and only 
if the translation action on $\cL_\psi^1(\fk)$ is periodic, which,  
in view of $\alpha_1 \xi = \psi^{-1} \xi$, is equivalent to the order 
of $\psi$ being finite. This gives a geometric interpretation 
for a preference of finite order automorphisms for the contructions 
of double extensions.

It remains to be explored how the representation theory 
of $\hat\cL_\phi(\fk)$ 
changes for other finite order (or even general) automorphisms of $\fk$. 
Is it possible to classify the semibounded 
representations of $\hat\cL_\phi(\fk)$ for any automorphism 
$\phi$ of finite order? 
The present paper covers the case $\phi = \id$ and the 
three involutions which lead to the three twisted locally affine root systems. 
\end{prob} 

\begin{prob} \mlabel{prob:2} The proof of Theorem~\ref{thm:5.9} shows that, 
for every $d$-minimal integral weight $\lambda \in \ft_\g^*$ (continuous or not), 
we have a unitary extremal weight representation $(\pi_\lambda, L(\lambda))$ of $\g_\C$ 
generated by a vector $v_\lambda$ annihilated by $\fp^-$. 
Then the representation $(\rho_\lambda, V_\lambda)$ on the minimal $d$-eigenspace 
$V_\lambda := \ker(\pi_\lambda(d) - \lambda(d)\1)$ is an extremal weight 
representation of the Lie algebra $Z_\g(d)$ and for the orthogonal projection 
$p_V \: L(\lambda) \to V_\lambda$, we obtain a positive definite linear map 
\[ \beta \: U(\g_\C) \to \End(V_\lambda), \quad D \mapsto p_V \pi_\lambda(D) p_V^*.\] 
However, if $\lambda$ is not continuous, then all these representations need not 
integrate to representations of $\hat\cL_\phi(K)$, resp., $Z_G(d)$. 

To deal with the global aspects of these representations, 
we need to pass from $Z_G(d)_0$ to a suitable central extension 
to integrate the representation $\rho_\lambda$ on the completion 
of $V_\lambda$. As the so obtained representation will not be bounded, we need a 
further refinement of the method of holomorphic induction to 
derive a corresponding unitary representation of $\hat\cL_\phi(K)$, or a suitable 
modification of this group, on the completion $\cH_\lambda$ of $L(\lambda)$ 
(cf.\ \cite{Ne11b}). A natural, but certainly not maximal, candidate for a group 
to which these representations may integrate is 
$\hat\cL_\phi(\U_1(\cH))$ if $K = \U_2(\cH)$. 
\end{prob}

\begin{prob} Suppose that $K$ is a $1$-connected simple Hilbert--Lie group. 
Is every irreducible positive energy representation of 
$G = \hat\cL_\phi(K)$ holomorphically induced? Using similar arguments 
as for semibounded representations (cf.\ Theorem~\ref{thm:5.4}), we obtain a representation 
$(\rho,V)$ of $H = Z_G(d)$ on the minimal eigenspace 
$V \not= \{0\}$ for 
$\dd\pi(d)$, but a priori we do not know if this representation 
is bounded, so that holomorphic induction of $(\rho,V)$ 
need not make sense. 

In this context it would be interesting if there are 
(irreducible) unitary representations of $H$ which are 
$d$-minimal in a suitable sense. These representations 
would be natural candidates for a ``holomorphic induction'' 
to a unitary representation of $G$ to make sense. The representations 
from Problem~\ref{prob:2} may lead to interesting examples. 
\end{prob}

\begin{prob} The group $\Diff(\bS^1)$ acts naturally by automorphisms 
on the group $\tilde\cL(K)$. Does it also act on the 
irreducible semibounded representations $(\pi_\lambda, \cH_\lambda)$? 
We expect a unitary representation of the Virasoro group 
which is positive/negative energy representation because we already 
have the action of the generator of the subgroup of rigid rotations. 

In this context it is important to observe that 
the restrictions $\tilde\pi_\lambda$ of the representations 
$\pi_\lambda$ to the codimension-one subgroup 
$\tilde\cL_\phi(K)$ remain irreducible because they are holomorphically 
induced from a bounded irreducible representation 
(cf.\ Remark~\ref{rem:5.9}).  

The philosophy is that the set $\{ [\tilde\pi_\lambda] 
\: \lambda \in \cP_d \}$ 
of equivalence classes of irreducible unitary representations 
of $\tilde\cL(K)$ should be ``discrete'' and therefore fixed 
pointwise under the action of $\Diff(\bS^1)_0$. One way to verify 
this is to observe that $\Diff(\bS^1)_0$ preserves the class of 
those representations which are ``semibounded'' in the sense that 
they extend to semibounded representations of a semidirect product 
with a compact circle subgroup of $\Diff(\bS^1)$. Then one can try 
to show that such representations are determined by their momentum 
sets, but here one looses information by restricting to 
$\tilde\cL_\phi(\fk)$ on which the representation $\pi_\lambda$ is not 
semibounded. 
\end{prob}

\appendix

\section{$H^1$-curves in Hilbert--Lie groups}  \mlabel{app:a}

In this appendix we briefly introduce the Banach--Lie algebra 
of $H^1$-curves with values in a Hilbert--Lie algebra and explain how this 
can be used to obtain Banach--Lie algebras 
$\tilde\cL^H_\phi(\fk)$ whose construction is based on $H^1$-curves instead 
of smooth ones. We also obtain corresponding Banach--Lie groups 
$\tilde\cL^H_\phi(K)$ in which the groups $\tilde \cL_\phi(K)$ are dense. 

\subsection{The group of $H^1$-curves} 

\begin{defn}
Let $I = [0,1] \subeq \R$ denote the unit interval. 
We write $H^1(I)$ for the space of absolutely continuous 
functions $f \: I \to \R$ with $f' \in L^2(I)$, endowed with the scalar product 
$$ \la f,g \ra := \int_0^1 f(t)g(t) + f'(t)g'(t)\, dt. $$
We recall from \cite[Cor.~9.7]{Pa68} that $H^1(I)$ is a Hilbert space, 
that the inclusion $H^1(I) \to C(I,\R)$ is continuous, and that 
$H^1(I)$ is a Banach algebra with respect to the pointwise product. 

If $\cH$ is a Hilbert space, then we write $H^1(I,\cH) := H^1(I) \hat\otimes \cH$ 
for the tensor product of Hilbert spaces. 
\end{defn}

\begin{rem}  \mlabel{rem:6.1} 
(a) Let $(e_i)_{i \in I}$ be an orthonormal basis of $\cH$ and 
$f \in H^1(I,\cH)$. Then $f = \sum_{i \in I} f_i e_i$ with 
$f_i \in H^1(I)$ satisfying 
$\|f\|^2 = \sum_i \|f_i\|^2 < \infty$. This implies that for 
each $t \in I$ we have 
$\sum_{i \in I} |f_i(t)|^2 < \infty,$
so that we obtain a well-defined function 
$$ f \: I \to \cH, \quad f(t) := \sum_{i \in I} f_i(t)e_i.$$ 
The sum on the right hand side is actually countable, so that we 
have a series expansion of $f$, where each finite sum is an 
$H^1$-function with values in some finite-dimensional subspace 
and the range of $f$ lies in a separable subspace. 
From the Dominated Convergence Theorem we now derive that 
$f$ is continuous. One can further show that $f$ is absolutely 
continuous and that $f' \: I \to \cH$ exists almost 
everywhere in such a way that the Fundamental Theorem holds 
(cf.\ \cite[Sect.~25]{Wl82} for a detailed treatment of 
$H^1$-spaces with values in a (separable) Hilbert space). 

(b) If, in addition, $\cH$ carries a continuous bilinear product 
$\cH \times \cH \to \cH$, then the product rule 
implies that $H^1(I,\cH) \subeq C(I,\cH)$ is a subalgebra. 
Now \cite[Lemma~A.2]{Ne02a} implies that the multiplication 
on $H^1(I,\cH)$ is continuous, hence turning it into a 
Banach algebra.
\end{rem}

Now let $K$ be a Hilbert--Lie group with Lie algebra $\fk$. Then 
$C(I,K)$ carries the structure of a Banach--Lie group 
with Lie algebra $C(I,\fk)$ and the inclusion 
$H^1(I,\fk) \into C(I,\fk)$
is a morphism of Banach--Lie algebras (cf.~Remark~\ref{rem:6.1}). 
We write 
$$ H^1(I,K) \subeq C(I,K) $$
for the corresponding integral subgroup of $C(I,K)$, so that 
$H^1(I,K)$ is a Banach--Lie group with Lie algebra $H^1(I,\fk)$ 
whose underlying set is the subgroup \break 
$\la \exp H^1(I,\fk) \ra 
\subeq C(I,K)$ (cf.\ \cite[Sect.~I.9]{Alb93}, \cite{Mais62}). 

Then $H^1(I,K)$ consists of paths $\gamma \: I \to K$ for which 
the left and right logarithmic derivative exists almost everywhere
and 
$\delta^r(\gamma), \delta^l(\gamma)\: I \to \fk$
are $L^2$-functions (it suffices to verify this for the 
exponential image of $H^1(I,\fk)$).

\begin{prop} \mlabel{prop:logder} For a Hilbert--Lie group $K$, 
the following assertions hold: 
\begin{description}
\item[\rm(i)] $\Ad_{L^2} \: H^1(I,K) \to \OO(L^2(I,\fk)), 
\Ad_{L^2}(f)(\xi)(t) = \Ad(f(t))\xi(t)$ defines a bounded 
representation of the Lie group $H^1(I,K)$ 
on the real Hilbert space $L^2(I,\fk)$. 
\item[\rm(ii)] The right logarithmic derivative 
$\delta^r \: H^1(I,K) \to L^2(I,\fk), f \mapsto f' \cdot f^{-1}$ 
is a smooth cocycle whose differential is the Lie algebra cocycle 
$\L(\delta^r)f = f'.$
\end{description}
\end{prop}

\begin{prf} (i) Since the multiplication map 
$H^1(I) \times L^2(I) \to L^2(I), (f,g) \mapsto fg$ 
is continuous, the Lie bracket 
induces a continuous bilinear map 
\begin{equation}
  \label{eq:h1brack}
 H^1(I,\fk) \times L^2(I,\fk) \to L^2(I,\fk), \quad 
(\xi,\eta) \mapsto [\xi,\eta], 
\end{equation}
defining a continuous representation of the Banach--Lie algebra 
$H^1(I,\fk)$ on the Hilbert space $L^2(I,\fk)$. This representation 
integrates to the morphism $\Ad_{L^2}$ of Banach--Lie groups. 

(ii) First we observe that the cocycle property follows from 
the product rule 
$$ \delta^r(fg) = \delta^r(f) + \Ad(f) \delta^r(g). $$
Since $\delta^r$ is a cocycle with values in the smooth $H^1(I,K)$-module 
$L^2(I,\fk)$, it defines a homomorphism 
$$ (\delta^r,\id) \: H^1(I,K) \to L^2(I,\fk)\rtimes H^1(I,K) $$
of Banach--Lie groups. Therefore its smoothness follows, once we 
have shown that it is continuous. 
As $\delta^r$ is a cocycle, it suffices to verify its continuity 
in an identity neighborhood, so that it suffices to show that the map 
$$ \delta^r \circ \exp_{H^1(I,K)} \: H^1(I,\fk) \to L^2(I,\fk), \quad 
f \mapsto \delta^r(\exp_K \circ f) $$
is continuous. 
Writing $\kappa^r_K$ for the right Maurer--Cartan form on $K$, 
we find 
$$ \delta^r(\exp_K \circ f) 
= (\exp_K \circ f)^*\kappa^r_K 
= f^*(\exp_K^*\kappa^r_K). $$
The $1$-form $\kappa_\fk := \exp_K^*\kappa^r_K \in \Omega^1(\fk,\fk)
\cong C^\infty(\fk,B(\fk))$ is explicitly given by the analytic function  
$$  F \: \fk \to B(\fk), \quad F(x) 
:= \frac{\1-e^{-\ad x}}{\ad x} 
= \sum_{n = 0}^\infty \frac{(-1)^n}{(n+1)!} (\ad x)^n,$$
and we have 
$$ \delta^r(\exp_K \circ f)(t)  
= F(f(t))(f'(t)). $$
The evaluation map $B(\fk) \times \fk \to \fk$ is continuous,  
it induces a continuous bilinear map 
$C(I,B(\fk)) \times L^2(I,\fk) \to L^2(I,\fk)$. 
Further, the map $H^1(I,\fk) \to L^2(I,\fk),  f \mapsto f',$ 
and the inclusion $H^1(I,\fk) \to C(I,\fk)$ are continuous. 
Therefore it remains to observe that the map 
$$ C(I,\fk) \to C(I,B(\fk)), \quad f \mapsto F \circ f $$
is continuous, 
because for each Banach space $X$, the topology on the space $C(I,X)$ 
defined by the sup-norm coincides with the compact open topology. 
This completes the proof of the smoothness of $\delta^r$. 

To calculate its derivative in $\1$, we note that for $s \in \R^\times$, 
we have 
$$ \frac{1}{s} \delta^r(\exp_K \circ (sf))
= \big(F\circ (s\cdot f)\big)(f'). $$
Therefore 
$\lim_{s \to 0} F\circ (s\cdot f) = F(0) = \id_\fk$
implies that $\L(\delta^r)f := T_\1(\delta^r)f = f'$. 
\end{prf}

\begin{lem} \mlabel{lem:trans} 
If a group $G$ acts by isometries on the metric space 
$(X,d)$, then each open $G$-orbit is also closed. 
In particular, the action is transitive if $X$ is connected 
and $G$ has an open orbit. 
\end{lem}

\begin{prf} (cf.\ \cite{Go03}).
Let $\cO = G x_0$ be an open orbit and suppose that the ball 
$B_\eps(x_0)$ is contained in $\cO$. 
We show that $\cO$ is also closed. 
Let $y \in \oline{\cal O}$. Then $B_\eps(y)$ intersects ${\cal O}$ 
in some point $gx_0$. Then 
$y \in B_\eps(gx_0) = gB_\eps(x_0) \subeq {\cal O}$ 
shows that $\cO$ is closed. 
\end{prf}

The following proposition is well known for the case where 
$K$ is a compact group (cf.\ \cite[p.~23]{Te89}). 
  
\begin{prop} \mlabel{prop:aff-act} The 
affine action of the normal subgroup 
\[ H^1(I,K)_* :=  \{ f \in H^1(I,K) \: f(0) = \1\} \]
of  $H^1(I,K)$ on $L^2(I,\fk)$ by 
$$ \tau_f(\xi) := \Ad(f)\xi - \delta^r(f) $$
is simply transitive. Each orbit map yields a diffeomorphism 
$H^1(I,K)_* \to L^2(I,\fk)$. In particular, 
$H^1(I,K)_*$ is contractible. 
\end{prop} 

\begin{prf} That $\tau$ defines a smooth group action  follows from the 
cocycle property and the smoothness of $\delta^r$ 
(Proposition~\ref{prop:logder}). Moreover, this 
action is isometric. Further, 
the derivative in $\1$ of the orbit map $\tau^0$ of $0$ is 
$$ H^1(I,\fk)_* \to L^2(I,\fk), \quad f \mapsto -f', $$
which is a topological linear isomorphism of Hilbert spaces. 
It follows from the Inverse Function Theorem 
that the orbit $\cO_0$ of $0$ is open 
and Lemma~\ref{lem:trans} implies that $\cO_0 = L^2(I,\fk)$. 

Since $\delta^r(f) = 0$ implies that $f$ is constant, the 
stabilizer of $0$ in $H^1(I,K)_*$ is trivial and the orbit map 
$$ \tau^0 \: H^1(I,K)_* \to L^2(\fk), \quad f \mapsto - \delta^r(f) $$
is a smooth equivariant bijection. Since its differential in $\1$ 
is a topological isomorphism, the equivariance implies that this is 
everywhere the case, and finally the Inverse Function Theorem  
shows that $\tau^0$ is a diffeomorphism. 
\end{prf}

\subsection{The $H^1$-version of twisted loop groups} 

To apply the method of holomorphic induction (cf.\ Appendix~\ref{app:c}) 
to the group $\hat\cL_\phi(K)$ constructed in Section~\ref{sec:3}, 
we need a Banach version of this group. Since we shall see that all 
semibounded representations of $\hat\cL_\phi(K)$ extend to various 
Banach completions of this  group (Remark~\ref{rem:5.9}), it makes 
sense to use one which is rather large. 

To this effect, we observe that 
the cocycle $\omega_D(\xi,\eta) = \la \xi',\eta\ra$ on $\cL_\phi(\fk)$ extends continuously 
to the Banach--Lie algebra $\cL_\phi^H(\fk)$ of twisted loops of class $H^1$, 
so that we obtain a central extension 
$\tilde\cL^H_\phi(\fk)$ which again is a Banach--Lie algebra. 
Below we show that this Lie algebra integrates to a $1$-connected 
Banach--Lie group $\tilde\cL^H_\phi(K)$ on which we have a continuous 
$\R$-action $\alpha$ defined by translations. 

The Lie algebra  $\cL^H(\fk)$ of $H^1$-loops is maximal with the 
property that the cocycle $D\xi := \xi'$ defined by the derivative 
defines a linear functional on $\cL(\fk)$ which is continuous with respect to the 
$L^2$-norm. This is crucial to define a corresponding cocycle by 
$\omega(\xi,\eta) = \la \xi',\eta\ra$. In particular, 
there are no non-trivial cocycles for the Lie algebra 
$C(\bS^1,\fk)$ of continuous loops (cf.\ \cite[Cor.~13, Thm.~16]{Mai02}).

\begin{defn} \mlabel{def:a1.1} For a Hilbert--Lie algebra $\fk$ and an automorphism 
$\phi \in \Aut(\fk)$ of order $N$, we write 
$\cL^H_\phi(\fk)$ for the Hilbert space of 
local $H^1$-maps $f \: \R \to \fk$ satisfying the 
condition 
\[ (\forall t \in \R)\ f\Big(t + \frac{2\pi}{N}\Big) = \phi^{-1}(f(t)), \] 
endowed with the Hilbert norm defined by 
\[ \|\xi\|_{H^1}^2 := \|\xi\|_2^2 + \|\xi'\|_2^2 
:=  \frac{1}{2\pi} \int_0^{2\pi} \|\xi(t)\|^2 + \|\xi'(t)\|^2\, dt.\] 
This defines on $\cL^H_\phi(\fk)$ the structure of a Banach--Lie algebra. 
It is NOT a Hilbert--Lie algebra in the sense of Definition~\ref{def:1.1} 
because the norm is not invariant under the adjoint action. 
Since the derivative 
defines a continuous map from $H^1$ to $L^2$, 
\[ \omega_D(\xi,\eta) 
= \la \xi',\eta \ra
= \frac{1}{2\pi}\int_0^{2\pi} \la \xi'(t), \eta(t)\ra\, dt \] 
defines a continuous $2$-cocycle on $\cL_\phi^H(\fk)$, and we thus obtain 
the centrally extended Banach--Lie algebra 
\[ \tilde\cL_\phi^H(\fk) := \R \oplus_{\omega_D} \cL_\phi^H(\fk), \quad 
[(z,\xi), (w,\eta)] := (\omega_D(\xi,\eta), [\xi, \eta]), \] 
containing the Fr\'echet--Lie algebra $\tilde\cL_\phi(\fk) = \R \oplus_{\omega_D} \cL_\phi(\fk)$. 
\end{defn}

\begin{rem} \mlabel{rem:homot-h1} Let $K$ be a connected Hilbert--Lie group 
and $\phi \in \Aut(K)$ be an automorphism of order $N$. 
To obtain similar information as in Remark~\ref{rem:homot} 
on the topology of the Banach--Lie group 
\begin{align*}
\cL^H_\phi(K)
&:= \Big\{ f \in H^1_{\rm loc}(\R,K) \: (\forall t \in \R)\ 
f\Big(t + \frac{2\pi}{N}\Big) = \phi^{-1}(f(t))\Big\} \\
&\cong \Big\{ f \in H^1([0,2\pi/N], K) \: 
f(2\pi/N) = \phi^{-1}(f(0))\Big\},
\end{align*}
we first claim that the inclusion 
$\cL(K) \to \cL^H(K)$ of untwisted loop groups induces isomorphisms
\[ \pi_k(\cL(K)) \to \pi_k(\cL^H(K)) \quad \mbox{ for } \quad k \in \N_0.\] 
Consider the commutative diagram 
\[ \matr{ 
\cL_\phi(K)_* & \to & \cL_\phi(K) & \smapright{\ev_0} & K \\ 
\mapdown{} & & \mapdown{}   &  & \mapdown{\id_K} \\ 
\cL^H_\phi(K)_* & \to & \cL^H_\phi(K) & \smapright{\ev_0} & K}\] 
in which both rows describe locally trivial fiber bundles. 
Let $I = [0,a]$ for $a := \frac{2\pi}{N}$ and 
$\Omega(I,K) \subeq H^1_*(I,K)$ denote the kernel 
of the evaluation map 
\[ \ev_a \: H^1_*(I,K) \to K, \quad f \mapsto f(a) \] 
in $a$. Then $\ev_a$ defines a locally trivial fiber bundle, 
so that the contractibility of $H^1_*(I,K)$ (Proposition~\ref{prop:aff-act}) 
implies the existence of natural isomorphisms 
\[ \pi_{k+1}(K) \to \pi_k(\Omega(I,K)), \quad k \in \N_0.\] 
Next we observe that restriction to  
$[0,\frac{2\pi}{N}]$ defines an isomorphism 
\[ \cL^H_\phi(K)_* \to \Omega([0,a],K).\] 
Since we also have natural isomorphisms 
\[ \pi_{k+1}(K) \to \pi_k(\cL_\phi(K)_*), \quad k \in \N_0\] 
(cf.\ \cite[Cor.~3.4]{NeWo09}), we conclude that the inclusion 
$\cL_\phi(K)_* \to \cL^H_\phi(K)_*$ 
induces isomorphisms of all homotopy groups. 
Applying the Five Lemma to the long exact homotopy sequence corresponding 
to the rows of the above diagram, we see that the inclusion 
$\cL_\phi(K) \to \cL_\phi^H(K)$ also induces isomorphisms of all 
homotopy groups (cf.\ \cite{Ne02a} for more details on this technique). 
\end{rem}

\begin{thm} \mlabel{thm:cent-ext-h1} 
The assertion of  Theorem~\ref{thm:cent-ext} remains true for 
the Banach--Lie algebra 
$\tilde\cL_\phi^H(\fk)$ defined by $H^1$-maps and the corresponding group 
$\tilde\cL_\phi^H(K)$. 
\end{thm}

\begin{prf}  From Remark~\ref{rem:homot-h1} it follows that 
the period homomorphism 
$\per_{\omega_D} \: \pi_2(\cL_\phi^H(K)) \to \R$ has the same range as the 
period homomorphism on $\pi_2(\cL_\phi(K))$, and  since 
$\cL^H_\phi(K)$ is also $1$-connected (Remark~\ref{rem:homot-h1}), 
\cite[Thm.~7.9]{Ne02b} applies as in the proof of 
Theorem~\ref{thm:cent-ext} the existence of a central $\T$-extension 
$\tilde\cL^H_\phi(K)$ of $\cL^H_\phi(K)$ which is compatible 
with the inclusion $\cL_\phi(K) \into \cL_\phi^H(K)$. 
\end{prf}

\begin{defn}\mlabel{def:3.5}
It is easy to see that the rotation action of $\R$ on 
$\cL_\phi^H(K)$ lifts uniquely to a continuous action on the central 
extension $\tilde\cL_\phi^H(K)$ (\cite[Thm.~V.9]{MN03}), 
but since the rotation action on 
$\cL^H_\phi(K)$ is not differentiable, the corresponding semidirect product group 
\[ \hat\cL_\phi^H(K) := \tilde\cL_\phi^H(K) \rtimes \R \] 
is a topological group but not a Lie group. This is the main difference 
to the smooth setting, where 
$\hat\cL_\phi(K) = \tilde\cL_\phi(K) \rtimes \R$
is a Fr\'echet--Lie group. 
\end{defn}

\begin{rem} \mlabel{rem:polar-h1} 
As in Remark~\ref{rem:polar}, we derive from the polar decomposition 
$K_\C = K \exp(i \fk)$ of the universal complexification of $K$ 
the existence of a central extension of complex 
Lie groups 
\[ \1 \to \C^\times \to \tilde\cL_\phi^H(K_\C) \to \cL_\phi^H(K_\C) \to \1 \] 
for which the inclusion 
$\tilde\cL_\phi^H(K) \into \tilde\cL_\phi^H(K_\C)$ is a universal complexification 
and a weak homotopy equivalence. 
\end{rem}

\section{Analytic operator-valued positive definite functions} 
\mlabel{app:b}

In this appendix we discuss operator-valued positive 
definite functions on Lie groups. The main result 
is Theorem~\ref{thm:charaposdef}, asserting that, for a Hilbert space $V$,  analytic 
$B(V)$-valued defined in a $\1$-neighborhood of a Fr\'echet--BCH 
Lie group $G$ are positive definite if the corresponding linear map 
$\beta \: U(\g_\C) \to B(V)$, defined by derivatives in $\1$, is positive definite. 

\begin{defn} \mlabel{def:a.1} 
Let $X$ be a set and $\cK$ be a Hilbert space.

\par (a)  A function $Q \: X \times X \to B(\cK)$ is called a 
{\it $B(\cK)$-valued kernel}. 
It is said to be {\it hermitian} if 
$Q(z,w)^* = Q(w,z)$ holds for all $z, w \in X$. 

\par (b) A hermitian $B(\cK)$-valued kernel $K$ on $X$ is said to be 
{\it positive definite} if 
for every finite sequence $(x_1, v_1), \ldots, (x_n,v_n)$ in $X \times \cK$
we have 
\[ \sum_{j,k = 1}^n \la Q(x_j, x_k)v_k, v_j \ra \geq 0. \] 

\par (c) If $(S,*)$ is an involutive semigroup, then a 
function $\phi \: S \to B(\cK)$ is called {\it positive definite} 
if the kernel $Q_\phi(s,t) := \phi(st^*)$ is positive definite. 

\par (d) Positive definite kernels can be characterized as those 
for which there exists a Hilbert space $\cH$ and a 
function $\gamma \: X \to B(\cH,\cK)$ such that 
\[ Q(x,y) = \gamma(x)\gamma(y)^* \quad \mbox{ for } \quad x,y \in X \] 
(cf.\ \cite[Thm.~I.1.4]{Ne00}). 
Here one may assume that the vectors 
$\gamma(x)^*v$, $x \in X, v \in \cK$, span a dense subspace of 
$\cH$. Then the pair $(\gamma,\cH)$ is called a {\it realization 
of $K$}. 
The map $\Phi \: \cH \to \cK^X, \Phi(v)(x) := \gamma(x)v$, 
then realizes $\cH$ as a Hilbert subspace of $\cK^X$ 
with continuous point evaluations $\ev_x \: \cH \to \cK$. 
It is the unique Hilbert subspace in $\cK^X$ with this property 
for which $Q(x,y) = \ev_x \ev_y^*$ for $x,y \in X$. 
We write $\cH_Q \subeq \cK^X$ for this subspace and call it 
the {\it reproducing kernel Hilbert space with kernel~$Q$}.
\end{defn}

\begin{defn} Let $\cK$ be a Hilbert space, $G$ be a group, and 
$U \subeq G$ be a subset. A function $\phi \: UU^{-1} \to B(\cK)$ 
is said to be {\it positive definite} if the kernel 
\[ Q_\phi \: U \times U \to B(\cK), \quad (x,y) \mapsto \phi(xy^{-1})\] 
is positive definite. 
\end{defn}

\begin{defn} A Lie group $G$ with Lie algebra $\g$ is said to be 
{\it locally exponential} 
if it has an exponential function for which there is an open $0$-neighborhood 
$U$ in $\g$ mapped diffeomorphically by $\exp_G$ onto an 
open subset of $G$. If, in addition, $G$ is analytic and 
the exponential function is an analytic local diffeomorphism in $0$, 
then $G$ is called a {\it BCH--Lie group} (for Baker--Campbell--Hausdorff). 
Then the Lie algebra $\g$ is a 
{\it BCH--Lie algebra}, i.e., there exists an open 
$0$-neighborhood $U \subeq \g$ such that for $x,y \in U$ the Hausdorff series 
$$x * y = x + y + \frac{1}{2}[x,y] + \cdots  $$
converges and defines an analytic function 
$U \times U \to \g, (x,y) \mapsto x * y$. 
The class of BCH--Lie groups contains in particular all Banach--Lie groups 
(\cite[Prop.~IV.1.2]{Ne06}). 
\end{defn}

\begin{thm} {\rm(Extension of local positive definite analytic functions)} 
\mlabel{thm:extension} (cf.\ \cite[Thm.~A.7]{Ne11c})
Let $G$ be a $1$-connected Fr\'echet--BCH--Lie group, 
$V \subeq G$ an open connected $\1$-neighborhood, 
$\cK$ be a Hilbert space and 
$\phi \: VV^{-1} \to B(\cK)$ be an analytic positive definite function.  
Then there exists a unique analytic positive definite 
function $\tilde\phi \: G \to B(\cK)$ extending~$\phi$. 
\end{thm}

\begin{defn} Let $U$ be an open subset of the Lie group $G$ and 
$E$ be a locally convex space. Then we obtain for each 
$x \in \g$ a differential operator on $C^\infty(U,E)$ by 
\[ (L_x f)(g) := \derat0 f(g \exp tx).\] 
These operators define a representation of the Lie algebra $\g$ on 
$C^\infty(U,E)$, so that we obtain a natural extension to a homomorphism 
\[ U(\g) \to \End(C^\infty(U,E)), \quad D \mapsto L_D. \] 
We likewise define 
\[ (R_x f)(g) := \derat0 f(\exp(tx) g)\] 
and note that $[R_x, R_y] = R_{[y,x]}$ for $x,y \in \g$. 
\end{defn}

\begin{thm} {\rm(Infinitesimal characterization of positive definite analytic functions)} 
\mlabel{thm:charaposdef} 
Let $G$ be a Fr\'echet--BCH--Lie group, 
$V \subeq G$ an open connected $\1$-neighborhood, 
$\cK$ be a Hilbert space and 
$\phi \: V \to B(\cK)$ be an analytic function satisfying 
$\phi(\1)= \1$.  
Then $\phi$ is positive definite on a $\1$-neighborhood in 
$G$ if and only if the corresponding linear map 
\[ \beta \: U(\g_\C) \to B(\cK), \quad \beta(D) := (L_D \phi)(\1) \] 
is a positive definite linear function on the $*$-algebra $U(\g_\C)$. 
\end{thm}

\begin{prf} ``$\Rarrow$'': Suppose first that $\phi$ is positive definite 
in a $\1$-neighborhood. Then Theorem~\ref{thm:extension} 
provides an extension of the germ of $\phi$ in $\1$ to an 
analytic positive definite function on all of $G$. We may 
therefore assume that $\phi$ is defined on $G$. 
Then the vector-valued GNS construction provides a unitary 
representation 
$(\pi, \cH)$ of $G$ on a Hilbert space $\cH$, containing 
$\cK$ as a closed subspace such that the orthogonal 
projection $p_\cK \: \cH \to \cK$ satisfies 
\[ \phi(g) = p_\cK \pi(g) p_\cK^* \quad \mbox{ for } \quad g \in G.\] 
This implies that $\cK$ consists of analytic vectors, and 
for $D \in U(\g_\C)$ we find the formula 
\[ \beta(D)  = p_\cK \dd\pi(D) p_\cK^*.\] 
Any function of this form is easily seen to be positive definite. 

``$\Larrow$'': Let $U_\g \subeq \g$ be an open symmetric 
$0$-neighborhood which is mapped by $\exp$ bianalytically 
to an open $\1$-neighborhood of $G$ and such that 
$\phi$ is defined on $\exp(U_\g)$. Then 
$\phi \circ \exp \: U_\g \to B(V)$ is also analytic, 
and, after shrinking $U_\g$, we may assume that 
\[ \phi(\exp x) = \sum_{n = 0}^\infty \phi_n(x),\] 
where $\phi_n \: \g \to B(V)$ is a continuous homogeneous polynomial 
function of degree $n$ (\cite{BS71}). Now the relation 
\[ \phi(\exp tx) 
= \sum_n \frac{t^n}{n!} (L_x^n \phi)(\1)
= \sum_n \frac{t^n}{n!} \beta(x^n) 
\quad \mbox{ for } \quad 
|t| < \eps \] 
implies that $\phi_n(x) = \frac{\beta(x^n)}{n!}$. 
In particular, 
\[ \phi(\exp x) = \sum_n \frac{1}{n!} \beta(x^n)
\quad \mbox{ for } \quad x \in U_\g, \] 
which implies that $\beta$ is {\it analytic} in the sense 
of \cite[Def.~3.2]{Ne11a}. 

Let $\beta_n(x_1, \ldots, x_n) := \beta(x_1\cdots x_n)$ 
and 
\[  \beta^s_n(x_1,\ldots, x_n) := \frac{1}{n!} \sum_{\sigma \in S_n} 
\beta(x_{\sigma(1)}, \ldots, x_{\sigma(n)}) \] 
be its symmetrization. For a continuous seminorm $p$ on $\g$, we then define 
\[  \|\beta_n^s\|_p := \sup \{ \|\beta_n^s(x_1,\ldots, x_n)\| \: 
x_1, \ldots, x_n \in \g, p(x_i)\leq 1\} 
\in [0,\infty]. \] 
From \cite[Prop.~3.4]{Ne11a} we now obtain the existence of a 
continuous seminorm $p$ on $\g$ with 
$\sum_n \frac{1}{n!} \|\beta_n^s\|_p  < \infty.$
In particular, there exists a constant $C > 0$ with 
\begin{equation}
  \label{eq:esti} 
\|\beta_n^s\|_p \leq C n! \quad \mbox{ for all } \quad n \in \N_0.
\end{equation}

Let 
$\cH \subeq \Hom(U(\g_\C),\cK)$ be the reproducing kernel 
Hilbert space corresponding to the positive definite function 
$\beta \: U(\g_\C) \to B(\cK)$. The corresponding positive definite 
kernel $Q$ and the corresponding evaluation maps $Q_D \: \cH \to \cK$ then satisfy 
\[ Q(D_1, D_2) = \beta(D_1 D_2^*) = Q_{D_1} Q_{D_2}^* 
\quad \mbox{ and } \quad 
Q_D f = f(D) \quad \mbox{ for } \quad f \in \cH. \] 
We have a $*$-representation of $U(\g_\C)$ on the dense subspace 
\[ \cH^0 := \Spann \{ Q_D^*v \: v \in \cK, D \in U(\g_\C)\} \] 
by 
\[ (\rho(D)f)(D') = f(D'D) \quad \mbox{ for }  \quad D,D' \in U(\g_\C).\] 
From  $\beta(\1) = Q_\1 Q_\1^* = \1$ we derive that we may identify 
$\cK$ with its image under the isometric embedding 
$Q_\1^* \: \cK \to \cH$. For $v \in \cK$ we then have 
\[ (\rho(D)v)(D') = v(D'D) = Q_{D'D} Q_\1^* v = \beta(D'D)v 
= Q_{D'} Q_{D^*}^* v = (Q_{D^*}^* v)(D'),\] 
so that 
\[ \rho(D)v = Q_{D^*}^* v \quad \mbox{ for }  \quad  D \in U(\g_\C).\] 

In view of $\|Q_D\|^2 = \|Q_D Q_D^*\| = \|\beta(DD^*)\|$, 
we find for the operators $Q_{x^n} \in B(\cH,\cK)$, $x \in \g$, the estimates 
\[ \frac{1}{n!} \|Q_{x^n}\| 
= \frac{1}{n!} \|\beta(x^{2n})\|^{1/2} 
\leq \frac{p(x)^n}{n!} \|\beta_{2n}^s\|_p^{1/2} 
\leq \frac{p(x)^n}{n!} \sqrt{C} \sqrt{(2n)!}.\] 
In view of $\lim_{n \to \infty} \frac{\sqrt{(2n+2)(2n+1)}}{n+1} = 2$, 
it follows that 
\[\sum_n \frac{1}{n!} \|Q_{x^n}\| < \infty \quad 
\mbox{ for } \quad p(x) < \shalf.\] 
We thus obtain an analytic function 
\[ \eta \: \{ x \in \g \: p(x) < \shalf \} \to B(\cH,\cK), \quad 
\eta(x) := \sum_n \frac{1}{n!} Q_{x^n}.\] 
Now let $W \subeq \{ x \in \g \:p(x) < \shalf\}$ be an open symmetric $0$-neighborhood such that 
all BCH products $x * y$ for $x,y \in W$ 
are defined and that we thus obtain an analytic function 
on $W \times W$ with values in the set $\{ z \in \g \: p(z) < \frac{1}{2}\}$. 
For $x,y \in W$ we finally derive 
\begin{align*}
 \phi(\exp x \exp(-y)) 
&= \phi(\exp (x * (-y))) = \sum_n \frac{1}{n!}\beta((x * (-y))^n)
= \sum_{k,\ell} \frac{1}{k!\ell!}\beta(x^k (-y)^\ell)\\
&= \sum_{k,\ell} \frac{1}{k!\ell!} Q_{x^k} Q_{y^\ell}^* 
= \eta(x) \eta(y)^*. 
\end{align*}
This factorization implies that $\phi$ is positive 
definite on $\exp W \exp W$. This completes the proof. 
\end{prf}

\begin{rem} If $\cK$ is one-dimensional, then a 
linear map $\beta \: U(\g_\C) \to B(\cK) \cong \C$ is 
positive definite if and only if it is a positive functional 
in the sense that $\beta(DD^*) \geq 0$ for every $D \in U(\g_\C)$. 
For general $\cK$, it is shown in \cite[Ex.~11.2.1, Thm.~11.2.2]{Sch90} that 
$\beta$ is positive definite if and only if it is 
{\it completely positive} in the sense that 
every induced map 
\[ M_n(\beta) \: M_n(U(\g_\C)) \to M_n(B(\cK)) \cong B(\cK^n) \] 
obtained by applying $\beta$ to all matrix entries maps positive elements 
to positive elements. 
\end{rem}

\section{Holomorphic induction for BCH--Lie groups} 
\mlabel{app:c} 

Let $G$ be a Lie group and $M = G/H$ be a homogeneous space of $G$ which carries the structure of a 
complex manifold so that $G$ acts analytically by holomorphic maps. 
In \cite{Ne12a} we have developed a theory of holomorphic induction for  bounded 
unitary representations of $H$ in the context where $G$ is a Banach--Lie group. 
To deal with semibounded representations of Fr\'echet--Lie groups such as 
the double extension $\hat\cL_\phi(K)$ of the Fr\'echet--Lie group 
$\cL_\phi(K)$ of smooth $\phi$-twisted loops, 
we need an extension of this theory to certain classes of 
Fr\'echet--Lie groups. In this appendix we explain which 
properties of Banach--Lie groups were used in \cite[Sects.~1,2]{Ne12a} and why 
$\hat\cL_\phi(K)$ also has these properties.

Let $G$ be a connected Fr\'echet--BCH--Lie group with Lie algebra $\g$. We further assume that 
there exists a complex BCH--Lie group $G_\C$ with Lie algebra $\g_\C$ and a natural 
map $\eta \: G \to G_\C$ for which $\L(\eta)$ is the inclusion $\g \into \g_\C$. 
Let $H \subeq G$ be a Lie subgroup for which $M := G/H$ carries the structure of a smooth manifold 
with a smooth $G$-action and $\fh \subeq \g$ be its Lie algebra. 
We also assume the existence of closed $\Ad(H)$-invariant subalgebras 
$\fp^\pm \subeq \g_\C$ with $\oline{\fp^\pm} = \fp^\mp$ for which we have a topological direct 
sum decomposition 
\begin{equation}
  \label{eq:splitcond}
 \g_\C = \fp^+ \oplus \fh_\C \oplus \fp^-.\tag{SC}
\end{equation}
We put 
\[ \fq := \fp^+ \rtimes \fh_\C \quad \mbox{ and } \quad \fp := \g \cap (\fp^+ \oplus \fp^-),\]
so that $\g = \fh \oplus \fp$ is a topological direct sum. We 
assume that there exist open symmetric convex $0$-neighborhoods 
$U_{\g_\C} \subeq \g_\C$, 
$U_\fp \subeq \fp \cap U_{\g_\C}, U_\fh \subeq \fh \cap U_{\g_\C}, U_{\fp^\pm} \subeq\fp^\pm  \cap U_{\g_\C}$ and  
$U_\fq \subeq \fq \cap U_{\g_\C}$ such that the 
BCH-product is defined and holomorphic on 
$U_{\g_\C} \times U_{\g_\C}$, and the following maps are analytic diffeomorphisms onto an open subset: 
\begin{description}
\item[\rm(A1)] $U_{\fp} \times U_\fh \to \g, (x,y) \mapsto x * y$. 
\item[\rm(A2)] $U_{\fp} \times U_\fq \to \g_\C, (x,y) \mapsto x * y$. 
\item[\rm(A3)] $U_{\fp^-} \times U_\fq \to \g_\C, (x,y) \mapsto x * y$. 
\end{description}

Then (A1) implies the existence of a smooth manifold structure on $M = G/H$ for which $G$ acts 
analytically. 
Condition (A2) implies the existence of a complex manifold structure on $M$ which is $G$-invariant 
and for which $T_{\1 H}(M) \cong \g_\C/\fq$. Finally, 
(A3) makes the proof of \cite[Thm.~1.6]{Ne12a} work, so that we can associate to every 
bounded unitary representation $(\rho,V)$ of $H$ 
a holomorphic Hilbert bundle $\bV := G \times_H V$ over the complex  $G$-manifold $M$ 
by defining $\beta \: \fq \to \gl(V)$ by $\beta(\fp^+) = \{0\}$ and $\beta\res_{\fh} = \dd\rho$.  
Now it is easy to check that all results in Sections~$1$ and $2$ of \cite{Ne12a} remain valid. 

\begin{defn}  \mlabel{def:d.1} 
We write $\Gamma(\bV)$ for the space of holomorphic sections 
of the holomorphic Hilbert bundle $\bV \to M = G/H$ on which the group $G$ acts by 
holomorphic bundle automorphisms. 
A unitary representation $(\pi, \cH)$ of $G$ is said to be 
{\it holomorphically induced from $(\rho,V)$} 
if there exists a $G$-equivariant linear injection 
$\Psi \: \cH \to \Gamma(\bV)$ such that the 
adjoint of the evaluation map $\ev_{\1 H} \: \cH \to V = \V_{\1 H}$ 
defines an isometric embedding $\ev_{\1 H}^* \: V \into \cH$. 
If a unitary representation $(\pi, \cH)$ holomorphically induced 
from $(\rho,V)$ exists, then it is uniquely determined 
(\cite[Def.~2.10]{Ne12a}) and we call $(\rho,V)$ {\it (holomorphically)  
inducible}. 

This concept of inducibility involves a choice of sign. 
Replacing $\fp^+$ by $\fp^-$ changes the complex structure on $G/H$ 
and thus leads to a different class of holomorphically inducible 
representations of of $H$. 
\end{defn}

\begin{thm} \mlabel{thm:c.1}
If the unitary representation $(\pi, \cH)$ of $G$ is holomorphically 
induced from the bounded $H$-representation $(\rho,V)$, then the 
following assertions hold: 
\begin{description}
\item[\rm(i)] $V \subeq \cH^\omega$ consists of analytic vectors. 
\item[\rm(ii)] $R \: \pi(G)' \to \rho(H)', A \mapsto A\res_V$ is an isomorphism of 
von Neumann algebras. 
\end{description}
\end{thm}

\begin{prf} (i) follows from \cite[Lemma~2.5]{Ne12a} and (ii) from \cite[Thm.~2.12]{Ne12a}. 
\end{prf}

\begin{thm} {\rm(\cite[Thm.~2.17]{Ne12a})} \mlabel{thm:c.3}
Suppose that $(\pi, \cH)$ is a unitary representation of $G$ and 
$V \subeq \cH$ is an $H$-invariant closed subspace such that 
\begin{description}
\item[\rm(HI1)] The representation $(\rho,V)$ of $H$ on $V$ is bounded. 
\item[\rm(HI2)] $V \cap (\cH^\infty)^{\fp^-}$ is dense in $V$. 
\item[\rm(HI3)] $\pi(G)V$ spans a dense subspace of $\cH$. 
\end{description}
Then $(\pi, \cH)$ is holomorphically induced from $(\rho,V)$. 
\end{thm}

\begin{exs} \mlabel{ex:c.4} 
(a) Let $G$ be  a simply connected Banach--Lie group for which $\g_\C$ also is the 
Lie algebra of a Banach--Lie group and $M = G/H$ is a Banach homogeneous space. 
If the subalgebras $\fp^\pm \subeq \g_\C$ satisfy the splitting condition \eqref{eq:splitcond}, 
then (A1-3) follows directly from the Inverse Function Theorem. This is the context of \cite{Ne12a}. 

(b) Let $G_B$ be a Banach--Lie group with Lie algebra $\g_B$, $H_B \subeq G_B$ and 
$M_B = G_B/H_B$ etc.\ as in (a). We assume that the splitting condition 
\eqref{eq:splitcond} is satisfied. 
In addition, let $\alpha \: \R \to \Aut(G_B)$ be a one-parameter group of automorphisms defining a continuous 
$\R$-action on $G_B$ and assume that the subalgebras 
$\fp_B^\pm$, $\fq_B$ and $\fh$ are $\alpha$-invariant. Then the subgroup 
\[ G := \{ g \in G_B \: \R \to G_B, t \mapsto \alpha_t(g) \ \mbox{ is smooth} \} \] 
of $G_B$ carries the structure of a Fr\'echet--BCH--Lie group with 
Lie algebra 
\[ \g := \{ x \in \g_B \: \R \to \g_B, t \mapsto \L(\alpha_t)x \ \mbox{ is smooth} \}, \] 
the Fr\'echet space of smooth vectors for the continuous $\R$-action on the Banach--Lie algebra $\g_B$. 
Likewise $H := G \cap H_B$ is a  Lie subgroup of $G$ for which 
$M := G/H$ is a smooth manifold consisting of the elements of $M_B = G_B/H_B$ with smooth 
orbit maps with respect to the one-parameter group of diffeomorphisms induced by 
$\alpha$ via $\alpha_t(gH_B) = \alpha_t(g)H_B$. 

Since the automorphisms $\L(\alpha_t)$ of $\g$ resp., $\g_\C$ are compatible with the 
BCH multiplication, it is easy to see with Lemma~\ref{lem:smooth} below 
that conditions (A1-3) are inherited by the closed subalgebras 
\[ \fh = \fh_B \cap \g, \quad \fp^\pm = \fp^\pm_B \cap \g_\C \quad \mbox{ and } \quad 
\fq = \fq_B \cap \g_\C.\] 
\end{exs}

\begin{lem} \mlabel{lem:smooth} Let $V_1$ and $V_2$ be Banach spaces and 
$(\alpha^1_t)$, resp., $(\alpha^2_t)$ define continuous $\R$-actions on $V_1$, resp., $V_2$. 
If $U \subeq V_1$ is an open invariant subset and $F \: U \to V_2$ an equivariant smooth map, 
then the induced map 
\[ F^\infty \: U^\infty := U \cap V_1^\infty  \to V_2^\infty, \quad 
v \mapsto F(v) \] 
is a smooth map on the open subset $U^\infty$ of the Fr\'echet space $V_1^\infty$. 
\end{lem}

\begin{prf} Let $D_j := \alpha_j'(0)$ denote the infinitesimal generator of $\alpha_j$.
Then we have to verify that all maps
$F_k \: U^\infty \to V_2, x \mapsto D_2^k F(x)$ 
are smooth. Since $\alpha^1$ defines a smooth $\R$-action on $U^\infty$, the map 
\[ \Phi \: \R \times U^\infty \to V_2, \quad (t,x) \mapsto 
F(\alpha_t^1(x))  = \alpha_t^2F(x) \] 
is smooth. Hence the map 
$F_k(x) = \frac{\partial^k}{\partial t^k}|_{t = 0} \Phi(t,x)$ is also smooth.
\end{prf}

From (A1-3) we derive the existence of open convex symmetric $0$-neighborhoods 
$U_\pm \subeq \fp^\pm$ and $U_0 \subeq \fh_\C$ 
for which the BCH-multiplication map
\[ U_+ \times U_0 \times U_- \to \g_\C, \quad 
(x_+,x_0,x_-) \mapsto x_+ * x_0 * x_- \] 
is biholomorphic onto an open $0$-neighborhood $U$ of $\g_\C$. 
For a bounded representation $(\rho,V)$ of $H_0$ we then define a holomorphic 
map 
\[ F_\rho \: U \to B(V), \quad 
F_\rho(x_+*x_0*x_-) := e^{\dd\rho(x_0)}.\] 

For the Banach case the equivalence of (i) and (ii) in the following theorem 
can also be found in \cite[Thm.~B]{Ne11c}. Its proof also 
works without change in our context. We include it for the sake of completeness. 

\begin{thm} \mlabel{thm:c.2} 
For a bounded representation $(\rho,V)$ of $H$, the following are equivalent: 
  \begin{description}
  \item[\rm(i)] $(\rho, V)$ is holomorphically inducible. 
  \item[\rm(ii)] $f_\rho(\exp x) 
:= F_\rho(x)$ defines a positive definite analytic 
function on 
a $\1$-neighborhood of $G$. 
  \item[\rm(iii)] The corresponding linear map 
$\beta \: U(\g_\C) \to B(V), \beta(D) = (L_D f_\rho)(\1)$ 
is positive definite. It is characterized by the property that
$\fp^+ U(\g_\C) + U(\g_\C) \fp^- \subeq \ker \beta$ and 
$\beta\res_{U(\fh_\C)} = \dd\rho$. 
  \end{description}
\end{thm}

\begin{prf} (i) $\Rarrow$ (ii): Let $(\pi, \cH)$ be the unitary 
representation of $G$ obtained by holomorphic induction from 
$(\rho, V)$. We identify $V$ with the corresponding closed subspace 
of $\cH$ and write $p_V\: \cH \to V$ for the corresponding orthogonal 
projection. For $v \in V \subeq (\cH^\omega)^{\fp^-}$ 
(Theorem~\ref{thm:c.1}), we let 
$f_\rho^v \: U_v \to G$ be a holomorphic 
map on an open convex $0$-neighborhood 
$U_v \subeq U$ satisfying $f_\rho^v(x) = \pi(\exp x)v$ for $x \in U_v \cap \g$. 
Then $\dd\pi(\fp^-)v = \{0\}$ implies that 
$L_z f_\rho^v = 0$ for $z \in \fp^-$. For $w \in V$ and $z \in \fp^+$, 
we also obtain 
\[ \la (R_z f_\rho^v)(x), w \ra 
= \la \dd\pi(z) f_\rho^v(x), w \ra
= \la f_\rho^v(x), \dd\pi(z^*) w \ra=0.\] 
This proves that $R_z (p_V \circ f_v) = 0$. 
We conclude that, for $x_\pm$ and $x_0$ sufficiently close to $0$, we  have 
\[ p_V f_\rho^v(x_+ * x_0 * x_-) = f_\rho^v(x_0) = e^{\dd\rho(x_0)}v = F_\rho(x_+ * x_0 * x_-)v.\]
Therefore $p_V \circ f_\rho^v$ extends holomorphically to $U$ and 
\[ \la \pi(\exp x)v,w \ra  = \la F_\rho(x)v,w \ra \quad \mbox{ for } \quad 
x \in U_v \cap \g, v,w \in V.\] 
We conclude that $F_\rho(x) = p_V \pi(\exp x) p_V$ 
holds for $x$ sufficiently close to $0$, and hence that 
$f_\rho(\exp x) = p_V \pi(\exp x) p_V$ defines a positive definite 
function on a $\1$-neighborhood of~$G$. 

(ii) $\Rarrow$ (i): From Theorem~\ref{thm:extension} it 
follows that some restriction of $f_\rho$ to a possibly smaller 
$\1$-neighborhood in $G$ extends to a global analytic positive 
definite function $\phi$. Then the vector-valued GNS construction 
yields a unitary representation of $G$ on the corresponding reproducing 
kernel Hilbert space $\cH_\phi \subeq V^G$ 
for which all the elements of $\cH_\phi^0 = \Spann(\phi(G)V)$ 
are analytic vectors. In particular, 
$V \subeq \cH_\phi^\omega$ consists of smooth vectors, 
and the definition of $f_\rho$ implies that 
$\dd\pi(\fp^-)V = \{0\}$. Therefore Theorem~\ref{thm:c.3} implies that the representation 
$(\pi, \cH_\phi)$ is holomorphically induced from 
$(\rho,V)$. 

(ii) $\Leftrightarrow$ (iii) follows from Theorem~\ref{thm:charaposdef}. 
The relation $U(\g_\C) \fp^- \subeq \ker \beta$ 
follows from the definition of $f_\rho$ which does not depend on the 
$x_-$-component. In view of $f_\rho(g^{-1}) = f_\rho(g)^*$, 
we have $\beta(D^*) = \beta(D)^*$ for $D \in U(\g_\C)$, 
and we thus also obtain $\fp^+ U(\g_\C) \subeq \ker \beta$, so that 
$\beta$ is determined by its restriction to $U(\fh_\C)$, where 
it coincides with $\dd\rho$. 
\end{prf}

\section{Finite order automorphisms of Hilbert--Lie algebras} \mlabel{app:d}

In this appendix we generalize some of the results on finite order automorphisms 
of complex, resp., compact 
semisimple Lie algebras (\cite[Sec.~X.5]{Hel78}) to Hilbert--Lie algebras.

Let $\fk$ be a Hilbert--Lie algebra and 
$\phi \in \Aut(\fk)$ be an automorphism of order $N$.

\begin{lem} \mlabel{lem:d.1} If $\fk$ is semisimple and non-zero, then $\fk^\phi \not=\{0\}$. 
\end{lem}

\begin{prf} We also write $\phi$ for the complex linear extension of 
$\phi$ to $\fk_\C$ and write 
\[ \fk_\C^n := \{ x \in \fk_\C \: \phi^{-1}(x) = e^{2\pi in/N} x \}. \] 
Assume that 
$\fk^\phi = \{0\}$. This means that $\fk_\C^0 = (\fk^\phi)_\C = \{0\}$. 
We show by induction that $\fk_\C^k = \{0\}$ for $k = 1,\ldots, N-1$. 
Assume $1 \leq k < N$ and that $\fk_\C^j = \{0\}$ holds for $j = 0, 1,\ldots, k-1$. 
Pick $x \in \fk_\C^k$. For each $j \in \Z$, there exists an $r \in \N$ such that 
$j + k r$ is congruent to one of the numbers $0,\ldots, k-1$ modulo $N$. We then 
obtain 
\[ (\ad x)^r \fk_\C^j \subeq \fk_\C^{j + rk} = \{0\},\] 
and conclude that $\ad x$ is nilpotent. Then $\ad x^*$ is also nilpotent. 
Moreover, $[x,x^*] \in [\fk_\C^k, \fk_\C^{-k}] \subeq \fk_\C^0 = \{0\}$ implies that 
$\ad x$ and $\ad x^* = (\ad x)^*$ 
commute. Thus $\ad x$ is a normal operator on the complex Hilbert space $\fk_\C$, 
and since it is nilpotent, we obtain $\ad x = 0$. Now $x \in \fz(\fk_\C) = \{0\}$ completes 
our inductive proof of $\fk_\C^k = \{0\}$ for $k = 0,\ldots, N-1$. This contradicts 
the assumption that $\fk$ is non-zero. 
\end{prf}

\begin{lem} \mlabel{lem:d.2} If $\ft \subeq \fk^\phi$ is maximal abelian, then 
$\ft_\fk := \fz_\fk(\ft)$ is maximal abelian in $\fk$. 
\end{lem}

\begin{prf} Clearly, $\ft_\fk$ is a closed subalgebra of $\fk$ invariant under $\phi$, hence a 
Hilbert--Lie algebra, endowed with a finite order automorphism $\phi\res_{\ft_\fk}$. 
Let $\fs = \fz(\ft_\fk)^\bot \cap \ft_\fk$ denote the commutator algebra of $\ft_\fk$. Then 
$\fs$ is also $\phi$-invariant and semisimple. If $\fs$ is non-zero, 
then Lemma~\ref{lem:d.1} implies that $\fs^\phi$ is non-zero, but 
this leads to the contradiction $\fs^\phi \subeq \fz_{\fk^\phi}(\ft) = \ft$. 
\end{prf}

Lemmas~\ref{lem:d.1} and \ref{lem:d.2} imply in particular, that there exists a maximal 
abelian subalgebra of $\fk$ which is $\phi$-invariant. 
According to \cite{Sch61}, $\fk_\C$ decomposes into an orthogonal sum of 
$\ft_\fk$-root spaces,
and this implies that $\fk_\C$ decomposes into $\ft$-weight spaces $\fk_\C^\alpha$, 
$\alpha \in \ft'$. 

Let $\Delta := \Delta(\fk,\ft) := \{ \alpha \in \ft'\: \fk_\C^\alpha \not= \{0\}\}$ 
denote the $\ft$-weight set of $\fk$. As $\ad \ft$ and $\phi$ commute, the weight spaces 
$\fk_\C^\alpha$ are $\phi$-invariant, so that we obtain a simultaneous diagonalization of 
$\ft$ and $\phi$ by the spaces 
\[ \fk_\C^{(\alpha,n)} := \fk_\C^\alpha \cap \fk_\C^n, \quad n \in \Z, \alpha \in \Delta.\] 

For $x,y \in \fk_\C^{(\alpha,n)}$ we then have 
$[x,y^*] \in \fk_\C^{(0,0)} = \ft_\C$, and for $h \in \ft_\C$ 
\[ \la h,[x,y^*]\ra = \la [h,y],x \ra
= \alpha(h) \la y,x \ra
= \la h, \la x, y \ra \alpha^\sharp\ra,\] 
where $\alpha^\sharp \in \ft_\C$ is the unique element satisfying 
$\la h,\alpha^\sharp\ra = \alpha(h)$ for $h \in \ft_\C$. This leads to 
\begin{equation}
  \label{eq:brarelx}
[x,y^*] = \la x,y\ra \alpha^\sharp \quad \mbox{ for } \quad x,y \in \fk_\C^{(\alpha,n)}.
\end{equation}
For $\|x\| = 1$ we obtain in particular 
$[x,x^*] = \alpha^\sharp$ and thus 
\[ \alpha([x,x^*]) = \alpha(\alpha^\sharp)  = \|\alpha^\sharp\|^2 > 0 
\quad \mbox{ for } \quad 0 \not= \alpha.\] 
We conclude that 
\[ \fk(\alpha,n) := \Spann_\R \{ x - x^*, i(x + x^*), i[x,x^*] \} \cong \su_2(\C) \] 
(cf.\ Lemma~\ref{lem:e.1}). For $y \bot x$ in $\fk_\C^{(\alpha,n)}$, we obtain 
$[x,y^*] = 0$ by \eqref{eq:brarelx} and thus 
\[ 0 \leq \la [x,y], [x,y] \ra 
= \la [x^*,[x,y]], y\ra = \la [[x^*,x],y], y \ra 
= -\alpha(\alpha^\sharp) \|y\|^2 \leq 0,\] 
so that $y = 0$, which means that 
\begin{equation}
  \label{eq:dim1}
\dim \fk_\C^{(\alpha,n)} = 1.
\end{equation}

\begin{lem}
  \mlabel{lem:d.3} 
If $\fk$ is simple, then the weight set $\Delta^\times := \Delta \setminus \{ 0\}$ 
does not decompose into two mutually orthogonal non-empty subsets. 
\end{lem}

\begin{prf} Suppose that $\Delta = \Delta_1 \dot\cup \Delta_2$ is a decomposition into 
mutually orthogonal subsets. Then, for $\alpha \in \Delta_1$, $\beta \in \Delta_2$, 
we have $\alpha + \beta \not\in \Delta$, so that $[\fk_\C^\alpha, \fk_\C^\beta] = \{0\}$. 
Therefore the subalgebra $\fk_1$ generated by the weight spaces $\fk_\C^\alpha$, 
$\alpha \in \Delta_1$, is invariant under brackets with all root spaces and with 
$\ft_\fk$, hence an ideal. As $\fk$ is simple and $\Delta_1 \not=\eset$, it follows that 
$\fk = \fk_1$, and this leads to $\Delta_2 = \eset$. 
\end{prf}


\begin{thebibliography}{aaaaaaa} 

\bibitem[AH78]{AH78} Albeverio, S., and R. J. H\o{}egh-Krohn, 
{\it The energy representation of Sobolev--Lie groups}, 
Composition Math. {\bf 36:1} (1978), 37--51 

\bibitem[Alb93]{Alb93} Albeverio, S., R. J. H\o{}egh-Krohn, J. A. Marion,  
D. H. Testard, and B. S. Torresani, ``Noncommutative Distributions
-- Unitary representations of Gauge Groups and Algebras,'' Pure and
Applied Mathematics {\bf 175}, Marcel Dekker, New York, 1993 

\bibitem[AP83]{AP83} Atiyah, M. F., and A. N. Pressley, 
{\it Convexity and loop groups}, in ``Arithmetic and Geometry: 
Papers dedicated to I. R. Shafarevich on the occasion of his 
60th birthday, Vol. II: Geometry,'' Birkh\"auser, Basel, 1983 

\bibitem[Bak07]{Bak07} Bakalov, B., N. M. Nikolov, 
K.-H. Rehren, I. Todorov, 
{\it Unitary positive energy representations 
of scalar bilocal quantum fields}, Comm. Math. Phys. 
{\bf 271:1} (2007), 223--246 

\bibitem[Ba69]{Ba69} Balachandran, V.\ K., {\it Simple $L^*$-algebras of classical
type}, Math.\ Ann.\ {\bf 180} (1969), 205--219 

\bibitem[Be79]{Be79} Berezin, F. A., {\it 
Representations of the continuous direct product of universal coverings 
of the group of motions of the complex ball}, 
Trans. Moscow Math. Soc. {\bf 2} (1979), 281--289 

\bibitem[BS71]{BS71} Bochnak, J., Siciak, J., 
{\it Analytic functions in topological vector spaces}, 
Studia Math. \textbf{39}  (1971), 77--112


\bibitem[Ca83]{Ca83} Carey, A. L., {\it Infinite Dimensional Groups and 
Quantum Field Theory}, Act. Appl. Math. {\bf 1} (1983), 321--333 

\bibitem[CR87]{CR87} Carey, A., and S. 
N. M. Ruijsenaars, {\it On fermion gauge
groups, current algebras, and Kac--Moody algebras}, Acta
Appl. Math. {\bf 10} (1987), 1--86 

\bibitem[CP86]{CP86} Chari, V., and A. Pressley, 
{\it New unitary representations 
of loop groups}, Math. Ann. {\bf 275} (1986), 87--104 

\bibitem[CP87]{CP87} ---, {\it Unitary representations of the maps $S^1 \to su(N,1)$}, 
Math. Proc. Camb. Phil. Soc. {\bf 102}(1987), 259--272 

\bibitem[CGM90]{CGM90} Cuenca Mira, J.\ A., A.\ Garcia Martin, and C.\ Martin
Gonzalez, {\it Structure theory of $L^*$-algebras}, Math.\ Proc.\
Camb.\ Phil.\ Soc.\ {\bf 107} (1990), 361--365


\bibitem[FH05]{FH05} Fewster, Chr., and S. Hollands, {\it 
Quantum energy inequalities in two-dimensional conformal field theory}, 
Rev. Math. Phys.  {\bf 17:5} (2005), 577--612 


\bibitem[GN12]{GN12} Gl\"ockner, H., and K.-H. Neeb, ``Infinite dimensional 
Lie groups, Vol. I, Basic Theory and Main Examples,'' book in preparation 

\bibitem[Go03]{Go03} Goertsches, O., ``Variationally complete and 
hyperpolar actions on compact symmetric spaces,'' Diploma thesis, K\"oln, 2003

\bibitem[GW84]{GW84} Goodman, R., and N. R. Wallach, {\it Structure and unitary
cocycle representations of loop groups and the group of
diffeomorphisms of the circle}, J. reine ang. Math. {\bf 347} (1984), 69--133

\bibitem[Hel78]{Hel78} Helgason, S., ``Differential Geometry, Lie Groups, and Symmetric 
Spa\-ces,'' Acad. Press, London, 1978 


\bibitem[HoMo98]{HoMo98} Hofmann, K.\ H., and S.\ A.\ Morris, ``The Structure of
Compact Groups,'' Studies in Math., de Gruyter, Berlin, 1998 

\bibitem[HN12]{HN12} ---, {\it On convex hulls of orbits of Coxeter groups 
and Weyl groups}, Preprint, 2012 

\bibitem[JK85]{JK85} Jakobsen, H.~P., and V.\ Kac, 
{\it A new class of unitarizable
highest weight representations of infinite-di\-men\-sio\-nal Lie algebras}, 
in ``Non--linear equations in classical and quantum field theory,''
N.\ Sanchez ed., Springer Verlag, Lecture Notes in Physics {\bf 226} (1985), 1--20

\bibitem[JK89]{JK89} ---, {\it A new class of unitarizable
highest weight representationsof infinite-di\-men\-sio\-nal Lie algebras,
II}, J.\ Funct.\ Anal. {\bf 82} (1989), 69--90 

\bibitem[JW10]{JW10} Janssens, B., and C. Wockel, {\it Universal Central Extensions 
of Gauge Algebras and Groups}, arXiv:math.DG:1010.3569 

\bibitem[Ka90]{Ka90} Kac, V., ``Infinite-dimensional Lie Algebras," Cambridge University
Press, $3^{rd}$ printing, 1990  

\bibitem[KP84]{KP84} Kac, V. G., and D. H. Peterson, {\it Unitary structure in 
representations of infinite dimensional groups and a convexity theorem}, 
Invent. Math. {\bf 76} (1984), 1--14 


\bibitem[Kue06]{Kue06} K\"uhn, K., {\it Direct limits of diagonal chains 
of type O, U and Sp, and their homotopy groups}, Comm. Alg. 
{\bf 34} (2006), 75--87 

\bibitem[LN04]{LN04} Loos, O., and E. Neher, ``Locally finite root systems,'' 
Memoirs of the Amer. Math. Soc., Vol. 171, {\bf 811}, 2004 

\bibitem[Mai02]{Mai02} Maier, P., {\it Central extensions of topological
current algebras}, in 
``Geometry and Analysis on Finite-
and Infinite-Dimensional Lie Groups,'' Eds. A.~Strasburger et al, 
Banach Center Publications {\bf 55}, Warszawa, 2002; 61--76 

\bibitem[MN03]{MN03} Maier, P., and K.-H. Neeb, 
{\it Central extensions of current groups},
Math. Annalen {\bf 326:2} (2003), 367--415 

\bibitem[Mais62]{Mais62} Maissen, B., {\it Lie-Gruppen mit Banachr\"aumen als
Parameterr\"aume}, Acta Math. {\bf 108} (1962), 229--269 

\bibitem[MR85]{MR85} A. Medina and P. Revoy, {\it Alg\`ebres de Lie et produit scalaire 
invariant}, Ann. scient. \'Ec. Norm. Sup. $4^e$ s\'erie {\bf 18}(1985), 
533-561


\bibitem[Mi89]{Mi89} Mickelsson, J., ``Current Algebras and Groups,'' Plenum Press,
New York, 1989 

\bibitem[MY06]{MY06} Morita, Y., and Y.\ Yoshii, {\it 
Locally extended affine Lie algebras}, J. Algebra {\bf 301} (2006), 59--81 

\bibitem[Ne98]{Ne98} Neeb, K.-H., {\it Holomorphic highest weight representations
of infinite dimensional complex classical groups}, 
J.\ reine angew. Math.\ {\bf 497} (1998), 171--222  

\bibitem[Ne00]{Ne00} ---, ``Holomorphy and Convexity in Lie Theory,'' 
Expositions in Mathematics {\bf 28}, de Gruyter Verlag, Berlin, 2000  


\bibitem[Ne02a]{Ne02a} ---, 
{\it Classical Hilbert--Lie groups, their extensions and their
homotopy groups}, in ``Geometry and Analysis on Finite-
and Infinite-Dimensional Lie Groups,'' 
Eds. A.~Strasburger et al., 
Banach Center Publications {\bf 55}, Warszawa 2002; 87--151 

\bibitem[Ne02b]{Ne02b} Neeb, K.-H., {\it Central extensions of 
infinite-dimensional Lie groups}, Ann. Inst. Fourier (Grenoble) {\bf 52}:5 
(2002), 1365--1442

\bibitem[Ne06]{Ne06} ---, {\it Towards a Lie theory of locally convex 
groups}, Jap. J. Math. 3rd ser. {\bf 1:2} (2006), 291--468 

\bibitem[Ne10a]{Ne10a} ---, {\it Unitary highest weight modules of 
locally affine Lie al
gebras}, in ``Quantum Affine Algebras, Extended Affine Lie Algebras
and their Applications'', Eds. Y. Gao et al, Contemporary Math. 
{\bf  506}, Amer. Math. Soc., 2010; 227--262

\bibitem[Ne10b]{Ne10b} ---, {\it Semibounded representations and invariant 
cones in infinite dimensional Lie algebras}, Confluentes Math. {\bf 2:1} 
(2010), 37--134

\bibitem[Ne11a]{Ne11a} ---, {\it On analytic vectors for 
unitary representations of infinite dimensional Lie groups}, 
Ann. Inst. Fourier  {\bf 61:5} (2011), 1441--1476 

\bibitem[Ne11b]{Ne11b} ---, {\it 
Projective semibounded representations of doubly extended 
Hilbert--Lie groups}, in preparation 

\bibitem[Ne11c]{Ne11c} ---, {\it Semibounded representations 
of hermitian Lie groups}, arXiv.math.RT:\\ 1104.2234v2, 20 May 2011 

\bibitem[Ne12a]{Ne12a} ---, {\it Holomorphic realization of unitary 
representations of Banach--Lie groups}, in Progress in Mathematics, 
A.\ Huckleberry et al eds., to appear; arXiv:math.RT.1011.1210v1, 4 Nov 2010

\bibitem[NS11]{NS11} Neeb,  K.-H., and H. Sepp\"anen, {\it 
Borel--Weil theory for groups over commutative Banach algebras}, 
J. reine angew. Math. {\bf 655} (2011), 165–187 

\bibitem[NeWo09]{NeWo09} Neeb, K.-H., and Chr.\ Wockel, {\it Central 
extensions of  groups of section}, 
Annals of Global Analysis and Geometry {\bf 36:4} (2009), 381--418


\bibitem[Neh93]{Neh93} Neher, E., {\it Generators and relations for $3$-graded Lie
algebras}, Journal of Algebra {\bf 155} (1993), 1--35 

\bibitem[Ot95]{Ot95} Ottesen, J. T., 
``Infinite Dimensional Groups and Algebras in Quantum 
Physics,'' Springer-Verlag, Lecture Notes in 
Physics {\bf m 27}, 1995 


\bibitem[Pa68]{Pa68} Palais, R., ``Foundations of Global Non-Linear Analysis,'' 
W. A. Benjamin, Inc., New York, Amsterdam, 1968 

\bibitem[PS86]{PS86} Pressley, A., and G. Segal, ``Loop Groups," Oxford University Press, 
Oxford, 1986



\bibitem[Sch60]{Sch60} Schue, J.\ R., {\it Hilbert space methods in the theory of Lie
algebras}, Transactions of the Amer.\ Math.\ Soc. {\bf 95} (1960),
69--80. 

\bibitem[Sch61]{Sch61} ---, {\it Cartan decompositions for $L^{*}$-algebras}, Trans.\
Amer.\ Math.\ Soc.\ {\bf 98} (1961), 334--349

\bibitem[Sch90]{Sch90} Schm\"udgen, K., 
``Unbounded Operator Algebras and
 Representation Theory,'' Operator Theory: Advances and Applications
 {\bf 37}, Birkh\"auser Verlag, Basel, 1990 

\bibitem[SeG81]{SeG81} Segal, G., {\it Unitary representations of some
infinite-dimensional groups}, Comm.\ Math.\ Phys. {\bf 80} (1981), 301--342 

\bibitem[Se58]{Se58} Segal, I.E., {\it Distributions in Hilbert spaces and canonical
systems of operators}, Trans.\ Amer.\ Math.\ Soc.\ 
{\bf 88} (1958), 12--41 

\bibitem[Se78]{Se78} ---, {\it The Complex-wave Representation of the Free
Boson Field}, in ``Topics in Funct. Anal.,'' Adv.\ in Math. Suppl. 
Studies {\bf 3} (1978), 321--343 

\bibitem[St99]{St99} Stumme, N., ``Locally Finite Split Lie Algebras,''
Ph.D.\ thesis, Darmstadt University of Technology, 1999


\bibitem[St01]{St01} ---, {\it Automorphisms and conjugacy 
of compact real forms of the classical infinite dimensional matrix 
Lie algebras}, Forum Math. {\bf 13:6}  (2001),  817--851


\bibitem[Te89]{Te89} Terng, C.-L., {\it Proper Fredholm submanifolds of 
Hilbert space}, J. Diff. Geom. {\bf 29} (1989), 9--47 

\bibitem[To87]{To87} Torresani, B. S., {\it Unitary positive energy 
representations of the gauge group}, 
Letters in Math. Physics {\bf 13} (1987), 7--15 

\bibitem[TL99]{TL99} Toledano Laredo, V., {\it Positive energy 
representations of the loop groups of non-simply connected Lie groups}, 
Comm. Math. Phys. {\bf 207:2} (1999),  307--339

\bibitem[VGG74]{VGG74} Vershik, A. M., Gelfand, I. M., and 
M. I. Graev, {\it Irreducible representations of the group 
$G^X$ and cohomology}, Funct. Anal. and its Appl. {\bf 8:2} (1974), 67--69

\bibitem[VGG80]{VGG80} ---, 
{\it Representations of the group of functions taking values in a compact Lie group}, 
Compositio Math. {\bf 42:2} (1980/81), 217--243 
.
\bibitem[Wl82]{Wl82} Wloka, J., {\it Partielle Differentialgleichungen}, 
B. G. Teubner Stuttgart, 1982 

\bibitem[Wo06]{Wo06} Wockel, Chr., {\it Smooth extensions and spaces of smooth 
and holomorphic mappings}, J. Geom. Symmetry Phys. {\bf 5} (2006), 118–-126 

\bibitem[YY10]{YY10} Y. Yoshii, {\it Locally extended affine root 
systems}, in this volume: ``Quantum affine algebras, extended affine Lie 
algebras and applications'', Eds. Y. Gao et. al., Contemp. Math.  
{\bf  506}, Amer. Math. Soc., 2010; 285--302 

\end{thebibliography}
\end{document}